\documentclass{article}
\usepackage{amsmath, amssymb}
\usepackage{amssymb,pb-diagram}
\usepackage{amscd}

\usepackage{fancybox}
\DeclareMathOperator\irr{Irr}
\DeclareMathOperator\ltop{\sigma_\text{top}}
\DeclareMathOperator\STop{\Sigma_\text{top}}

\newcommand{\ZZ}{\mathbb Z}
\newcommand{\PP}{\mathbb P}
\newcommand{\QQ}{\mathbb Q}
\newcommand{\CC}{\mathbb C}

\newcommand{\mcD}{\mathcal D}
\newcommand{\mcQ}{\mathcal Q}

\newcommand{\tor}{{\mathrm {tor}}}
\newcommand{\hor}{{\mathrm {hor}}}
\newcommand{\ver}{{\mathrm {ver}}}
\newcommand{\red}{{\mathrm {red}}}

\newcommand{\MW}{\mathop {\rm MW}\nolimits}
\newcommand{\Gal}{\mathop {\rm Gal}\nolimits}

\newcommand{\NS}{\mathop {\rm NS}\nolimits}

\newcommand{\Red}{\mathop {\rm Red}\nolimits}
\newcommand{\rank}{\mathop {\rm rank}\nolimits}
\newcommand{\Pic}{\mathop {\rm Pic}\nolimits}

\newcommand{\Supp}{\mathop {\rm Supp}\nolimits}
\newcommand{\Sing}{\mathop {\rm Sing}\nolimits}

\def\miya#1{\mathrel{\mathop{\rightarrow}\limits^#1}}

\newcommand{\mcC}{\mathcal C}

\newtheorem{thm}{Theorem}[section]
\newtheorem{cor}{Corollary}[section]
\newtheorem{prop}{Proposition}[section]
\newtheorem{lem}{Lemma}[section]
\newtheorem{defin}{Definition}[section]
\newtheorem{exmple}{Example}[section]
\newtheorem{rem}{Remark}[section]
\newtheorem{qz}{Question}[section]

\newcommand{\I}{\mathop {\rm I}\nolimits}
\newcommand{\II}{\mathop {\rm II}\nolimits}
\newcommand{\III}{\mathop {\rm III}\nolimits}
\newcommand{\IV}{\mathop {\rm IV}\nolimits}
\newcommand{\qed}{\hfill $\Box$}

\newcommand{\proof}{\noindent{\textsl {Proof}.}\hskip 3pt}
\newcommand{\proofend}{\qed \par\smallskip\noindent}

\renewcommand{\thesubparagraph}{\theparagraph.\@arabic\c@subparagraph}
\begin{document}
  
  \begin{center}
  
 {\bf  \Large 
 Elliptic dihedral covers in dimension $2$, \\ 
 geometry of sections of elliptic surfaces,   \\ and \\
 Zariski pairs for line-conic arrangements
 }
\bigskip

\bigskip
\large 
Hiro-o TOKUNAGA\footnote{Research partly supported by the research grant 22540052
from JSPS}

\end{center}
\normalsize



\bigskip

{\large \bf Introduction}

\bigskip

In this article, all varieties are defined over the field of complex numbers, ${\mathbb C}$.  
Let $X$ and $Y$ be normal projective varieties. We call $X$ a dihedral cover of $Y$ if 
there exists a finite surjective morphism $\pi : X \to Y$ such that
the induced field extension of the rational function fields $\CC(X)/\CC(Y)$ is a Galois
extension whose Galois group is isomorphic to a dihedral group. 

 Let $D_{2n}$ be the dihedral group of order $2n$. In order to present $D_{2n}$, we use
 the notation
 \[
 D_{2n} = \langle \sigma, \tau \mid \sigma^2 = \tau^n = (\sigma\tau)^2 = 1\rangle.
 \]
 By $D_{2n}$-covers, we mean  a Galois cover whose Galois group is isomorphic to $D_{2n}$. 
 Given a $D_{2n}$-cover, we obtain a double cover, $D(X/Y)$, canonically by considering the
 $\CC(X)^{\tau}$-normalization of $Y$, where $\CC(X)^{\tau}$ denotes the fixed field 
 of the subgroup generated by $\tau$.  We denote these covering morphisms by
 $\beta_1(\pi) : D(X/Y) \to Y$ and $\beta_2(\pi) : X \to D(X/Y)$, respectively. 
 
 In \cite{tokunaga04}, we introduce a notion of an elliptic $D_{2n}$-cover, whose definition
 is as follows:
 
 \begin{defin}\label{def:dihed-elliptic}{A $D_{2n}$-cover $\pi : X \to Y$ is called an
  elliptic $D_{2n}$-cover if it it satisfies the following condition:
  
  $\bullet$ $D(X/Y)$ has a structure of an elliptic fiber space $\varphi : D(X/Y) \to 
  S$ over a projective variety $S$ with a section $O : S \to D(X/Y)$.
  
  $\bullet$ The covering transformation $\sigma_{\beta_1(\pi)}$  coincides with
  the inversion with respect to the group law on the generic fiber $D(X/Y)_{\eta}$.
  Here the group law on $D(X/Y)_{\eta}$ is given by regarding $O$ as the zero element.
  }
  \end{defin}
  
  In this article, as a continuation of \cite{tokunaga04}, we study an elliptic $D_{2p}$-cover ($p$: odd prime) of a rational ruled
  surface $\Sigma_d$ ($d$: even).   
   Our main results are Theorems~\ref{thm:main} and~\ref{thm:main1}.   As an application, we study some
   Zariski pairs  of 
   degree $7$ for line-conic arrangements.  Let us recall the definition of a Zariski pair.
   
   \begin{defin}\label{def:zpair}{\rm 
 A pair  $(B_1, B_2)$ of reduced plane curves $B_i$ 
 $(i = 1, 2)$ of degree $n$ is called a 
 Zariski pair of degree $n$ if it satisfies the following condition:
   
   \begin{enumerate}
   \item[(i)] $B_i$ $(i = 1, 2)$ are curves of degree $n$. The combinatorial type (see Definition~\ref{def:combinatorics} below) of $B_1$ is the same as that of $B_2$
   
   
   \item[(ii)] $(\PP^2, B_1)$ is not homeomorphic to $(\PP^2, B_2)$.
   \end{enumerate}
   }
   \end{defin}
  
  \begin{defin}\label{def:combinatorics}
  {\rm (\cite{act}) The \emph{combinatorial type} of a curve $B$ is given by a 7-tuple
$$(\irr(B), \deg, \Sing(B), \STop(B), \ltop, \{B(P)\}_{P\in \Sing(B)},
\{\beta_P\}_{P\in \Sing(B)}),$$ where:
\begin{itemize}
\item $\irr(B)$ is the set of irreducible components of $B$ and 
$\deg:\irr(B)\to\ZZ$ assigns to each irreducible component its degree.
\item $\Sing(B)$ is the set of singular points of $B$, 
$\STop(B)$ is the set of topological types of $\Sing(B)$, and
$\ltop:\Sing(B)\to\STop(B)$ assigns to each singular point its topological type.
\item $B(P)$ is the set of local branches of $B$ at 
$P\in \Sing(B)$, (a local branch can be seen as an arrow in the dual graph 
of the minimal resolution of $B$ at $P$, see~\cite[Chapter~II.8]{Eisenbud-Neumann} for
details) and $\beta_P : B(P) \to \irr(B)$ assigns to each local branch the global
irreducible component containing it.
\end{itemize}
Two curves $B_1$ and $B_2$ is said to have the 
\emph{same combinatorial type}
(or simply the \emph{same combinatorics}) if their data of combinatorial types
$$(\irr(B_i), \deg_i, \Sing(B_i), \STop(B_i), \ltop_i,
\{\beta_{i,P}\}_{P\in \Sing(B_i)}, \{B_i(P)\}_{P\in \Sing(B_i)}), \quad i=1,2,$$
are equivalent, that is, if $\STop(B_1)=\STop(B_2)$, 
and there exist bijections $\varphi_{\Sing} :\Sing(B_1)\to \Sing(B_2)$, 
$\varphi_P : B_1(P)\to B_2(\varphi_{\Sing}(P))$ (restriction of a bijection
of dual graphs) for each $P\in \Sing(B_1)$, and 
$\varphi_{\irr}:\irr(B_1)\to \irr(B_2)$ such 
that $\deg_2 \circ \varphi_{\irr}=\deg_1$, $\ltop_2 \circ \varphi_{\Sing}=\ltop_1$, and
$\beta_{2,\varphi_{\Sing}(P)} \circ \varphi_P= \varphi_{\irr} \circ \beta_{1,P}$.
}
\end{defin}

Note that when  $B_i$ $(i = 1, 2)$ are irreducible, $B_1$ and $B_2$ have the same 
combinatorics
if they have the same degree and the same local topological types for singularities. 
On the other extreme, for line arrangements,  $B_1$ and $B_2$ have the same
combinatorial type if they have the same set 
of incidence relations. 
The first example of a Zariski pair is given by Zariski (\cite{zariski29, zariski37}), which is 
as follows:

\begin{exmple}\label{eg:zariski}{\rm 
Let $(B_1, B_2)$ be a pair of irreducible sextics such that (i) both of $B_1$ and $B_2$
have six cusps as their singularities, and (ii) the six cusp of $B_1$ are on a conic, while
no such conic for $B_2$. Then $(B_1, B_2)$ is a Zariski pair.
}
\end{exmple}

 For these twenty years, Zariski pairs have been studied by many mathematicians and many 
 examples have been found (see \cite{act} and its reference).  Among them, Zariski pairs for
 line arrangements of degrees $9$ and $11$ are considered by Artal Bartolo, Carmonoa Ruber, Cogolludo
 Agustin and Marco
 Buzunariz (\cite{Artal-Carmona-ji-Marco-real, Artal-Carmona-ji-Marco-ryb}), Rybnikov
 (\cite{rybnikov}) and
 those for conic arrangements of degree $8$ are considered by Namba and 
 Tsuchihashi (\cite{namba-tsuchi}).  In this article, we study Zariski pairs for
 line-conic arrangement.

As we explain in \cite{act}, the study of Zariski pairs, in general, consists of two parts:

\begin{enumerate}
 \item[(I)]  To give curves $B_1$ and $B_2$ having the same combinatorics, but some ^^ ^^ {\it different property}," e.g., the location of singularities as in Example~\ref{eg:zariski}.
 
 \item[(II)] To show $(\PP^2, B_1)$ is not homeomorphic to $(\PP^2, B_2)$.
 
 \end{enumerate}
 One of our goals in this article is to add another method to find two curves with the same combinatorics. Namely we make use of 
 geometry of sections of the Mordell-Weil group of an elliptic surface as follows:
 
  Let $\varphi : S \to \PP^1$ be an elliptic surface with a section $O$ and let 
 \[
\begin{CD}
S' @<{\mu}<< S \\
@V{f'}VV                 @VV{f}V \\
\Sigma_d @<<{q}< \widehat{\Sigma}_d.
\end{CD}
\]
 be the double cover diagram for $S$ (see \ref{subsec:double-cover}).  Let $\Delta_1$ and $\Delta_2$ be sections of
 $\Sigma_d$ with $\Delta_i^2 = d$  $(i = 1, 2)$. Suppose that $(q\circ f)^*(\Delta_i)$ consists of two sections $s_{\Delta_i}^{\pm}$
 for each $i$ and $\widehat{\Sigma}_d$ can be blow down to $\PP^2$, which we denote by $\overline{q} :
 \widehat{\Sigma}_d \to \PP^2$. Let $[2]s_{\Delta_i}^+$ be the duplication of $s_{\Delta_i}^+$ in $\MW(S)$. 
 In order to produce two reduced curves $B_1$ and $B_2$ with
 the same combinatorics, we use $\overline{q}\circ f(s_{\Delta_i}^+)$ $( i = 1, 2)$, 
 $\overline{q}\circ f([2]s_{\Delta_i})$, 
 and $\overline{q}(\Delta(S/\widehat{\Sigma}_d))$, where $\Delta(S/\widehat{\Sigma}_d)$ is the branch locus of $f$. The author hopes that this method add  a new viewpoint to the study of elliptic surfaces and their Mordell-Weil groups.

 As for (II), we also make use of theory of dihedral covers and elementary arithmetic on the Mordell-Weil group of an
 elliptic surface as in our previous papers (\cite{tokunaga98, tokunaga99, tokunaga04}).

Now let us explain 
line-conic arrangements of 
degree $7$ considered in this article.

\medskip

\noindent{\bf Line-conic arrangement 1}

 Let $C_i$ $(i =1, 2)$  be smooth conics and let
$L_j$ $(i = 1, 2, 3, 4)$ be lines as follows:

\begin{enumerate}

\item[(i)] Both $L_1$ and $L_2$ meet $C_1$ transversely. We put
$C_1\cap L_1 = \{P_1, P_2\}, \, C_1\cap L_2 = \{P_3, P_4\}$.

\item[(ii)] $C_2$ is tangent to $C_1$ at two distinct points $\{Q_1, Q_2\}$ or
at one point $\{Q\}$.  We call the former type $(a)$ and the latter type $(b)$.

\item [(iii)]The tangent lines at $C_1\cap C_2$ do not pass through $L_1\cap L_2$.

\item[(iv)] $C_2$ is tangent to $L_1$ and $L_2$.

\item[(v)] $L_3$ passes through $P_1$ and $P_3$.

\item[(vi)] $L_4$ passes through $P_1$ and $P_4$.
\item[(vii)] Both $L_3$ and $L_4$ meet $C_2$ transversely.

\end{enumerate}

We put $B_1: = C_1 + C_2 + L_1 + L_2 + L_3$ and $B_2 := C_1 + C_2 + L_1+ L_2 + L_4$.
Then $B_1$ and $B_2$ have the same combinatorics.

\begin{center}
\unitlength 0.1in
\begin{picture}( 42.4300, 19.0800)(  4.0000,-24.0600)
%
\special{pn 13}%
\special{ar 1662 1580 558 416  0.0000000 6.2831853}%
%
\special{pn 13}%
\special{ar 1668 1574 558 832  0.8238768 6.2831853}%
\special{ar 1668 1574 558 832  0.0000000 0.7853982}%
%
\special{pn 13}%
\special{pa 420 1812}%
\special{pa 2148 730}%
\special{fp}%
%
\special{pn 13}%
\special{pa 400 1664}%
\special{pa 2218 2188}%
\special{fp}%
%
\special{pn 13}%
\special{pa 1910 602}%
\special{pa 2102 2374}%
\special{fp}%
%
\special{pn 13}%
\special{ar 4010 1580 556 416  0.0000000 6.2831853}%
%
\special{pn 13}%
\special{ar 4016 1574 558 832  0.8238768 6.2831853}%
\special{ar 4016 1574 558 832  0.0000000 0.7853982}%
%
\special{pn 13}%
\special{pa 2768 1812}%
\special{pa 4496 730}%
\special{fp}%
%
\special{pn 13}%
\special{pa 2750 1664}%
\special{pa 4566 2188}%
\special{fp}%
%
\special{pn 13}%
\special{pa 4458 652}%
\special{pa 3222 2188}%
\special{fp}%
\put(11.3000,-9.4700){\makebox(0,0)[lb]{$C_1$}}%
\put(15.9700,-13.1200){\makebox(0,0)[lb]{$C_2$}}%
\put(21.9200,-6.7800){\makebox(0,0)[lb]{$L_1$}}%
\put(22.7500,-22.8400){\makebox(0,0)[lb]{$L_2$}}%
\put(21.5000,-24.7000){\makebox(0,0)[lb]{$L_3$}}%
\put(20.2600,-9.2800){\makebox(0,0)[lb]{$P_1$}}%
\put(8.7400,-13.3700){\makebox(0,0)[lb]{$P_2$}}%
\put(21.4100,-21.3100){\makebox(0,0)[lb]{$P_3$}}%
\put(9.7000,-19.9600){\makebox(0,0)[lb]{$P_4$}}%
\put(22.9400,-15.7400){\makebox(0,0)[lb]{$Q_1$}}%
\put(11.6200,-16.7000){\makebox(0,0)[lb]{$Q_2$}}%
\put(45.4700,-7.8000){\makebox(0,0)[lb]{$L_1$}}%
\put(43.9400,-9.4700){\makebox(0,0)[lb]{$P_1$}}%
\put(33.5700,-10.4300){\makebox(0,0)[lb]{$C_1$}}%
\put(32.2900,-13.7600){\makebox(0,0)[lb]{$P_2$}}%
\put(41.3100,-13.6300){\makebox(0,0)[lb]{$C_2$}}%
\put(32.8600,-16.7600){\makebox(0,0)[lb]{$Q_2$}}%
\put(46.4300,-16.6400){\makebox(0,0)[lb]{$Q_1$}}%
\put(32.2200,-19.7700){\makebox(0,0)[lb]{$P_4$}}%
\put(45.5400,-21.3700){\makebox(0,0)[lb]{$P_3$}}%
\put(44.3200,-23.4200){\makebox(0,0)[lb]{$L_2$}}%
\put(30.3000,-23.6800){\makebox(0,0)[lb]{$L_4$}}%
\put(4.9000,-9.3000){\makebox(0,0)[lb]{$B_1$}}%
\put(29.1000,-9.1000){\makebox(0,0)[lb]{$B_2$}}%
\end{picture}%

\medskip

Line-conic arrangement 1 of type $(a)$ 
\end{center}

\noindent{\bf Line-conic arrangement 2} 
\begin{center}
\unitlength 0.1in
\begin{picture}( 54.5000, 21.3400)(  5.6000,-22.8300)
%
\special{pn 13}%
\special{ar 1872 1206 512 512  0.0000000 6.2831853}%
%
\special{pn 13}%
\special{ar 1886 1212 1326 512  0.0000000 6.2831853}%
%
\special{pn 13}%
\special{ar 1866 1212 512 1064  0.0000000 6.2831853}%
%
\special{pn 13}%
\special{pa 2314 368}%
\special{pa 2326 2090}%
\special{fp}%
%
\special{pn 13}%
\special{ar 4672 1214 512 512  0.0000000 6.2831853}%
%
\special{pn 13}%
\special{ar 4686 1220 1326 512  0.0000000 6.2831853}%
%
\special{pn 13}%
\special{ar 4666 1220 512 1064  0.0000000 6.2831853}%
%
\special{pn 13}%
\special{pa 5392 478}%
\special{pa 3952 1982}%
\special{fp}%
\put(6.7200,-17.0200){\makebox(0,0)[lb]{$C_1$}}%
\put(16.4000,-4.1400){\makebox(0,0)[lb]{$C_2$}}%
\put(15.4400,-10.1400){\makebox(0,0)[lb]{$C_3$}}%
\put(23.2000,-22.8600){\makebox(0,0)[lb]{$L_1$}}%
\put(23.8400,-6.9400){\makebox(0,0)[lb]{$P_1$}}%
\put(11.7600,-6.7000){\makebox(0,0)[lb]{$P_2$}}%
\put(23.8400,-19.1000){\makebox(0,0)[lb]{$P_3$}}%
\put(12.4800,-18.6200){\makebox(0,0)[lb]{$P_4$}}%
\put(44.0000,-4.7800){\makebox(0,0)[lb]{$C_2$}}%
\put(53.6800,-7.6600){\makebox(0,0)[lb]{$P_1$}}%
\put(39.7600,-7.1800){\makebox(0,0)[lb]{$P_2$}}%
\put(43.1200,-10.3800){\makebox(0,0)[lb]{$C_3$}}%
\put(34.4800,-16.9400){\makebox(0,0)[lb]{$C_1$}}%
\put(37.8400,-18.4600){\makebox(0,0)[lb]{$P_4$}}%
\put(51.6800,-18.9400){\makebox(0,0)[lb]{$P_3$}}%
\put(36.7200,-21.7400){\makebox(0,0)[lb]{$L_2$}}%
\put(5.8000,-4.4000){\makebox(0,0)[lb]{$B_1$}}%
\put(33.9000,-4.3000){\makebox(0,0)[lb]{$B_2$}}%
\put(11.6000,-11.8000){\makebox(0,0)[lb]{$Q_1$}}%
\put(24.8000,-11.8000){\makebox(0,0)[lb]{$Q_2$}}%
\put(16.5000,-6.4000){\makebox(0,0)[lb]{$Q_3$}}%
\put(16.5000,-19.8000){\makebox(0,0)[lb]{$Q_4$}}%
\put(44.9000,-19.8000){\makebox(0,0)[lb]{$Q_4$}}%
\put(39.3000,-11.8000){\makebox(0,0)[lb]{$Q_1$}}%
\put(52.7000,-11.7000){\makebox(0,0)[lb]{$Q_2$}}%
\put(45.0000,-6.4000){\makebox(0,0)[lb]{$Q_3$}}%
\end{picture}%

\medskip

Line-conic arrangement 2 of type $(a)$ 
\end{center}

Let $C_1, C_2$ and $C_3$ be smooth conics and
$L_1$ and $L_2$ be lines as follows:

\begin{enumerate}

\item[(i)] $C_1$ and $C_2$ meet transversely. We put 
$C_1\cap C_2 = \{P_1, P_2, P_3, P_4\}$.

\item[(ii)] $C_3$ is tangent to both $C_1$ and $C_2$ such that the intersection multiplicities at 
intersection points are all
even. By exchanging $C_1$ and $C_2$ if necessary, we may assume that there
are three possibilities: 

\noindent $(a)$ $C_3 \cap C_1 = \{Q_1, Q_2\}, C_3\cap C_2 = \{Q_3, Q_4\}$, 

\noindent $(b)$
$C_3\cap C_1 = \{Q_1\}, C_3\cap C_2 = \{Q_2, Q_3\}$ or 

\noindent $(c)$ $C_3\cap C_1 = \{Q_1\}, C_3\cap C_2 = \{Q_2\}$. 

\item[(iii)] No tangent line at $Q_i$ is bitangent of $C_1 + C_2$.

\item[(iv)] $L_1$ passes through $P_1$ and $P_3$.

\item[(v)] $L_2$ passes through $P_1$ and $P_4$.
\item[(vi)] Both of $L_1$ and $L_2$ meet $C_3$ transversely.

\end{enumerate}

  We put $B_1: = C_1 + C_2 + C_3 + L_1$ $B_2 := C_1 + C_2 + C_3+ L_2$.
Then $B_1$ and $B_2$ have the same combinatorics.




\begin{thm}\label{thm:zpair}{ 
\begin{enumerate}
\item[(i)] Let $(B_1, B_2)$ be the pair of Line-conic arrangement 1. 
Then $(B_1, B_2)$ is a Zariski pair.

\item[(ii)] Let $C_1$ and $C_2$ be conics intersecting four distinct points,
$P_1, P_2, P_3$ and $P_4$ and let $L_0, L_1$ and $L_2$ be lines
through $\{P_1, P_2\},  \{P_1, P_3\}$ and $\{P_1, P_4\}$, respectively.
Choose a point $z_o$ on $C_1$ such that the tangent line at $z_o$ to $C_1$ is not
tangent to $C_2$. Then there exist just three conics $C_3^{(0)}, C_3^{(1)}$ and $C_3^{(2)}$ satisfying
the following conditions:
\begin{itemize}

\item $z_o \in C_3^{(i)}$ for each $i$,
\item For each $i$, $C_3^{(i)}$ is tangent to both $C_1$ and $C_2$ such that
the intersection multiplicities $I_x(C_3^{(i)}, C_j)$ $(j = 1, 2)$ at $\forall
 x \in C_3^{(i)}\cap C_j$ $(j = 1, 2)$  are all even.
\item For $i = 1, 2$, if both of $C_1 + C_2 + C_3^{(i)} + L_1$ and 
$C_1 + C_2 + C_3^{(i)} + L_2$ have the combinatoric for Line-conic arrangement 2 of the
same type, then $(C_1 + C_2 + C_3^{(i)} + L_1, C_1 + C_2 + C_3^{(i)} + L_2)$ is a Zariski pair.
\end{itemize}
\end{enumerate}
 }
\end{thm}

 \begin{rem}{\rm The triple $(C_1 + C_2 + C_3^{(i)} + L_0, C_1 + C_2 + C_3^{(i)} + L_1, C_1 + C_2 + C_3^{(i)} + L_2)$ may be a candidate for a Zariski triple. Our method in this article,
 however, does not work to see whether it is or not.
 }
 \end{rem}


  This article consists of $5$ sections. In \S 1 and \S 2, we summarize some facts and results
    for theory of elliptic surfaces and $D_{2n}$-covers, which we need to 
    prove our theorem. We prove Theorem~\ref{thm:main}  in \S 3 and Theorem~\ref
    {thm:main1} in \S 4.  In \S5, we prove Theorem~\ref{thm:zpair} and give another
    example of a Zariski pair.
  
  \bigskip

  {\bf Acknowledgement.} Some part of this article was done during author's visit   to Universidad de Zaragoza and Ruhr Universit\"at Bochum in 
  September 2011. He thanks for
  Professors E. Artal Bartolo, J.-I. Cogolludo and P. Heinzner for their
  hospitality.

 

\section{$D_{2n}$-covers}

In this section, we summarize some facts on Galois covers.
We refer to \cite{tokunaga94} and  \cite[\S 3]{act} for details.

We start with terminology on Galois covers. Let $X$ and $Y$ be normal projective 
varieties with finite morphism $\pi : X \to Y$.  We say $X$ is a Galois cover of $Y$ if
the induced field extension $\CC(X)/\CC(Y)$ by $\pi^*$ is Galois, where
$\CC(\bullet)$ means the rational function field of $\bullet$. 
Note that the Galois group acts on $X$ such that $Y$ is obtained as the quotient
space with respect to this action ({\it cf.} \cite[\S 1]{tokunaga97}).
 If the Galois group 
$\Gal(\CC(X)/\CC(Y))$ is isomorphic to a finite group $G$, we call $X$ a $G$-cover of $Y$.
The branch locus of $\pi : X \to Y$, which we denote by $\Delta_{\pi}$ or $\Delta(X/Y)$,
is the subset of $Y$ consisting of points $y$ of $Y$, over which $\pi$ is not locally
isomorphic. It is well-known that $\Delta_{\pi}$ is an algebraic subset of pure codimension $1$
if $Y$ is smooth (\cite{zariski}).  

 Suppose that $Y$ is smooth. Let $B$ be a reduced divisor on $Y$ with irreducible
 decomposition $B = \sum_{i=1}^rB_i$. A $G$-cover $\pi : X \to Y$ is said to be branched
 at $\sum_{i=1}^r e_iB_i$ if $(i)$ $\Delta_{\pi} = B$ (here we identify $B$ with its support) and $(ii)$ the ramification index
 along $B_i$ is $e_i$ for each $i$, where  the ramification index mean the
 one along the smooth part of $B_i$ for each $i$.  Note that the study of $G$-covers is related
 to that of Zariski pairs, since we have the following proposition (see \cite{act} for details): 
 
 \begin{prop}\label{prop:fund}{(\cite[Proposition 3.6]{act}) 
 Let $\gamma_i$ be a meridian around $B_i$, and
$[\gamma_i]$ denote its class in the topological fundamental group 
$\pi_1(Y\setminus B, p_o)$.
Let $Y$ be a smooth projective variety and let $B = B_1 + \cdots + B_r$ be 
the decomposition into irreducible components of a reduced divisor $B$ on $Y$.
If there exists a $G$-cover $\pi : X \to Y$ branched at $e_1B_1 + \cdots + e_rB_r$, 
then there exists a normal subgroup $H_{\pi}$ of $\pi_1(Y\setminus B, p_o)$ 
such that:

\begin{enumerate}
\item[(i)] $[\gamma_i]^{e_i} \in H_{\pi},  [\gamma_i]^k \not\in H_{\pi}, (1 \le k \le e_i -1)$, and
\item[(ii)]  $\pi_1(Y\setminus B, p_o)/H_{\pi} \cong G$.
\end{enumerate}

Conversely, if there exists a normal subgroup $H$ of $\pi_1(Y\setminus B, p_o)$ satisfying the above two conditions for $H_{\pi}$, then there exists a $G$-cover
$\pi_H : X_H \to Y$ branched at $e_1B_1 + \cdots + e_rB_r$.
}
\end{prop}

 We keep our notation for $D_{2n}$-covers in Introduction.  
 Here are two propositions for later use.
  
 \begin{prop}\label{prop:dihed-suff}{Let $n$ be an odd integer $\ge 3$. Let $Z$ be a 
 smooth double cover of a smooth projective variety $Y$. We denote its covering morphism
 and covering transformation by $f$ and $\sigma_f$, respectively.
 Let $D$ be an effective divisor on $Z$ satisfying the following condition:
 \medskip
 
 $(i)$ $D$ and $\sigma_f^*D$ have no common component.
 
 \medskip
 
 $(ii)$ If $D = \sum_ia_iD_i$ denotes its irreducible decomposition, then
  $\gcd(a_i, n) = 1$ for every $i$.
  
  \medskip
  
  $(iii)$ $D - \sigma^*_fD$ is $n$-divisible in $\Pic(Z)$.
  
  \medskip
 Then there exists a $D_{2n}$-cover $\pi : X \to Y$ such that
 
 \medskip
 
 $(a)$ $\beta_2(\pi)$ is branched at $n((D + \sigma_f^*D)_{\red})$.
 
 \medskip
 
 $(b)$ $D(X/Y) = Z$ and $f = \beta_2(\pi)$.
 \medskip

 }
 \end{prop}
 
 \proof By \cite[Proposition~0.4]{tokunaga94}, our statements except the ramification
 indices are straightforward. As for the ramification indices, it follows from the last line
 of the proof of \cite[Proposition~0.4]{tokunaga94}. 
 \proofend

 \begin{prop}\label{prop:dihed-nec}{Let $n$ be an odd integer $\ge 3$. Let 
 $\pi : X \to Y$ be a $D_{2n}$-cover such that both $Y$ and $D(X/Y)$ are smooth.
 Let $\sigma_{\beta_1}$ be the covering transform of $\beta_1(\pi)$. If $\beta_2(\pi)$
 is branched at $n{\mathcal D}$ for some non-zero reduced divisor ${\mathcal D}$ on
 $D(X/Y)$, then there exists an effective divisor $D$, whose irreducible decomposition
 is $\sum_ia_iD_i$ satisfying the following conditions:
 
 \medskip
 
 $(i)$ $D$ and $\sigma_{\beta_1}^*D$ have no common component.
 
 \medskip
 
 $(ii)$ $D - \sigma_{\beta_1}^*D$ is $n$-divisible in $\Pic(D(X/Y))$.
 \medskip
 
 $(iii)$ For every $i$, 
 $\gcd(a_i, n) = 1$.
 
 \medskip

 $(iv)$ ${\mathcal D} = (D + \sigma_{\beta_1}^*D)_{\red}$.
 }
 \end{prop}
 
 \proof The statement essentially follows from Proposition~0.5 and its proof in
  \cite{tokunaga94}. We, however, give another simple proof based on the idea of versal
  $D_{2n}$-covers (see \cite{tokunaga06, tsuchihashi} for versal Galois covers).  
  By \cite{tsuchihashi}, 
  there exists an element $\xi \in \CC(X)$ such that the action of $D_{2n}$ on 
  $\xi$ is given in such a way that:
  \[
  \left \{
  \begin{array}{ccl}
  \xi^{\sigma} & = & \frac 1{\xi} \\
  \xi^{\tau} & = & \zeta_n\xi, \quad \zeta_n = \exp(\frac{2\pi i}n).
  \end{array}
  \right.
  \]
  By using $\xi$, we have $\CC(D(X/Y)) = \CC(Y)(\xi^n), \CC(X) = \CC(Y)(\xi)$.  Put 
  $\theta = \xi^n \in \CC(D(X/Y))$. Let $(\theta), (\theta)_0$ and $(\theta)_{\infty}$ be the divisor of $\theta$,
  the zero and polar divisors of $\theta$, respectively.
 Write $(\theta)_0$ in such a way that
  \[
  (\theta)_0 = \sum_ia_iD_i + nD',
\]
where $D_i$'s are irreducible divisor on $D(X/Y)$ with $1 \le a_i < n$ and $D'$ is
an effective divisor on $D(X/Y)$. Since $\sigma$ induces $\sigma_{\beta_1}$ on 
$D(X/Y)$ and $\theta^{\sigma}(= \theta^{\sigma_{\beta_1}})
= 1/\theta$, we have equalities of divisors:
\begin{eqnarray*}
(\theta)_{\infty} & = & \sum_ia_i\sigma_{\beta_1}^*D_i + n\sigma_{\beta_1}^*D' \\
(\theta) & = & (\varphi)_0 - (\varphi)_{\infty} \\
            & = & \sum_ia_i(D_i - \sigma_{\beta_1}^*D_i) + n(D' - \sigma_{\beta_1}^*D').
\end{eqnarray*}
Now we put $D = \sum_ia_iD_i$.
Since  we may assume that $(\theta)_0$ and $(\theta)_{\infty}$ have no common
components, our statements $(i)$ and $(ii)$ follow. As $X$ is the 
$\CC(D(X/Y))(\sqrt[n]{\theta})$-normalization of $D(X/Y)$ and the ramification index
along $D_i$ is $n/\gcd(a_i, n)$, our statements $(iii)$ and $(iv)$ follow. 

\proofend
 
 \begin{cor}\label{cor:split}{Under the same assumption of 
 Proposition~\ref{prop:dihed-nec}, if $D$ is an irreducible divisor
 on $Y$ such that $(\beta_1(\pi))^{-1}(D) \subset \Delta_{\beta_2(\pi)}$,
 then $\beta_1(\pi)^*D$ consists of two irreducible components. 
 In particular, in the case of $\dim Y = 2$, $I_x(D, \Delta_{\beta_1(\pi)})$
 is even for $\forall x \in D\cap \Delta_{\beta_1(\pi)}$
 }
 \end{cor}
 
 \proof The first statement is immediate from Proposition~\ref{prop:dihed-nec}. For the second statement, let $\widetilde{D}$ be the normalization of
 $D$. If there exists $x \in D\cap \Delta_{\beta_1(\pi)}$ such that
 $I_x(D, \Delta_{\beta_1(\pi)})$ is odd, $\beta_1(\pi)$ induces
 a branched double cover of $\widetilde{D}$. This means $\beta_1(\pi)^*D$
 is irreducible.
 \proofend


\section{Elliptic surfaces}

\subsection{General Facts}
We first summarize some facts on the theory of elliptic surfaces. As for details, we refer to 
\cite{kodaira}, \cite{miranda-basic}, \cite{miranda-persson} and \cite{shioda90}.

 In this article, the term, an {\it elliptic surface}, always means a smooth projective surface $S$ equipped with 
 a structure of a fiber space $\varphi : S \to C$ over a smooth projective curve, $C$, as follows:
 \begin{enumerate}
  \item[(i)] There exists a non-empty finite subset, $\Sing(\varphi)$, of $C$ such
  that $\varphi^{-1}(v)$ is a smooth curve of genus $1$ for $v \in C\setminus \Sing(\varphi)$,
while $\varphi^{-1}(v)$ is not a smooth curve of genus $1$ for $v \in \Sing(\varphi)$.
 \item[(ii)]  $\varphi$ has a section
 $O : C \to S$ (we identify $O$ with its image). 
 \item[(iii)]  There is no exceptional
 curve of the first kind in any fiber. 
 \end{enumerate}
 
 For $v \in \Sing(\varphi)$, we put $F_v = \varphi^{-1}(v)$. 
 We denote its irreducible decomposition by 
 \[
 F_v = \Theta_{v, 0} + \sum_{i=1}^{m_v-1} a_{v,i}\Theta_{v,i}
 \]
 where $m_v$ is the number of irreducible components of $F_v$ and $\Theta_{v,0}$ is the
 irreducible component with $\Theta_{v,0}O = 1$. We call $\Theta_{v,0}$ the {\it identity
 component}.  The types of singular fibers are classified as follows (\cite{kodaira}):

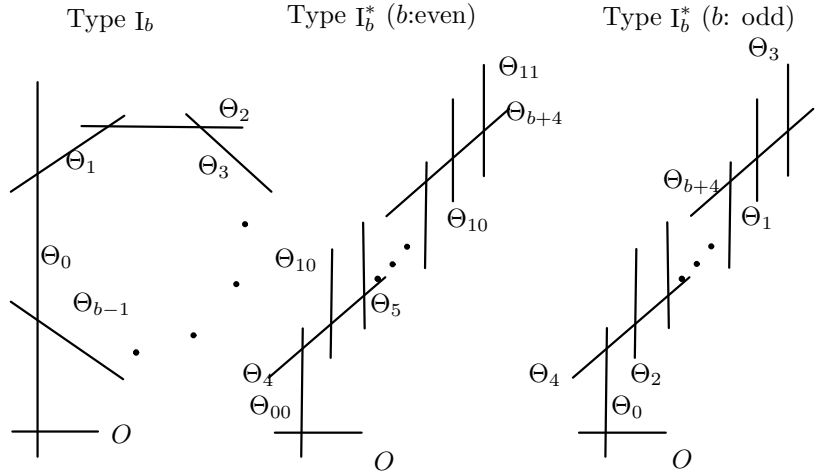
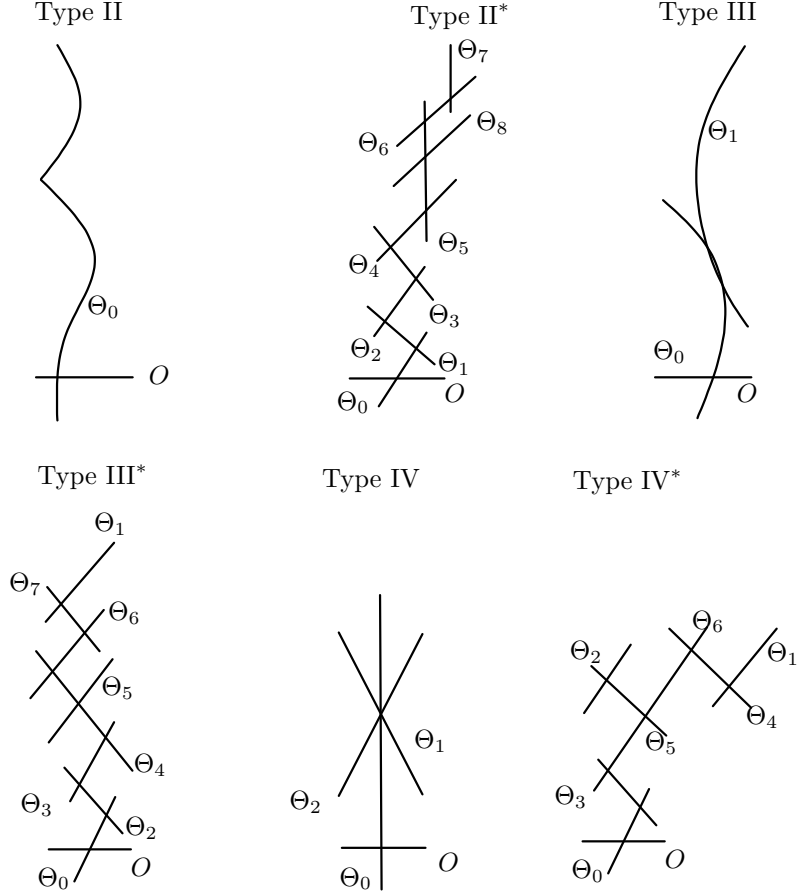
 \begin{center}
 
\unitlength 0.1in
\begin{picture}( 42.0000, 24.1000)(  2.2000,-27.7000)
%
\special{pn 13}%
\special{pa 360 808}%
\special{pa 360 2770}%
\special{fp}%
%
\special{pn 13}%
\special{pa 228 2638}%
\special{pa 676 2638}%
\special{fp}%
%
\special{pn 13}%
\special{pa 220 1382}%
\special{pa 814 984}%
\special{fp}%
%
\special{pn 13}%
\special{pa 590 1038}%
\special{pa 1432 1046}%
\special{fp}%
%
\special{pn 13}%
\special{pa 1138 976}%
\special{pa 1584 1382}%
\special{fp}%
%
\special{pn 13}%
\special{pa 220 1958}%
\special{pa 808 2364}%
\special{fp}%
%
\special{pn 13}%
\special{sh 1}%
\special{ar 1446 1552 10 10 0  6.28318530717959E+0000}%
\special{sh 1}%
\special{ar 1400 1866 10 10 0  6.28318530717959E+0000}%
\special{sh 1}%
\special{ar 1176 2134 10 10 0  6.28318530717959E+0000}%
\special{sh 1}%
\special{ar 876 2224 10 10 0  6.28318530717959E+0000}%
\special{sh 1}%
\special{ar 876 2224 10 10 0  6.28318530717959E+0000}%
%
\special{pn 13}%
\special{pa 1748 2098}%
\special{pa 1740 2770}%
\special{fp}%
%
\special{pn 13}%
\special{pa 1602 2644}%
\special{pa 2056 2644}%
\special{fp}%
%
\special{pn 13}%
\special{pa 1570 2356}%
\special{pa 2180 1830}%
\special{fp}%
%
\special{pn 13}%
\special{pa 1902 1684}%
\special{pa 1894 2252}%
\special{fp}%
%
\special{pn 13}%
\special{pa 2064 1544}%
\special{pa 2072 2098}%
\special{fp}%
%
\special{pn 13}%
\special{pa 2396 1228}%
\special{pa 2388 1782}%
\special{fp}%
%
\special{pn 13}%
\special{pa 2186 1510}%
\special{pa 2828 948}%
\special{fp}%
%
\special{pn 13}%
\special{pa 2534 898}%
\special{pa 2534 1432}%
\special{fp}%
%
\special{pn 13}%
\special{pa 2696 718}%
\special{pa 2696 1300}%
\special{fp}%
%
\special{pn 13}%
\special{sh 1}%
\special{ar 2294 1670 10 10 0  6.28318530717959E+0000}%
\special{sh 1}%
\special{ar 2218 1762 10 10 0  6.28318530717959E+0000}%
\special{sh 1}%
\special{ar 2140 1838 10 10 0  6.28318530717959E+0000}%
\special{sh 1}%
\special{ar 2140 1838 10 10 0  6.28318530717959E+0000}%
%
\special{pn 13}%
\special{pa 3340 2096}%
\special{pa 3332 2770}%
\special{fp}%
%
\special{pn 13}%
\special{pa 3194 2644}%
\special{pa 3648 2644}%
\special{fp}%
%
\special{pn 13}%
\special{pa 3162 2356}%
\special{pa 3772 1830}%
\special{fp}%
%
\special{pn 13}%
\special{pa 3494 1682}%
\special{pa 3486 2252}%
\special{fp}%
%
\special{pn 13}%
\special{pa 3656 1542}%
\special{pa 3664 2096}%
\special{fp}%
%
\special{pn 13}%
\special{pa 3988 1228}%
\special{pa 3980 1780}%
\special{fp}%
%
\special{pn 13}%
\special{pa 3780 1508}%
\special{pa 4420 948}%
\special{fp}%
%
\special{pn 13}%
\special{pa 4126 898}%
\special{pa 4126 1430}%
\special{fp}%
%
\special{pn 13}%
\special{pa 4288 716}%
\special{pa 4288 1298}%
\special{fp}%
%
\special{pn 13}%
\special{sh 1}%
\special{ar 3886 1668 10 10 0  6.28318530717959E+0000}%
\special{sh 1}%
\special{ar 3810 1760 10 10 0  6.28318530717959E+0000}%
\special{sh 1}%
\special{ar 3732 1838 10 10 0  6.28318530717959E+0000}%
\special{sh 1}%
\special{ar 3732 1838 10 10 0  6.28318530717959E+0000}%
\put(21.2000,-28.4000){\makebox(0,0)[lb]{$O$}}%
\put(3.7400,-17.8100){\makebox(0,0)[lb]{$\Theta_0$}}%
\put(5.4400,-20.4100){\makebox(0,0)[lb]{$\Theta_{b-1}$}}%
\put(5.0600,-12.8400){\makebox(0,0)[lb]{$\Theta_1$}}%
\put(13.0700,-10.1100){\makebox(0,0)[lb]{$\Theta_2$}}%
\put(11.9900,-13.1300){\makebox(0,0)[lb]{$\Theta_3$}}%
\put(5.1200,-5.4800){\makebox(0,0)[lb]{Type $\I_b$}}%
\put(7.4500,-27.0000){\makebox(0,0)[lb]{$O$}}%
\put(14.6900,-25.6700){\makebox(0,0)[lb]{$\Theta_{00}$}}%
\put(16.0000,-18.0000){\makebox(0,0)[lb]{$\Theta_{10}$}}%
\put(25.0000,-15.9000){\makebox(0,0)[lb]{$\Theta_{10}$}}%
\put(27.7000,-7.9000){\makebox(0,0)[lb]{$\Theta_{11}$}}%
\put(14.3100,-23.9900){\makebox(0,0)[lb]{$\Theta_4$}}%
\put(21.0100,-20.3400){\makebox(0,0)[lb]{$\Theta_5$}}%
\put(28.1000,-10.4000){\makebox(0,0)[lb]{$\Theta_{b+4}$}}%
\put(33.6200,-25.6600){\makebox(0,0)[lb]{$\Theta_{0}$}}%
\put(34.6200,-23.9800){\makebox(0,0)[lb]{$\Theta_{2}$}}%
\put(40.4100,-15.7800){\makebox(0,0)[lb]{$\Theta_{1}$}}%
\put(40.7900,-6.9500){\makebox(0,0)[lb]{$\Theta_3$}}%
\put(36.3300,-13.9600){\makebox(0,0)[lb]{$\Theta_{b+4}$}}%
\put(29.3000,-24.0000){\makebox(0,0)[lb]{$\Theta_{4}$}}%
\put(16.6600,-5.4000){\makebox(0,0)[lb]{Type $\I_b^*$ ($b$:even)}}%
\put(33.1600,-5.4700){\makebox(0,0)[lb]{Type $\I_b^*$ ($b$: odd)}}%
\put(36.8000,-28.3000){\makebox(0,0)[lb]{$O$}}%
\end{picture}%

 \bigskip

\unitlength 0.1in
\begin{picture}( 37.5100, 22.5000)(  3.9000,-25.9000)
%
\special{pn 13}%
\special{pa 2040 2370}%
\special{pa 2536 2370}%
\special{fp}%
\put(23.5500,-5.3800){\makebox(0,0)[lb]{Type $\mbox{II}^*$}}%
%
\special{pn 13}%
\special{pa 2442 2130}%
\special{pa 2192 2516}%
\special{fp}%
%
\special{pn 13}%
\special{pa 2134 1994}%
\special{pa 2484 2296}%
\special{fp}%
%
\special{pn 13}%
\special{pa 2168 2148}%
\special{pa 2432 1786}%
\special{fp}%
%
\special{pn 13}%
\special{pa 2168 1576}%
\special{pa 2476 1958}%
\special{fp}%
%
\special{pn 13}%
\special{pa 2184 1754}%
\special{pa 2596 1330}%
\special{fp}%
%
\special{pn 13}%
\special{pa 2432 920}%
\special{pa 2442 1650}%
\special{fp}%
%
\special{pn 13}%
\special{pa 2270 1362}%
\special{pa 2670 992}%
\special{fp}%
%
\special{pn 13}%
\special{pa 2288 1152}%
\special{pa 2698 790}%
\special{fp}%
%
\special{pn 13}%
\special{pa 2568 624}%
\special{pa 2568 974}%
\special{fp}%
%
\special{pn 13}%
\special{pa 398 2364}%
\special{pa 904 2364}%
\special{fp}%
%
\special{pn 13}%
\special{pa 510 624}%
\special{pa 528 654}%
\special{pa 544 684}%
\special{pa 560 714}%
\special{pa 576 744}%
\special{pa 590 774}%
\special{pa 604 802}%
\special{pa 614 832}%
\special{pa 622 862}%
\special{pa 628 892}%
\special{pa 630 920}%
\special{pa 630 950}%
\special{pa 626 978}%
\special{pa 620 1008}%
\special{pa 610 1036}%
\special{pa 600 1064}%
\special{pa 586 1094}%
\special{pa 570 1122}%
\special{pa 554 1150}%
\special{pa 534 1178}%
\special{pa 516 1208}%
\special{pa 494 1236}%
\special{pa 472 1264}%
\special{pa 452 1292}%
\special{pa 428 1320}%
\special{pa 424 1326}%
\special{sp}%
%
\special{pn 13}%
\special{pa 424 1330}%
\special{pa 452 1356}%
\special{pa 478 1382}%
\special{pa 504 1410}%
\special{pa 530 1436}%
\special{pa 554 1462}%
\special{pa 578 1488}%
\special{pa 600 1514}%
\special{pa 620 1542}%
\special{pa 640 1568}%
\special{pa 658 1594}%
\special{pa 672 1622}%
\special{pa 684 1650}%
\special{pa 694 1676}%
\special{pa 702 1704}%
\special{pa 706 1732}%
\special{pa 706 1760}%
\special{pa 704 1788}%
\special{pa 698 1816}%
\special{pa 690 1846}%
\special{pa 680 1874}%
\special{pa 668 1904}%
\special{pa 654 1934}%
\special{pa 638 1964}%
\special{pa 622 1994}%
\special{pa 608 2024}%
\special{pa 592 2054}%
\special{pa 578 2084}%
\special{pa 564 2114}%
\special{pa 552 2146}%
\special{pa 542 2176}%
\special{pa 534 2208}%
\special{pa 526 2240}%
\special{pa 520 2272}%
\special{pa 516 2302}%
\special{pa 512 2334}%
\special{pa 510 2366}%
\special{pa 508 2398}%
\special{pa 508 2430}%
\special{pa 508 2462}%
\special{pa 508 2496}%
\special{pa 508 2528}%
\special{pa 508 2560}%
\special{pa 510 2590}%
\special{sp}%
%
\special{pn 13}%
\special{pa 3638 2364}%
\special{pa 4142 2364}%
\special{fp}%
%
\special{pn 13}%
\special{pa 3680 1436}%
\special{pa 3704 1458}%
\special{pa 3730 1480}%
\special{pa 3754 1502}%
\special{pa 3776 1524}%
\special{pa 3800 1546}%
\special{pa 3822 1570}%
\special{pa 3844 1594}%
\special{pa 3864 1618}%
\special{pa 3884 1642}%
\special{pa 3902 1668}%
\special{pa 3920 1696}%
\special{pa 3936 1722}%
\special{pa 3950 1752}%
\special{pa 3962 1782}%
\special{pa 3972 1812}%
\special{pa 3982 1842}%
\special{pa 3990 1874}%
\special{pa 3996 1906}%
\special{pa 4000 1938}%
\special{pa 4004 1970}%
\special{pa 4006 2002}%
\special{pa 4006 2036}%
\special{pa 4004 2068}%
\special{pa 4002 2100}%
\special{pa 4000 2132}%
\special{pa 3994 2162}%
\special{pa 3988 2194}%
\special{pa 3982 2226}%
\special{pa 3974 2256}%
\special{pa 3966 2288}%
\special{pa 3958 2318}%
\special{pa 3948 2348}%
\special{pa 3938 2378}%
\special{pa 3926 2410}%
\special{pa 3914 2440}%
\special{pa 3904 2470}%
\special{pa 3892 2500}%
\special{pa 3880 2530}%
\special{pa 3866 2560}%
\special{pa 3860 2578}%
\special{sp}%
%
\special{pn 13}%
\special{pa 4108 630}%
\special{pa 4090 658}%
\special{pa 4072 686}%
\special{pa 4056 714}%
\special{pa 4038 742}%
\special{pa 4022 770}%
\special{pa 4006 798}%
\special{pa 3990 826}%
\special{pa 3974 854}%
\special{pa 3960 882}%
\special{pa 3946 912}%
\special{pa 3932 940}%
\special{pa 3920 970}%
\special{pa 3908 1000}%
\special{pa 3898 1030}%
\special{pa 3888 1060}%
\special{pa 3880 1090}%
\special{pa 3872 1120}%
\special{pa 3866 1152}%
\special{pa 3862 1184}%
\special{pa 3858 1214}%
\special{pa 3856 1246}%
\special{pa 3854 1278}%
\special{pa 3854 1310}%
\special{pa 3854 1344}%
\special{pa 3856 1376}%
\special{pa 3858 1408}%
\special{pa 3862 1440}%
\special{pa 3866 1472}%
\special{pa 3870 1506}%
\special{pa 3876 1538}%
\special{pa 3882 1570}%
\special{pa 3890 1600}%
\special{pa 3898 1632}%
\special{pa 3906 1664}%
\special{pa 3916 1694}%
\special{pa 3926 1726}%
\special{pa 3936 1756}%
\special{pa 3948 1786}%
\special{pa 3960 1816}%
\special{pa 3974 1846}%
\special{pa 3988 1874}%
\special{pa 4002 1904}%
\special{pa 4016 1932}%
\special{pa 4032 1960}%
\special{pa 4048 1988}%
\special{pa 4066 2014}%
\special{pa 4084 2042}%
\special{pa 4102 2068}%
\special{pa 4120 2094}%
\special{pa 4124 2098}%
\special{sp}%
\put(3.9000,-5.2000){\makebox(0,0)[lb]{Type II}}%
\put(36.5500,-5.2000){\makebox(0,0)[lb]{Type III}}%
\put(9.8700,-24.0000){\makebox(0,0)[lb]{$O$}}%
\put(6.6400,-20.5500){\makebox(0,0)[lb]{$\Theta_0$}}%
\put(25.3500,-24.9100){\makebox(0,0)[lb]{$O$}}%
\put(19.7100,-25.5400){\makebox(0,0)[lb]{$\Theta_0$}}%
\put(40.6500,-24.9100){\makebox(0,0)[lb]{$O$}}%
\put(36.2800,-22.8900){\makebox(0,0)[lb]{$\Theta_0$}}%
\put(39.1000,-11.3400){\makebox(0,0)[lb]{$\Theta_1$}}%
\put(25.0000,-23.5000){\makebox(0,0)[lb]{$\Theta_1$}}%
\put(20.4000,-22.8000){\makebox(0,0)[lb]{$\Theta_2$}}%
\put(24.5000,-21.0000){\makebox(0,0)[lb]{$\Theta_3$}}%
\put(20.3000,-18.4000){\makebox(0,0)[lb]{$\Theta_4$}}%
\put(25.0000,-17.2000){\makebox(0,0)[lb]{$\Theta_5$}}%
\put(20.8000,-12.0000){\makebox(0,0)[lb]{$\Theta_6$}}%
\put(26.0000,-7.1000){\makebox(0,0)[lb]{$\Theta_7$}}%
\put(27.0000,-10.9000){\makebox(0,0)[lb]{$\Theta_8$}}%
\end{picture}%

\bigskip

\unitlength 0.1in
\begin{picture}( 40.1800, 22.7700)(  3.1000,-25.6000)
%
\special{pn 13}%
\special{pa 518 2350}%
\special{pa 948 2350}%
\special{fp}%
\put(4.5800,-4.6300){\makebox(0,0)[lb]{Type $\mbox{III}^*$}}%
%
\special{pn 13}%
\special{pa 866 2076}%
\special{pa 652 2518}%
\special{fp}%
%
\special{pn 13}%
\special{pa 600 1922}%
\special{pa 904 2266}%
\special{fp}%
%
\special{pn 13}%
\special{pa 630 2098}%
\special{pa 860 1684}%
\special{fp}%
\put(9.4700,-24.9000){\makebox(0,0)[lb]{$O$}}%
\put(4.5800,-25.6000){\makebox(0,0)[lb]{$\Theta_0$}}%
%
\special{pn 13}%
\special{pa 452 1312}%
\special{pa 956 1938}%
\special{fp}%
%
\special{pn 13}%
\special{pa 518 1790}%
\special{pa 852 1354}%
\special{fp}%
%
\special{pn 13}%
\special{pa 806 1096}%
\special{pa 422 1558}%
\special{fp}%
%
\special{pn 13}%
\special{pa 510 976}%
\special{pa 784 1312}%
\special{fp}%
%
\special{pn 13}%
\special{pa 502 1158}%
\special{pa 860 744}%
\special{fp}%
\put(7.5400,-7.0300){\makebox(0,0)[lb]{$\Theta_1$}}%
\put(9.1700,-22.9400){\makebox(0,0)[lb]{$\Theta_2$}}%
\put(3.7000,-21.7400){\makebox(0,0)[lb]{$\Theta_3$}}%
\put(9.6200,-19.6500){\makebox(0,0)[lb]{$\Theta_4$}}%
\put(7.9200,-15.6400){\makebox(0,0)[lb]{$\Theta_5$}}%
\put(8.2800,-11.7100){\makebox(0,0)[lb]{$\Theta_6$}}%
\put(3.1000,-10.1800){\makebox(0,0)[lb]{$\Theta_7$}}%
%
\special{pn 13}%
\special{pa 2052 2342}%
\special{pa 2490 2342}%
\special{fp}%
%
\special{pn 13}%
\special{pa 2252 1018}%
\special{pa 2260 2560}%
\special{fp}%
%
\special{pn 13}%
\special{pa 2036 1214}%
\special{pa 2474 2062}%
\special{fp}%
%
\special{pn 13}%
\special{pa 2474 1222}%
\special{pa 2036 2070}%
\special{fp}%
\put(25.5600,-24.6200){\makebox(0,0)[lb]{$O$}}%
\put(20.4000,-25.6000){\makebox(0,0)[lb]{$\Theta_0$}}%
\put(24.2900,-18.3800){\makebox(0,0)[lb]{$\Theta_1$}}%
\put(17.7700,-21.6000){\makebox(0,0)[lb]{$\Theta_2$}}%
%
\special{pn 13}%
\special{pa 3312 2308}%
\special{pa 3742 2308}%
\special{fp}%
%
\special{pn 13}%
\special{pa 3660 2034}%
\special{pa 3446 2476}%
\special{fp}%
%
\special{pn 13}%
\special{pa 3394 1880}%
\special{pa 3698 2224}%
\special{fp}%
\put(37.4200,-24.4800){\makebox(0,0)[lb]{$O$}}%
\put(32.5300,-25.1800){\makebox(0,0)[lb]{$\Theta_0$}}%
\put(42.7600,-13.8200){\makebox(0,0)[lb]{$\Theta_1$}}%
\put(32.4500,-13.7500){\makebox(0,0)[lb]{$\Theta_2$}}%
\put(31.6400,-21.3200){\makebox(0,0)[lb]{$\Theta_3$}}%
\put(41.6500,-17.1200){\makebox(0,0)[lb]{$\Theta_4$}}%
\put(36.3800,-18.5200){\makebox(0,0)[lb]{$\Theta_5$}}%
\put(38.8300,-12.1400){\makebox(0,0)[lb]{$\Theta_6$}}%
%
\special{pn 13}%
\special{pa 3964 1186}%
\special{pa 3372 2042}%
\special{fp}%
%
\special{pn 13}%
\special{pa 3358 1390}%
\special{pa 3750 1754}%
\special{fp}%
%
\special{pn 13}%
\special{pa 3564 1278}%
\special{pa 3320 1636}%
\special{fp}%
%
\special{pn 13}%
\special{pa 3766 1194}%
\special{pa 4194 1608}%
\special{fp}%
%
\special{pn 13}%
\special{pa 4328 1194}%
\special{pa 3994 1600}%
\special{fp}%
\put(19.4800,-4.7800){\makebox(0,0)[lb]{Type $\mbox{IV}$}}%
\put(32.6000,-4.9200){\makebox(0,0)[lb]{Type $\mbox{IV}^*$}}%
\end{picture}%

\end{center} 
 Note that every
 smooth irreducible component is a rational curve with self-intersection number $-2$.
 
 We also define a subset of $\Sing(\varphi)$ by
 $\Red(\varphi) := \{v \in \Sing(\varphi)\mid \mbox{$F_v$ is reducible}\}$. Let 
 $\MW(S)$ be the set of sections of $\varphi : S \to  C$. From our assumption, 
 $\MW(S) \neq  \emptyset$. By regarding $O$ as the zero element of 
 $\MW(S)$ and 
 considering fiberwise addition (see \cite[\S 9]{kodaira} or \cite[\S 1]{tokunaga11} for
 the addition on singular fibers), 
 $\MW(S)$ becomes an abelian group. We denote its addition by $\dot{+}$.
 
 Also for $k \in \ZZ$ and $s \in \MW(S)$, we write
 \[
 [k]s := \left \{ \begin{array}{l}
                  \mbox{$k$-times addition of $s$ if $k \geq 0$} \\
                  \mbox{$k$-times addition of the inverse of $s$ if $k < 0$}.
                  \end{array} \right .
\]
Let $\NS(S)$ be the N\'eron-Severi group of $S$ and let $T_{\varphi}$ be the 
subgroup of $\NS(S)$ generated by $O, F$ and  $\Theta_{v,i}$ $(v \in \Red(\varphi)$,
$1 \le i \le m_v-1)$. Then we have the following theorems:

\begin{thm}\label{thm:shioda-basic0}{\cite[Theorem~1.2]{shioda90}) Under our assumption,
$\NS(S)$ is torsion free.
}
\end{thm}

\begin{thm}\label{thm:shioda-basic}{(\cite[Theorem~1.3]{shioda90}) Under our assumption,
there is a natural map $\tilde{\psi} : \NS(S) \to \MW(S)$ which induces an isomorphisms of 
groups
\[
\psi : \NS(S)/T_{\varphi} \cong \MW(S).
\]
In particular, $\MW(S)$ is a finitely generated abelian group.
}
\end{thm}

In the following,  by the rank, $\rank\MW(S)$, of $\MW(S)$, we mean that of the free part 
of $\MW(S)$.
By Theorem~\ref{thm:shioda-basic}, we may regard the addition of two sections in
$\NS(S)$ as that in $\MW(S)$. We, however, use the notation $\dot{+}$ in order to 
distinguish them. For a divisor on $S$, we  put $s(D) = \psi(D)$. Then we have

\begin{lem}\label{lem:fund-relation}{(\cite[Lemma~5.1]{shioda90})
$D$ is uniquely written in the form:
\[
D \approx s(D) + (d-1)O + nF + \sum_{v\in \Red(\varphi)}\sum_{i=1}^{m_v-1}b_{v,i}\Theta_{v,i},
\]
where $\approx$ denotes the algebraic equivalence of divisors, and $d, n$ and $b_{v,i}$
are integers defined as follows:
\[
d = DF \qquad n = (d-1)\chi({\mathcal O}_S) + OD - s(D)O,
\]
and
\[
\left ( \begin{array}{c}
          b_{v,1} \\
          \vdots \\
          b_{v, m_v-1} \end{array} \right ) = A_v^{-1}\left (\begin{array}{c}
                                                                          D\Theta_{v,1} - s_D\Theta_{v,1} \\
                                                                          \vdots \\
                                                                      D\Theta_{v,m_v-1} - s_D\Theta_{v, m_v-1}
                                                                      \end{array} \right )
\]
Here $A_v$ is the intersection matrix $(\Theta_{v,i}\Theta_{v, j})_{1\le i, j \le m_v-1}$.  
}
\end{lem}

For a proof, see \cite{shioda90}.

\medskip

Put $\NS_{\QQ}:= \NS(S)\otimes \QQ$ and $T_{\varphi, \QQ} := T_{\varphi}\otimes \QQ$.
Since $\NS(S)$ is torsion free under our setting, there is no harm in considering $\NS_{\QQ}$.
By using the intersection pairing, we have the orthogonal decomposition $\NS_{\QQ}
= T_{\varphi, \QQ}\oplus(T_{\varphi, \QQ})^{\perp}$.   In \cite{shioda90}, the homomorphism
$\phi: \MW(S) \to (T_{\varphi, \QQ})^{\perp}\subset\NS_{\QQ}$ is defined as follows:
\begin{eqnarray*}
\phi : \MW(S) \ni s &\mapsto & s - O - (sO + \chi({\mathcal O}_S))F \\
&&  +\sum_{v \in \Red(\varphi)}
(\Theta_{v,1}, \ldots, \Theta_{v, m_v-1})(-A_v)^{-1}
\left ( \begin{array}{c}
        s\Theta_{v,1} \\
        \vdots \\
        s\Theta_{v, m_v-1}
        \end{array} \right ) \in (T_{\varphi, \QQ})^{\perp}.
\end{eqnarray*}
 Also,  in \cite{shioda90},  a $\QQ$-valued bilinear form $\langle \, , \, \rangle$ on
 $\MW(S)$  is defined by
$\langle s_1, s_2 \rangle := - \phi(s_1)\phi(s_2)$, where the right hand side means the intersection
pairing in $\NS_{\QQ}$.    Here are two basic properties of $\langle \, , \, \rangle$:
\begin{itemize}

\item $\langle s, \, s \rangle \ge 0$ for $\forall s \in \MW(S)$ and the equality holds if and 
only if $s$ is an element of finite order in $\MW(S)$. 

\item An explicit formula for $\langle s_1,
s_2\rangle$ ($s_1, s_2 \in \MW(S)$) is given as follows:
\[
\langle s_1, s_2 \rangle = \chi({\mathcal O}_S) + s_1O + s_2O - s_1s_2 - \sum_{v \in \Red(\varphi)}
\mbox{Contr}_v(s_1, s_2),
\]
where $\mbox{Contr}_v(s_1, s_2)$ is given by
\[
\mbox{Contr}_v(s_1, s_2) = (s_1\Theta_{v,1}, \ldots, s_1\Theta_{v, m_v-1})(-A_v)^{-1}
\left ( \begin{array}{c}
        s_2\Theta_{v,1} \\
        \vdots \\
        s_2\Theta_{v, m_v-1}
        \end{array} \right ).
\]
 As for explicit values of 
$\mbox{Contr}_v(s_1, s_2)$, we refer to \cite[(8.16)]{shioda90}.
\end{itemize}

\subsection{Double cover construction of an elliptic surface}\label{subsec:double-cover}

For details about this subsection, see \cite[Lectures III and IV]{miranda-basic}.
Let $\varphi: S \to C$ be an elliptic surface. By our assumption, the generic fiber of $\varphi$ can be
considered as an elliptic curve over $\CC(C)$, the rational function field of $C$.
The inverse morphism with respect to
the group law induces an involution $[-1]_{\varphi}$ on $S$. Let $S/\langle [-1]_{\varphi}
\rangle$ be the quotient by $[-1]_{\varphi}$. The quotient surface $S/\langle[-1]_{\varphi}\rangle$ 
is known to be smooth and $S/\langle [-1]_{\varphi} \rangle$ can be blown down to its relatively minimal
model $W$ over $C$ satisfying the following condition:

Let us denote
      \begin{itemize}
       \item $f: S \to S/\langle [-1]_{\varphi}\rangle$: the quotient morphism, 
       \item $q: S/\langle [-1]_{\varphi}\rangle \to W$: the blow down, and 
      \item $S \miya{\mu} S' \miya{{f'}} W$: the Stein factorization of $q\circ f$. 
      \end{itemize}
      Then we have:

\begin{enumerate}

 \item  The branch locus $\Delta_{f'}$ of $f'$ consists of a section $\Delta_0$ and the
 trisection $T$ such that its singularities are at most simple singularities (see \cite[Chapter II, \S 8]{bpv} for simple singularities
 and their notation) and $\Delta_0\cap T = \emptyset$. 
 
 \item $\Delta_0 + T$ is $2$-divisible in $\Pic(W)$.
 
 \item     
   The morphism  $\mu$ is obtained by contracting all the irreducible components of 
 singular fibers not meeting $O$.

 \end{enumerate}

Conversely, if $\Delta_0$ and $T$ on $W$ satisfy the above condition, we
obtain an elliptic surface $\varphi : S \to {\mathbb P}^1$, as the canonical
resolution of a double cover $f' : S' \to W$ with $\Delta_{f'} = \Delta_0
+ T$, and the diagram (see \cite{horikawa} for the canonical resolution):
\[
\begin{CD}
S' @<{\mu}<< S \\
@V{f'}VV                 @VV{f}V \\
W@<<{q}< \widehat{W}.
\end{CD}
\]
Here $q$ is a composition of blowing-ups so that $\widehat{W}
= S/\langle [-1]_{\varphi}\rangle$. 
Hence any elliptic surface  is obtained as above.
In the following, we call the diagram above
{\it the double cover diagram for $S$}.

In the case of $C = {\mathbb P}^1$, $W$ is the Hirzebruch surface, 
 $\Sigma_d$, of degree 
 $d = 2\chi({\mathcal O}_S)$ and $\Delta_{f'}$ is of the form $\Delta_0 + T$, where
 $\Delta_0$ is a section with $\Delta_0^2 = -d$ and $T \sim 3(\Delta_0 + d{\mathfrak f})$, 
 ${\mathfrak f}$ being a fiber of the ruling $\Sigma_d \to {\mathbb P}^1$.
\begin{rem}\label{rem:double-cover-construction}{\rm

\begin{itemize}
  \item For each $v \in \Sing(\varphi)$, the type of $\varphi^{-1}(v)$ is determined by the type of singularity of $T$
  on ${\mathfrak f}_v$ and the relative position between ${\mathfrak f}_v$ and $T$ (see \cite[Table 6.2]{miranda-persson}).
  
  \item Note that the covering transformation,
  $\sigma_f$, of $f$ coincides with  $[-1]_{\varphi}$. Also, by \cite[Theorem~9.1]{kodaira},
   the action 
  of $\sigma_f$ on irreducible components of singular fibers is described as follows:
  
  \begin{center}
  \begin{tabular}{|c|c|} \hline
  Type of a singular fiber & The action on irreducible component \\ \hline
  $\I_n$ & $\displaystyle{\begin{array}{lc}
                                     \Theta_0 \mapsto \Theta_0  & \\
                                     \Theta_i \mapsto \Theta_{n-i} & i = 1, \ldots, n-1 
                                     \end{array}}$ \\ \hline
   $\I_n^*$ $(n: even)$ &    $\displaystyle{\begin{array}{lc}
                                     \Theta_i \mapsto \Theta_i  & \forall i \\
                                     \Theta_{ij} \mapsto \Theta_{ij} & \forall i, j 
                                     \end{array}}$ \\ \hline         
   $\I_n^*$ $(n: odd)$ &     $\displaystyle{\begin{array}{cc}
                                     \Theta_i \mapsto \Theta_i  &  i \neq 1, 3 \\
                                     \Theta_1 \mapsto \Theta_3 & \Theta_3 \mapsto \Theta_1
                                     \end{array}}$ \\ \hline                                                            
  $\II,\, \II^*, \,  \III, \, \III^*$ & $\displaystyle{\begin{array}{cc}
                                     \Theta_i \mapsto \Theta_i  &  \forall i \\
                                     \end{array}}$ \\ \hline 
  $\IV$    &            $\displaystyle{\begin{array}{cc}
                                     \Theta_0 \mapsto \Theta_0 &    \\
                                     \Theta_1 \mapsto \Theta_2 & \Theta_2 \mapsto \Theta_1
                                     \end{array}}$ \\ \hline 
  $\IV^*$    &            $\displaystyle{\begin{array}{lc}
                                     \Theta_i \mapsto \Theta_i &    i = 0, 3, 6 \\
                                     \Theta_1 \mapsto \Theta_2 & \Theta_2 \mapsto \Theta_1 \\
                                     \Theta_4 \mapsto \Theta_5 & \Theta_5 \mapsto \Theta_4 \\
                                     \end{array}}$ \\ \hline               
  \end{tabular}   
  \end{center}                                
  \end{itemize}
  }
  \end{rem}


\section{Elliptic $D_{2p}$-cover over rational ruled surface}

Let $\varphi : S \to \PP^1$ be an elliptic surface over $\PP^1$. Let
$S \to \widehat{\Sigma}_d$ be the double cover appearing in the double cover diagram 
for $S$  in
\S\ref{subsec:double-cover}.

We first note that any elliptic $D_{2p}$-cover ($p$: odd prime)
$\pi_p : X_p \to \widehat{\Sigma}_d$ satisfies the following conditions:

\begin{itemize}
 \item $S = D(X_p/\widehat{\Sigma}_d)$ and $\beta_1(\pi_p) = f$.
 
 \item The branch locus of $\beta_2(\pi_p)$ is of the form
 \[
 \mcD + \sigma_f^*\mcD + \Xi + \sigma_f^*\Xi
 \]
 where 
   \begin{enumerate}
   \item all irreducible components of $\mcD$ are horizontal and there is no common component between
   $\mcD$ and $\sigma_f^*\mcD$, and
   
    \item all irreducible component of $\Xi$ are vertical and there is no common component between $\Xi$ and
    $\sigma_f^*\Xi$.
   \end{enumerate}
 \end{itemize}
 
 \begin{rem}{\rm

 Under the above notation, the case when $\mcD = \emptyset$ (resp. $=$ a section) is considered in the author's previous work(\cite{tokunaga94, tokunaga97, tokunaga98, tokunaga99}) 
(resp. \cite{tokunaga04}).

}
\end{rem}

In the following, we always assume that
\begin{center}

$(\ast)$ $\mcD \neq \emptyset$,.

\end{center}

The proposition below, which is a generalization of
 \cite[Propositions 4.1 and 4.2]{tokunaga04}, plays an 
important role in this article:

\begin{thm}\label{thm:main}{ Let $p$ be an odd prime. Let $C_1, \ldots, C_r$ be irreducible horizontal divisors 
on $S$ such that $\sum_{i=1}^r C_i$ and $\sum_{i=1}^r\sigma_f^*C_i$ have no common component. Then 
(I) and (II) in the below are equivalent:

\begin{enumerate}

\item[(I)] Put ${\mathcal C} = \sum_{i=1}^rC_i$. There exists an elliptic $D_{2p}$-cover $\pi_p : X_p \to \widehat{\Sigma}_d$ such
that 
\begin{itemize}
 \item $D(X_p/\widehat{\Sigma}_d) = S$ and $\beta_1(\pi_p) = f$.
 
 \item
  \[
 \Delta_{\beta_2(\pi_p)} = \Supp(\mcC + \sigma_f^*\mcC + \Xi + \sigma_f^*\Xi)
 \]
 where all irreducible components of $\Xi$ are vertical and there is no common component between $\Xi$ and
    $\sigma_f^*\Xi$.
\end{itemize}

\item[(II)] Let $s(C_i) = \tilde{\psi}(C_i)$ $(i = 1,\ldots, r)$. There exist integers $a_i$ $(i = 1,\ldots, r)$ such that
\begin{itemize}
  \item $1 \le a_i < p$ $(i = 1, \ldots, r)$ and
  \item 
  \[
  \sum_{i=1}^r [a_i]s(C_i) \in [p]\MW(S) := \{[p]s \mid s \in \MW(S)\}.
  \]
\end{itemize}  
\end{enumerate}
}
\end{thm}

\proof $(I) \Rightarrow (II)$ Let $D$ be the effective divisor in Proposition~\ref{prop:dihed-nec}. We put $D = D_{\hor} + D_{\ver}$, where the irreducible components of $D_{\hor}$ are all horizontal, while those of $D_{\ver}$ are all in fibers of 
$\varphi$.  By Proposition~\ref{prop:dihed-nec} $(iv)$, 
$(D_{\hor} + \sigma_f^*D_{\hor})_{\red} = \sum_{i=1}^r C_i + \sum_{i=1}^r\sigma_f^*C_i$.  

\medskip

{\bf Claim.} We may assume that
$D_{\hor} = \sum_{i=1}^r a_iC_i$. 

\medskip

{\sl Proof of Claim.} If $\sigma_f^*C_i$ is an irreducible component of $D_{\hor}$, then we consider
\[
D'_{\hor}:= D_{\hor} + (p - a_i)C_i  - a_i\sigma_f^*C_i,
\]
and put $D' = D'_{\hor} + D_{\ver}$. 
Then we infer that $D'$ also satisfies all  four conditions in
 Proposition~\ref{prop:dihed-nec}. After repeating this process finitely many times, we may assume that
 $D_{\hor} = \sum_{i=1}^r a_iC_i$.
 
 \medskip

 By Claim and Proposition~\ref{prop:dihed-nec} $(iii)$, there exists ${\mathcal L} \in \Pic(S)$ such that 
 \[
 \sum_{i=1}^r a_i(C_i - \sigma_f^*C_i) + D_{\ver} - \sigma_f^*D_{\ver} \sim  p{\mathcal L},
 \]
 where $\sim$ means linear equivalence of divisors (Note that linear
 equivalence coincides with algebraic equivalence on $S$). This implies
 \[
 \tilde{\psi}(\sum_{i=1}^r a_i(C_i - \sigma_f^*C_i)) = [p]\tilde{\psi}({\mathcal L}) \quad
 \mbox{ in $\MW(S)$}.
 \]
 As $\tilde{\psi}(\sigma_f^*C_i) = [-1]s(C_i)$, we have
 \[
 \tilde{\psi}(\sum_{i=1}^r a_i(C_i - \sigma_f^*C_i))  = [2]([a_1]s(C_1) \dot{+}\ldots \dot{+} [a_r]s(C_r)).
 \]
 Since $p$ is an odd prime, we infer that $[a_1]s(C_1) \dot{+}\ldots \dot{+} [a_r]s(C_r) \in [p]\MW(S)$.
 
 \medskip
 
 $(II) \Rightarrow (I)$ Our proof is similar to that of \cite[Proposition~4.2]{tokunaga04}. 
 By Lemma~\ref{lem:fund-relation}, we  have 
 \[
 C_i \sim s(C_i) + (d_i - 1)O + n_i F + \sum_{v \in \Red(\varphi)}\sum_{j=1}^{m_v-1}b_{v,j}^{(i)}\Theta_{v,i}.
 \]
 This implies
 \[
 \sum_{i=1}^ra_iC_i \sim \sum_{i=1}^ra_is(C_i) +\sum_{i=1}^ra_i
 \left ((d_i - 1)O + n_i F + \sum_{v \in \Red(\varphi)}\sum_{j=1}^{m_v-1}b_{v,j}^{(i)}\Theta_{v,j}\right ).
 \]
By our assumption, there exists $s_o$ such that $\sum_{i=1}^r[a_i]s(C_i) = [p]s_o$ in
 $\MW(S)$. By Theorem~\ref{thm:shioda-basic}, this implies that
\[
\sum_{i=1}^ra_is(C_i) \sim ps_o + \left (-p +\sum_{i=1}^ra_i  \right )O + n_oF + 
 \sum_{v \in \Red(\varphi)}\sum_{j=1}^{m_v-1}c_{v,j}\Theta_{v,j}
 \]
 for some integers $n_o, c_{v, j}$.  Hence we have
 \begin{eqnarray*}
 \sum_{i=1}^ra_iC_i  &\sim & ps_o + \left (-p + \sum_{i=1}^r a_id_i  \right )O + \left (n_o + \sum_ia_in_i \right ) F \\
 &+&
 \sum_{v \in \Red(\varphi)}\sum_{j=1}^{m_v-1}\left (c_{v,j} + 
 \sum_{i=1}^ra_ib_{v,j}^{(i)}\right )\Theta_{v,j}, 
 \end{eqnarray*}
 and put
  \[
 D' := \sum_{i=1}^r a_iC_i +\sum_{v \in \Red(\varphi)}\sum_{j=1}^{m_v-1}
 \left (c_{v,j} + \sum_{i=1}^ra_ib_{v,j}^{(i)}\right)\sigma_f^*\Theta_{v,j} .
 \]
 Then we have

 \[
 D' - \sigma_f^*D' \sim p\left (s_o - \sigma_f^*s_o \right ).
 \]
The left hand side of the above equivalence contains some redundancy in the sum for $\Theta_{v,i}$ and $\sigma_f^*\Theta_{v,i}$.
By taking the action of $\sigma_f$ on $\Theta_{v,i}$'s (see 
Remark~\ref{rem:double-cover-construction}) into account, 
we can find a divisor $D = \sum_{i=1}^ra_iC_i + \sum_j k_j\Xi_j$ and $\Xi'$ on $S$ such that

\begin{enumerate}
 \item[(i)] all $\Xi_j$ and  all irreducible components  of $\Xi'$
are those in fibers not meeting $O$, 

\item[(ii)] $D$ and $\sigma_f^*D$ have no common component, 

\item[(iii)] $1 \le k_j < p$, and 

\item[(iv)] $D' - \sigma_f^*D' = D - \sigma_f^*D + p\Xi'$.

\end{enumerate}
Now we easily infer that $D$ satisfies the three conditions in Proposition~\ref{prop:dihed-suff}
for $p$. 
\proofend
 
 \bigskip
 
 Fix an isomorphism $\MW(S) = M_o \oplus \MW_{\tor}$, $M_o \cong \ZZ^{\oplus r}$, 
 $r = \rank \MW(S)$. 
 By Theorem~\ref{thm:main}, we have the following proposition:
 
 \begin{prop}\label{prop:2-sections}{Choose $s \in M_o$ such that 
 $M_o/\ZZ s$ is free. For any finite number of odd prime numbers
 $p_1,\ldots, p_l$, there exists a section $s_{p_1,\ldots, p_l}$ satisfying the following conditions:
 
 \begin{enumerate}
  \item[(i)] $\langle s_{p_1,\ldots,p_l},  s_{p_1,\ldots, p_l} \rangle = (p_1\cdots p_l)^2\langle s, s \rangle$.
  
  \item[(ii)] For any prime $p \not\in \{ p_1, \ldots, p_l\}$, there exists an elliptic $D_{2p}$-cover $\pi_p : X_p \to \widehat{\Sigma}_d$
  such that 
     \begin{itemize}
       \item $D(X_p/\widehat{\Sigma}_d) = S, \, \beta_1(\pi_p) = f$, and 
       \item $\beta_2(\pi_p)$ is branched at $p(s + s_{p_1,\ldots, p_l} + \sigma_f^*(s + s_{p_1,\ldots, p_l}) + \Xi_o)$, where all irreducible components
       of $\Xi_o$ are those in fibers not meeting $O$.
      \end{itemize}
  \item[(iii)] For $p \in  \{p_1, \ldots, p_l\}$, there exists no elliptic $D_{2p}$-cover $\pi_p : X_p \to \widehat{\Sigma}_d$ as
  in (ii)
  \item[(iv)] $\{s_{p_1,\ldots, p_l}, [-1]s_{p_1,\ldots, p_l}\}$ is unique up to torsion elements.
  \end{enumerate}
 }
 \end{prop}
 
 \proof Define $s_{p_1, \ldots, p_l}:= [\Pi_{i=1}^rp_i]s$. By Theorem~\ref{thm:main}, our statements (i), (ii) and (iii) are
 immediate.  Suppose that $s' \in \MW(S)$ satisfies the statements (i), (ii) and (iii).  Put $s' = s'_o + t'_o$, 
 $s'_o \in M_o, t'_o \in \MW_{\tor}$. 
 By Theorem~\ref{thm:main}, for $p \not\in \{p_1, \ldots, p_l \}$, 
 there exists an integer $k$ ($1 \le k < p$) such that
 \[
 s \dot{+} [k] s'_o = 0 
 \]
 in $M_o/pM_o$. This implies that there exist integers $l$ and $l'$ with $\gcd(l, l') = 1$ such that $[l]s  \dot{+}  [l']s'_o = 0$
 in $\MW(S)$. Choose integers $m, m'$ with $ml + m'l' = 1$. Then
 \[
 0 = [m]([l] s \dot{+} [l]' s'_o) = s \dot{+} [-m'l']s \dot{+} [ml']
 s'_o = s \dot{+} [l']([-m']s + [m]s'_o).
 \]
 Since $M_o/\ZZ s$ is free, we infer that $l' = 1$ and $s'_o = [-l]s$. Thus 
 \[
 \langle s'_o, \,  s'_o \rangle = l^2\langle s_o, \, s_o \rangle = (p_1\cdots p_l)^2\langle s, \, s \rangle.
 \]
 Since $\langle s, s \rangle \neq 0$ by the
 basic properties of $\langle \, , \, \rangle$ (see \S 1), we have $l = \pm p_1\cdots p_l$. 
 Hence $s'$ is equal to $[\pm 1]s_{p_1,\ldots, p_l}$ up to torsion elements.
 \proofend

 

\section{Applications}

 Let $\varphi : S \to {\mathbb P}^1$ be an elliptic surface and we keep our notation for the double cover construction
 for $S$  in \S\ref{subsec:double-cover}. 
 We fix an isomorphism $\MW(S) \cong M_o\oplus \MW_{\mathrm {tor}}$
 $M_o \cong {\mathbb Z}^{\oplus r}$, $r = \rank \MW(S)$.
 Choose $s_1, s_2 \in M_o$ so that $s_1$ and $s_2$ are a part of a basis
 of ${\mathbb Z}^{\oplus r}$, i.e., $M_o/\ZZ s_1 + \ZZ s_2$ is free of rank $r-2$. 
 Put $s_3 :=[2]s_1$.
 
 \begin{thm}\label{thm:main1}{For any odd prime $p$, there exists an elliptic $D_{2p}$-cover 
 $\pi_p : X_p \to \widehat{\Sigma}_d$  such that the horizontal part of the branch locus, $\Delta_{\beta_2(\pi_p)}$, of 
 $\beta_2(\pi_p)$ is of the form
 \[
 s_1 + s_3 + \sigma_{f}^*(s_1 + s_3),
 \]
  while there exists no elliptic $D_{2p}$-cover 
 $\pi_p : X_p \to \widehat{\Sigma}_d$  such that the horizontal part of the branch locus of 
 $\beta_2(\pi_p)$ is $s_2 + s_3 + \sigma_{f}^*(s_2 + s_3)$ 
 }
 \end{thm}
 
 \proof Since $[p-2]s_1 + s_3 \in [p]\MW(S)$, we have the first statement by Theorem~\ref{thm:main}.
 Since $s_1$ and $s_2$ are a part of basis, their image  in $M_o/pM_o$ are
 linearly independent over $\ZZ/p\ZZ$. Hence our second statement follows from Theorem~\ref{thm:main}.
 \proofend
 
 \bigskip
 
 By Proposition~\ref{prop:fund} and Theorem~\ref{thm:main}, we  have
 
 \begin{cor}\label{cor:zpair}{Let $T$ be the trisection on $\Sigma_d$ appearing 
 in the double cover diagram for $S$.
 Put $\Delta_i := q\circ f(s_i)$ $(i = 1, 2, 3)$. Then there
 exists a $D_{2p}$-cover of $\Sigma_d$ branched at $2(\Delta_0 + T) + p(\Delta_1+ \Delta_3)$, while there exists no $D_{2p}$-cover of $\Sigma_d$
 branched at  $2(\Delta_0 + T) + p(\Delta_2 + \Delta_3)$.
   }
  \end{cor}
  
 We end this section by considering the case when $S$ is a rational elliptic surface.
 In this case, we have the double cover diagram for $S$ as follows:
 \[
\begin{CD}
S' @<{\mu}<< S \\
@V{f'}VV                 @VV{f}V \\
\Sigma_2 @<<{q}< \widehat{\Sigma}_2.
\end{CD}
\]
Write $q:= q_1\circ \cdots \circ q_r : \widehat{\Sigma}_2 =
\Sigma_2^{(r)} \to \cdots \to \Sigma_2^{(1)} \to \Sigma_2^{(0)} = \Sigma_2$, 
where $q_i$ is a blowing up at a point at $\Sigma_2^{(i-1)}$. 
Put $\Delta_{f'} = \Delta_0 + T$. In the following, we assume that 
\medskip
\begin{center}
$T$ has a node $x_o$.
\end{center}
\medskip

Note that this is equivalent to the fact that $S$ has a singular fiber of type $\I_2$ or $\II$
by \cite[Table 6.2]{miranda-persson}.  
We may assume that  $q_1$ is a blowing up at $x_o$. Let $E_1$ be 
the exceptional divisor of $q_1$ and let $\overline{\mathfrak f}_o$  and $\overline{T}$ be the proper transforms 
of a fiber,  ${\mathfrak f}_o$, through $x_o$ and $T$, respectively. Then we have the following picture:

\begin{center}

\unitlength 0.1in
\begin{picture}( 38.9400, 23.5000)(  1.5000,-23.7000)
%
\special{pn 13}%
\special{pa 432 2210}%
\special{pa 1072 2210}%
\special{fp}%
%
\special{pn 13}%
\special{pa 732 1154}%
\special{pa 732 2370}%
\special{fp}%
%
\special{pn 13}%
\special{pa 444 1978}%
\special{pa 1072 1974}%
\special{fp}%
%
\special{pn 13}%
\special{pa 508 1392}%
\special{pa 974 1806}%
\special{fp}%
%
\special{pn 13}%
\special{pa 986 1392}%
\special{pa 478 1806}%
\special{fp}%
%
\special{pn 13}%
\special{pa 1642 1010}%
\special{pa 1096 1420}%
\special{fp}%
\special{sh 1}%
\special{pa 1096 1420}%
\special{pa 1160 1396}%
\special{pa 1138 1388}%
\special{pa 1136 1364}%
\special{pa 1096 1420}%
\special{fp}%
%
\special{pn 13}%
\special{pa 1700 1104}%
\special{pa 2340 1104}%
\special{fp}%
%
\special{pn 13}%
\special{pa 2000 50}%
\special{pa 2000 1264}%
\special{fp}%
%
\special{pn 13}%
\special{pa 1712 874}%
\special{pa 2340 866}%
\special{fp}%
%
\special{pn 13}%
\special{pa 1786 706}%
\special{pa 2604 20}%
\special{fp}%
%
\special{pn 13}%
\special{pa 2068 246}%
\special{pa 2368 550}%
\special{fp}%
%
\special{pn 13}%
\special{pa 2242 130}%
\special{pa 2530 418}%
\special{fp}%
%
\special{pn 8}%
\special{pa 2580 378}%
\special{pa 2600 404}%
\special{pa 2616 430}%
\special{pa 2630 460}%
\special{pa 2638 492}%
\special{pa 2644 526}%
\special{pa 2644 562}%
\special{pa 2642 598}%
\special{pa 2636 636}%
\special{pa 2628 674}%
\special{pa 2614 712}%
\special{pa 2600 748}%
\special{pa 2582 784}%
\special{pa 2562 818}%
\special{pa 2542 850}%
\special{pa 2518 880}%
\special{pa 2492 906}%
\special{pa 2464 928}%
\special{pa 2436 948}%
\special{pa 2406 962}%
\special{pa 2376 972}%
\special{pa 2344 976}%
\special{pa 2314 976}%
\special{pa 2282 970}%
\special{sp}%
%
\special{pn 8}%
\special{pa 450 1306}%
\special{pa 426 1326}%
\special{pa 402 1348}%
\special{pa 380 1372}%
\special{pa 358 1396}%
\special{pa 340 1422}%
\special{pa 322 1448}%
\special{pa 308 1478}%
\special{pa 296 1508}%
\special{pa 286 1540}%
\special{pa 278 1572}%
\special{pa 272 1606}%
\special{pa 270 1638}%
\special{pa 270 1672}%
\special{pa 270 1704}%
\special{pa 274 1738}%
\special{pa 280 1768}%
\special{pa 290 1800}%
\special{pa 300 1830}%
\special{pa 312 1858}%
\special{pa 326 1888}%
\special{pa 340 1916}%
\special{pa 356 1944}%
\special{pa 374 1972}%
\special{pa 392 1998}%
\special{pa 410 2026}%
\special{pa 430 2052}%
\special{pa 438 2066}%
\special{sp}%
%
\special{pn 13}%
\special{pa 2436 1018}%
\special{pa 3008 1438}%
\special{fp}%
\special{sh 1}%
\special{pa 3008 1438}%
\special{pa 2966 1382}%
\special{pa 2964 1406}%
\special{pa 2942 1414}%
\special{pa 3008 1438}%
\special{fp}%
%
\special{pn 13}%
\special{pa 2876 1840}%
\special{pa 2910 1844}%
\special{pa 2942 1846}%
\special{pa 2976 1850}%
\special{pa 3010 1850}%
\special{pa 3042 1852}%
\special{pa 3074 1850}%
\special{pa 3106 1848}%
\special{pa 3138 1844}%
\special{pa 3170 1838}%
\special{pa 3200 1830}%
\special{pa 3228 1820}%
\special{pa 3258 1806}%
\special{pa 3286 1788}%
\special{pa 3312 1770}%
\special{pa 3336 1748}%
\special{pa 3360 1724}%
\special{pa 3380 1698}%
\special{pa 3398 1670}%
\special{pa 3414 1642}%
\special{pa 3426 1612}%
\special{pa 3434 1582}%
\special{pa 3440 1550}%
\special{pa 3442 1520}%
\special{pa 3442 1488}%
\special{pa 3438 1456}%
\special{pa 3434 1422}%
\special{pa 3428 1390}%
\special{pa 3422 1358}%
\special{pa 3414 1324}%
\special{pa 3410 1310}%
\special{sp}%
%
\special{pn 13}%
\special{pa 3024 2002}%
\special{pa 3928 1264}%
\special{fp}%
%
\special{pn 13}%
\special{pa 3606 1282}%
\special{pa 3894 1570}%
\special{fp}%
%
\special{pn 13}%
\special{pa 3756 1196}%
\special{pa 4044 1484}%
\special{fp}%
\put(12.9100,-13.9700){\makebox(0,0)[lb]{$q_1$}}%
\put(1.5000,-17.3700){\makebox(0,0)[lb]{$T$}}%
\put(11.1800,-22.4300){\makebox(0,0)[lb]{$\Delta_0$}}%
\put(26.6100,-7.9800){\makebox(0,0)[lb]{$\overline{T}$}}%
\put(6.7000,-10.6000){\makebox(0,0)[lb]{${\mathfrak {f}}_o$}}%
\put(19.1000,-14.3000){\makebox(0,0)[lb]{$\overline{{\mathfrak {f}}}_o$}}%
\put(16.5300,-5.3800){\makebox(0,0)[lb]{$E_1$}}%
\put(34.5600,-9.9900){\makebox(0,0)[lb]{$\mcQ$}}%
\put(33.3600,-18.1700){\makebox(0,0)[lb]{$z_o$}}%
\put(29.1500,-20.9900){\makebox(0,0)[lb]{$L_{z_o}$}}%
%
\special{pn 8}%
\special{pa 3256 1288}%
\special{pa 3270 1258}%
\special{pa 3284 1230}%
\special{pa 3300 1202}%
\special{pa 3318 1176}%
\special{pa 3338 1152}%
\special{pa 3362 1130}%
\special{pa 3386 1110}%
\special{pa 3414 1092}%
\special{pa 3444 1076}%
\special{pa 3474 1064}%
\special{pa 3506 1054}%
\special{pa 3538 1046}%
\special{pa 3572 1040}%
\special{pa 3604 1036}%
\special{pa 3638 1034}%
\special{pa 3670 1036}%
\special{pa 3702 1038}%
\special{pa 3734 1042}%
\special{pa 3766 1048}%
\special{pa 3796 1054}%
\special{pa 3828 1062}%
\special{pa 3854 1068}%
\special{sp}%
\put(21.2000,-12.2400){\makebox(0,0)[lb]{$\Delta_0$}}%
\put(28.8000,-12.0000){\makebox(0,0)[lb]{$\overline{q}_1$}}%
\end{picture}%

\medskip

The case $(a)$

\bigskip

\unitlength 0.1in
\begin{picture}( 38.9400, 23.5000)(  1.5000,-23.7000)
%
\special{pn 13}%
\special{pa 432 2210}%
\special{pa 1072 2210}%
\special{fp}%
%
\special{pn 13}%
\special{pa 732 1154}%
\special{pa 732 2370}%
\special{fp}%
%
\special{pn 13}%
\special{pa 444 1978}%
\special{pa 1072 1974}%
\special{fp}%
%
\special{pn 13}%
\special{pa 1642 1010}%
\special{pa 1096 1420}%
\special{fp}%
\special{sh 1}%
\special{pa 1096 1420}%
\special{pa 1160 1396}%
\special{pa 1138 1388}%
\special{pa 1136 1364}%
\special{pa 1096 1420}%
\special{fp}%
%
\special{pn 13}%
\special{pa 1700 1104}%
\special{pa 2340 1104}%
\special{fp}%
%
\special{pn 13}%
\special{pa 2000 50}%
\special{pa 2000 1264}%
\special{fp}%
%
\special{pn 13}%
\special{pa 1712 874}%
\special{pa 2340 866}%
\special{fp}%
%
\special{pn 13}%
\special{pa 1786 706}%
\special{pa 2604 20}%
\special{fp}%
%
\special{pn 13}%
\special{pa 1910 430}%
\special{pa 2210 736}%
\special{fp}%
%
\special{pn 13}%
\special{pa 2242 130}%
\special{pa 2530 418}%
\special{fp}%
%
\special{pn 8}%
\special{pa 2580 378}%
\special{pa 2600 404}%
\special{pa 2616 430}%
\special{pa 2630 460}%
\special{pa 2638 492}%
\special{pa 2644 526}%
\special{pa 2644 562}%
\special{pa 2642 598}%
\special{pa 2636 636}%
\special{pa 2628 674}%
\special{pa 2614 712}%
\special{pa 2600 748}%
\special{pa 2582 784}%
\special{pa 2562 818}%
\special{pa 2542 850}%
\special{pa 2518 880}%
\special{pa 2492 906}%
\special{pa 2464 928}%
\special{pa 2436 948}%
\special{pa 2406 962}%
\special{pa 2376 972}%
\special{pa 2344 976}%
\special{pa 2314 976}%
\special{pa 2282 970}%
\special{sp}%
%
\special{pn 8}%
\special{pa 450 1306}%
\special{pa 426 1326}%
\special{pa 402 1348}%
\special{pa 380 1372}%
\special{pa 358 1396}%
\special{pa 340 1422}%
\special{pa 322 1448}%
\special{pa 308 1478}%
\special{pa 296 1508}%
\special{pa 286 1540}%
\special{pa 278 1572}%
\special{pa 272 1606}%
\special{pa 270 1638}%
\special{pa 270 1672}%
\special{pa 270 1704}%
\special{pa 274 1738}%
\special{pa 280 1768}%
\special{pa 290 1800}%
\special{pa 300 1830}%
\special{pa 312 1858}%
\special{pa 326 1888}%
\special{pa 340 1916}%
\special{pa 356 1944}%
\special{pa 374 1972}%
\special{pa 392 1998}%
\special{pa 410 2026}%
\special{pa 430 2052}%
\special{pa 438 2066}%
\special{sp}%
%
\special{pn 13}%
\special{pa 2436 1018}%
\special{pa 3008 1438}%
\special{fp}%
\special{sh 1}%
\special{pa 3008 1438}%
\special{pa 2966 1382}%
\special{pa 2964 1406}%
\special{pa 2942 1414}%
\special{pa 3008 1438}%
\special{fp}%
%
\special{pn 13}%
\special{pa 3024 2002}%
\special{pa 3928 1264}%
\special{fp}%
%
\special{pn 13}%
\special{pa 3606 1282}%
\special{pa 3894 1570}%
\special{fp}%
%
\special{pn 13}%
\special{pa 3756 1196}%
\special{pa 4044 1484}%
\special{fp}%
\put(12.9100,-13.9700){\makebox(0,0)[lb]{$q_1$}}%
\put(1.5000,-17.3700){\makebox(0,0)[lb]{$T$}}%
\put(11.1800,-22.4300){\makebox(0,0)[lb]{$\Delta_0$}}%
\put(26.6100,-7.9800){\makebox(0,0)[lb]{$\overline{T}$}}%
\put(6.7000,-10.6000){\makebox(0,0)[lb]{${\mathfrak {f}}_o$}}%
\put(19.1000,-14.3000){\makebox(0,0)[lb]{$\overline{{\mathfrak {f}}}_o$}}%
\put(16.5300,-5.3800){\makebox(0,0)[lb]{$E_1$}}%
\put(34.5600,-9.9900){\makebox(0,0)[lb]{$\mcQ$}}%
\put(33.3600,-18.1700){\makebox(0,0)[lb]{$z_o$}}%
\put(29.1500,-20.9900){\makebox(0,0)[lb]{$L_{z_o}$}}%
%
\special{pn 8}%
\special{pa 3256 1288}%
\special{pa 3270 1258}%
\special{pa 3284 1230}%
\special{pa 3300 1202}%
\special{pa 3318 1176}%
\special{pa 3338 1152}%
\special{pa 3362 1130}%
\special{pa 3386 1110}%
\special{pa 3414 1092}%
\special{pa 3444 1076}%
\special{pa 3474 1064}%
\special{pa 3506 1054}%
\special{pa 3538 1046}%
\special{pa 3572 1040}%
\special{pa 3604 1036}%
\special{pa 3638 1034}%
\special{pa 3670 1036}%
\special{pa 3702 1038}%
\special{pa 3734 1042}%
\special{pa 3766 1048}%
\special{pa 3796 1054}%
\special{pa 3828 1062}%
\special{pa 3854 1068}%
\special{sp}%
\put(21.2000,-12.2400){\makebox(0,0)[lb]{$\Delta_0$}}%
\put(28.8000,-12.0000){\makebox(0,0)[lb]{$\overline{q}_1$}}%
%
\special{pn 13}%
\special{pa 560 1350}%
\special{pa 584 1372}%
\special{pa 606 1396}%
\special{pa 628 1418}%
\special{pa 648 1444}%
\special{pa 668 1470}%
\special{pa 684 1498}%
\special{pa 700 1528}%
\special{pa 712 1558}%
\special{pa 720 1592}%
\special{pa 726 1624}%
\special{pa 726 1656}%
\special{pa 724 1688}%
\special{pa 716 1716}%
\special{pa 702 1744}%
\special{pa 686 1770}%
\special{pa 666 1796}%
\special{pa 642 1820}%
\special{pa 618 1844}%
\special{pa 594 1868}%
\special{pa 590 1870}%
\special{sp}%
%
\special{pn 13}%
\special{pa 590 1620}%
\special{pa 930 1610}%
\special{fp}%
%
\special{pn 13}%
\special{pa 3450 1660}%
\special{pa 3424 1680}%
\special{pa 3396 1698}%
\special{pa 3370 1718}%
\special{pa 3344 1738}%
\special{pa 3320 1758}%
\special{pa 3298 1780}%
\special{pa 3278 1804}%
\special{pa 3260 1828}%
\special{pa 3244 1854}%
\special{pa 3230 1882}%
\special{pa 3216 1910}%
\special{pa 3206 1940}%
\special{pa 3196 1970}%
\special{pa 3186 2002}%
\special{pa 3178 2034}%
\special{pa 3172 2066}%
\special{pa 3166 2100}%
\special{pa 3160 2132}%
\special{pa 3154 2166}%
\special{pa 3150 2190}%
\special{sp}%
%
\special{pn 13}%
\special{pa 3520 1340}%
\special{pa 3524 1374}%
\special{pa 3526 1410}%
\special{pa 3526 1442}%
\special{pa 3526 1476}%
\special{pa 3524 1508}%
\special{pa 3516 1538}%
\special{pa 3508 1566}%
\special{pa 3494 1592}%
\special{pa 3476 1618}%
\special{pa 3454 1640}%
\special{pa 3430 1662}%
\special{pa 3404 1684}%
\special{pa 3378 1704}%
\special{pa 3370 1710}%
\special{sp}%
\end{picture}%


\medskip

The case $(b)$
\end{center}

Note that  if ${\mathfrak f}_o$ meets both of the local branches of $T$ at
$x_o$ transversely, we have the case $(a)$, while if
${\mathfrak f}_o$ is tangent to one of the local branches of $T$ at $x_o$,
we have the case $(b)$.


Blow down  $\overline{\mathfrak f}_o$ and $\Delta_0$ in this order. Then 
the resulting surface is ${\mathbb P}^2$. We denote this composition of 
blowing downs by $\overline{q}_1 : \Sigma_2^{(1)} \to \PP^2$ and put
 $\mcQ:= \overline{q}_1(T)$. Then $\mcQ$ is a reduced quartic 
with the distinguished point $z_o := \overline{q}_1({\mathfrak f}_o\cup \Delta_0)$. 
Note that $\overline{q}_1(E_1)$ is the tangent line $L_{z_o}$ of $\mcQ$ at $z_o$.
Put $\overline{q} := \overline{q}_1\circ q_2\circ \cdots \circ q_r$ and 
we have the following diagram:

\[
\begin{CD}
S'' @<{\overline{\mu}}<< S \\
@V{f''}VV                 @VV{f}V \\
\PP^2  @<<{\overline{q}}< \widehat{\Sigma}_2.
\end{CD}
\]
Here $S''$ is the Stein factorization of $\overline{q}\circ f$. Note that $S''$ is a double 
cover with branch locus $\mcQ$ and that
the pencil of lines through $z_o$ gives rise to the elliptic fibration of 
$S$. Now we have the following proposition.

\begin{prop}\label{prop:res-1}{Let $s_1, s_2$ and $s_3$ be sections as in the beginning
of this section and put $C_i := \overline{q}(s_i) (i = 1, 2, 3)$. If $(\mcQ +  C_1 +  C_3)$ and
 $(\mcQ + C_2 +  C_3)$
have the same combinatorics, then 
$(\mcQ +  C_1 +  C_3, \mcQ +  C_2 + C_3)$ is a Zariski pair.
}
\end{prop}

 \proof Our statement is immediate from Proposition~\ref{prop:fund},  Theorem~\ref{thm:main1} and the following lemma.
 \proofend
 
\begin{lem}\label{lem:equiv}{Let $p$ be an odd prime. For $i = 1, 2$, there exists
a $D_{2p}$-cover 
$\varpi_p: {\mathcal X}_p \to \PP^2$ of $\PP^2$ branched at 
$2\mcQ+ p(C_i + C_3)$ if and only if 
there exists an elliptic $D_{2p}$-cover $\pi_p : X_p \to \widehat{\Sigma}_2$ of 
$\widehat{\Sigma}_2$ such that the horizontal part of $\Delta_{\beta_2(\pi_p)}$ is
$s_i + s_3 + \sigma_f^*(s_i + s_3)$.
}
\end{lem}

\proof  Suppose that there exists a $D_{2p}$-cover $\varpi_p : {\mathcal X}_p \to \PP^2$ 
branched at $2\mcQ + p(C_i + C_3)$. Let $\varpi_p^{(i)} : {\mathcal X}_p^{(i)} \to 
\Sigma_2^{(i)}$ be the induced $D_{2p}$-cover, i.e., ${\mathcal X}_p^{(i)}$ is 
the $\CC({\mathcal X}_p)$-normalization of $\Sigma_2^{(i)}$. Since
$D({\mathcal X}_p/\PP^2) = S''$ and $\beta_1(\varpi_p) = f''$, 
 $D({\mathcal X}_p^{(1)}/\Sigma_2^{(1)})$
  is the $\CC(S'')$-normalization of $\Sigma_2^{(1)}$.
 Hence $\Delta_{\beta_1(\varpi_p^{(1)})} = \Delta_0 + \overline{T}$ as 
 $\overline{q}_1^*\mcQ = \Delta_0 + \overline{T} + 2\overline{\mathfrak f}_o$.  
 This implies that $D({\mathcal X}_p^{(r)}/\widehat{\Sigma}_2) = S$ and
 $\beta_1(\varpi_p^{(r)}) = f$. As  $C_i = \overline{q}\circ f(s_i) (i = 1, 2, 3)$,
 $\varpi_p^{(r)} : {\mathcal X}_p^{(r)} \to \widehat{\Sigma}_2$ is an elliptic
 $D_{2p}$-cover such that the horizontal par of $\Delta_{\beta_2(\varpi_p^{(r)})}$ is
 $s_i + s_3 + \sigma_f^*(s_i + s_3)$.
 
 Conversely, suppose that there exists an elliptic $D_{2p}$-cover $\pi_p : X_p \to \widehat{\Sigma}_2$  such that the horizontal part of $\Delta_{\beta_2(\pi_p)}$ is
$s_i + s_3 + \sigma_f^*(s_i + s_3)$. 
Since $E_1$ gives rise to an irreducible component
``$\Theta_1$"  of singular fiber of type $\I_2$ or $\II$, the preimage of 
$E_1$ in $\widehat{\Sigma}_2$ is not contained
in the branch locus of $\pi_p$ by Corollary~\ref{cor:split} and Remark~\ref{rem:double-cover-construction}.
  Now let $\overline{X}_p$ be the
Stein factorization of $\overline{q}\circ \pi_p$. Then the induced
$D_{2p}$-cover $\overline{\pi}_p : \overline{X}_p \to \PP^2$ is branched at
$2\mcQ + p(C_i + C_3)$.  \proofend


\section{Proof of Theorem~\ref{thm:zpair}}

{\sl Proof  of Theorem~\ref{thm:zpair} (i).}
@Put $\mcQ = C_1 + L_1 + L_2$ and choose 
a point $z_o \in C_1\cap C_2$ as the distinguished point. Let $f''_{\mcQ} : S''_{\mcQ} \to \PP^2$ 
be a double cover
with branch locus $\mcQ$ and let $\varphi_{z_o} : S_{(\mcQ, z_o)} \to \PP^1$ be
the rational elliptic surface as in \S 4. By our construction of $S_{\mcQ,z_o}$,  both $L_3$ and 
$L_4$ give rise 
to sections, which we denote by $s_{L_i}^+$ and 
$s_{L_i}^- (= \sigma_f^*s_{L_i}^+ = [-1]s_{L_i})$ $(i = 3, 4)$, respectively.
By labeling singular fibers suitably, we may assume that $s_{L_i}^+$ $(i = 3, 4)$ and reducible 
singular fibers 
meet as in the following picture:

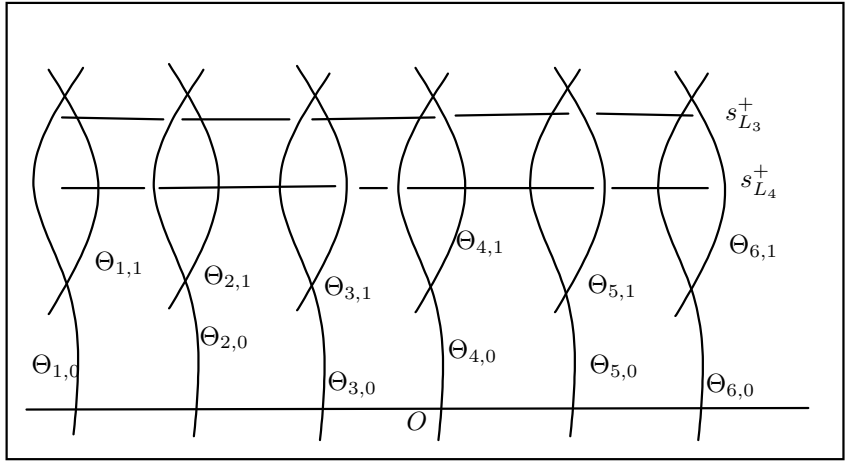
\begin{center}

\unitlength 0.1in
\begin{picture}( 43.8000, 23.9000)(  9.8000,-28.1000)
%
\special{pn 13}%
\special{pa 980 420}%
\special{pa 5360 420}%
\special{pa 5360 2810}%
\special{pa 980 2810}%
\special{pa 980 420}%
\special{fp}%
%
\special{pn 13}%
\special{pa 1084 2550}%
\special{pa 5176 2542}%
\special{fp}%
\put(47.4000,-10.9000){\makebox(0,0)[lb]{$s_{L_3}^+$}}%
\put(48.2000,-14.4000){\makebox(0,0)[lb]{$s_{L_4}^+$}}%
\put(30.7400,-26.5700){\makebox(0,0)[lb]{$O$}}%
\put(11.1000,-24.0000){\makebox(0,0)[lb]{$\Theta_{1,0}$}}%
\put(33.3000,-17.5000){\makebox(0,0)[lb]{$\Theta_{4,1}$}}%
\put(14.4000,-18.6000){\makebox(0,0)[lb]{$\Theta_{1,1}$}}%
\put(32.9000,-23.2000){\makebox(0,0)[lb]{$\Theta_{4,0}$}}%
%
\special{pn 13}%
\special{pa 2010 770}%
\special{pa 1992 800}%
\special{pa 1972 830}%
\special{pa 1952 858}%
\special{pa 1934 888}%
\special{pa 1916 916}%
\special{pa 1898 946}%
\special{pa 1880 976}%
\special{pa 1862 1004}%
\special{pa 1846 1034}%
\special{pa 1832 1064}%
\special{pa 1818 1094}%
\special{pa 1804 1122}%
\special{pa 1792 1152}%
\special{pa 1782 1182}%
\special{pa 1772 1210}%
\special{pa 1764 1240}%
\special{pa 1758 1270}%
\special{pa 1754 1300}%
\special{pa 1752 1330}%
\special{pa 1750 1360}%
\special{pa 1750 1388}%
\special{pa 1754 1418}%
\special{pa 1758 1448}%
\special{pa 1764 1478}%
\special{pa 1772 1508}%
\special{pa 1782 1538}%
\special{pa 1792 1568}%
\special{pa 1802 1598}%
\special{pa 1814 1628}%
\special{pa 1826 1660}%
\special{pa 1840 1690}%
\special{pa 1852 1720}%
\special{pa 1866 1750}%
\special{pa 1878 1782}%
\special{pa 1890 1812}%
\special{pa 1904 1842}%
\special{pa 1914 1874}%
\special{pa 1926 1904}%
\special{pa 1936 1936}%
\special{pa 1944 1968}%
\special{pa 1952 1998}%
\special{pa 1958 2030}%
\special{pa 1964 2062}%
\special{pa 1970 2092}%
\special{pa 1974 2124}%
\special{pa 1978 2156}%
\special{pa 1980 2188}%
\special{pa 1982 2220}%
\special{pa 1984 2252}%
\special{pa 1984 2284}%
\special{pa 1984 2316}%
\special{pa 1984 2348}%
\special{pa 1982 2380}%
\special{pa 1982 2412}%
\special{pa 1980 2444}%
\special{pa 1978 2478}%
\special{pa 1976 2510}%
\special{pa 1974 2542}%
\special{pa 1972 2574}%
\special{pa 1968 2606}%
\special{pa 1966 2638}%
\special{pa 1962 2672}%
\special{pa 1960 2690}%
\special{sp}%
%
\special{pn 13}%
\special{pa 1200 770}%
\special{pa 1218 800}%
\special{pa 1238 830}%
\special{pa 1256 858}%
\special{pa 1272 888}%
\special{pa 1290 918}%
\special{pa 1308 946}%
\special{pa 1324 976}%
\special{pa 1340 1006}%
\special{pa 1356 1034}%
\special{pa 1370 1064}%
\special{pa 1384 1094}%
\special{pa 1398 1122}%
\special{pa 1410 1152}%
\special{pa 1420 1182}%
\special{pa 1430 1212}%
\special{pa 1438 1240}%
\special{pa 1446 1270}%
\special{pa 1452 1300}%
\special{pa 1456 1330}%
\special{pa 1460 1360}%
\special{pa 1460 1388}%
\special{pa 1460 1418}%
\special{pa 1458 1448}%
\special{pa 1456 1478}%
\special{pa 1450 1508}%
\special{pa 1444 1538}%
\special{pa 1436 1566}%
\special{pa 1428 1596}%
\special{pa 1418 1626}%
\special{pa 1408 1656}%
\special{pa 1396 1686}%
\special{pa 1382 1716}%
\special{pa 1370 1746}%
\special{pa 1354 1776}%
\special{pa 1340 1806}%
\special{pa 1324 1836}%
\special{pa 1308 1866}%
\special{pa 1292 1896}%
\special{pa 1274 1926}%
\special{pa 1256 1956}%
\special{pa 1240 1986}%
\special{pa 1222 2016}%
\special{pa 1204 2046}%
\special{pa 1200 2050}%
\special{sp}%
%
\special{pn 13}%
\special{pa 1380 760}%
\special{pa 1362 790}%
\special{pa 1342 820}%
\special{pa 1322 848}%
\special{pa 1304 878}%
\special{pa 1286 906}%
\special{pa 1268 936}%
\special{pa 1250 966}%
\special{pa 1232 994}%
\special{pa 1216 1024}%
\special{pa 1202 1054}%
\special{pa 1188 1084}%
\special{pa 1174 1112}%
\special{pa 1162 1142}%
\special{pa 1152 1172}%
\special{pa 1142 1200}%
\special{pa 1134 1230}%
\special{pa 1128 1260}%
\special{pa 1124 1290}%
\special{pa 1122 1320}%
\special{pa 1120 1350}%
\special{pa 1120 1378}%
\special{pa 1124 1408}%
\special{pa 1128 1438}%
\special{pa 1134 1468}%
\special{pa 1142 1498}%
\special{pa 1152 1528}%
\special{pa 1162 1558}%
\special{pa 1172 1588}%
\special{pa 1184 1618}%
\special{pa 1196 1650}%
\special{pa 1210 1680}%
\special{pa 1222 1710}%
\special{pa 1236 1740}%
\special{pa 1248 1772}%
\special{pa 1260 1802}%
\special{pa 1274 1832}%
\special{pa 1284 1864}%
\special{pa 1296 1894}%
\special{pa 1306 1926}%
\special{pa 1314 1958}%
\special{pa 1322 1988}%
\special{pa 1328 2020}%
\special{pa 1334 2052}%
\special{pa 1340 2082}%
\special{pa 1344 2114}%
\special{pa 1348 2146}%
\special{pa 1350 2178}%
\special{pa 1352 2210}%
\special{pa 1354 2242}%
\special{pa 1354 2274}%
\special{pa 1354 2306}%
\special{pa 1354 2338}%
\special{pa 1352 2370}%
\special{pa 1352 2402}%
\special{pa 1350 2434}%
\special{pa 1348 2468}%
\special{pa 1346 2500}%
\special{pa 1344 2532}%
\special{pa 1342 2564}%
\special{pa 1338 2596}%
\special{pa 1336 2628}%
\special{pa 1332 2662}%
\special{pa 1330 2680}%
\special{sp}%
%
\special{pn 13}%
\special{pa 1830 740}%
\special{pa 1848 770}%
\special{pa 1868 800}%
\special{pa 1886 828}%
\special{pa 1902 858}%
\special{pa 1920 888}%
\special{pa 1938 916}%
\special{pa 1954 946}%
\special{pa 1970 976}%
\special{pa 1986 1004}%
\special{pa 2000 1034}%
\special{pa 2014 1064}%
\special{pa 2028 1092}%
\special{pa 2040 1122}%
\special{pa 2050 1152}%
\special{pa 2060 1182}%
\special{pa 2068 1210}%
\special{pa 2076 1240}%
\special{pa 2082 1270}%
\special{pa 2086 1300}%
\special{pa 2090 1330}%
\special{pa 2090 1358}%
\special{pa 2090 1388}%
\special{pa 2088 1418}%
\special{pa 2086 1448}%
\special{pa 2080 1478}%
\special{pa 2074 1508}%
\special{pa 2066 1536}%
\special{pa 2058 1566}%
\special{pa 2048 1596}%
\special{pa 2038 1626}%
\special{pa 2026 1656}%
\special{pa 2012 1686}%
\special{pa 2000 1716}%
\special{pa 1984 1746}%
\special{pa 1970 1776}%
\special{pa 1954 1806}%
\special{pa 1938 1836}%
\special{pa 1922 1866}%
\special{pa 1904 1896}%
\special{pa 1886 1926}%
\special{pa 1870 1956}%
\special{pa 1852 1986}%
\special{pa 1834 2016}%
\special{pa 1830 2020}%
\special{sp}%
%
\special{pn 13}%
\special{pa 2670 790}%
\special{pa 2652 820}%
\special{pa 2632 850}%
\special{pa 2612 878}%
\special{pa 2594 908}%
\special{pa 2576 936}%
\special{pa 2558 966}%
\special{pa 2540 996}%
\special{pa 2522 1024}%
\special{pa 2506 1054}%
\special{pa 2492 1084}%
\special{pa 2478 1114}%
\special{pa 2464 1142}%
\special{pa 2452 1172}%
\special{pa 2442 1202}%
\special{pa 2432 1230}%
\special{pa 2424 1260}%
\special{pa 2418 1290}%
\special{pa 2414 1320}%
\special{pa 2412 1350}%
\special{pa 2410 1380}%
\special{pa 2410 1408}%
\special{pa 2414 1438}%
\special{pa 2418 1468}%
\special{pa 2424 1498}%
\special{pa 2432 1528}%
\special{pa 2442 1558}%
\special{pa 2452 1588}%
\special{pa 2462 1618}%
\special{pa 2474 1648}%
\special{pa 2486 1680}%
\special{pa 2500 1710}%
\special{pa 2512 1740}%
\special{pa 2526 1770}%
\special{pa 2538 1802}%
\special{pa 2550 1832}%
\special{pa 2564 1862}%
\special{pa 2574 1894}%
\special{pa 2586 1924}%
\special{pa 2596 1956}%
\special{pa 2604 1988}%
\special{pa 2612 2018}%
\special{pa 2618 2050}%
\special{pa 2624 2082}%
\special{pa 2630 2112}%
\special{pa 2634 2144}%
\special{pa 2638 2176}%
\special{pa 2640 2208}%
\special{pa 2642 2240}%
\special{pa 2644 2272}%
\special{pa 2644 2304}%
\special{pa 2644 2336}%
\special{pa 2644 2368}%
\special{pa 2642 2400}%
\special{pa 2642 2432}%
\special{pa 2640 2464}%
\special{pa 2638 2498}%
\special{pa 2636 2530}%
\special{pa 2634 2562}%
\special{pa 2632 2594}%
\special{pa 2628 2626}%
\special{pa 2626 2658}%
\special{pa 2622 2692}%
\special{pa 2620 2710}%
\special{sp}%
%
\special{pn 13}%
\special{pa 2500 750}%
\special{pa 2518 780}%
\special{pa 2538 810}%
\special{pa 2556 838}%
\special{pa 2572 868}%
\special{pa 2590 898}%
\special{pa 2608 926}%
\special{pa 2624 956}%
\special{pa 2640 986}%
\special{pa 2656 1014}%
\special{pa 2670 1044}%
\special{pa 2684 1074}%
\special{pa 2698 1102}%
\special{pa 2710 1132}%
\special{pa 2720 1162}%
\special{pa 2730 1192}%
\special{pa 2738 1220}%
\special{pa 2746 1250}%
\special{pa 2752 1280}%
\special{pa 2756 1310}%
\special{pa 2760 1340}%
\special{pa 2760 1368}%
\special{pa 2760 1398}%
\special{pa 2758 1428}%
\special{pa 2756 1458}%
\special{pa 2750 1488}%
\special{pa 2744 1518}%
\special{pa 2736 1546}%
\special{pa 2728 1576}%
\special{pa 2718 1606}%
\special{pa 2708 1636}%
\special{pa 2696 1666}%
\special{pa 2682 1696}%
\special{pa 2670 1726}%
\special{pa 2654 1756}%
\special{pa 2640 1786}%
\special{pa 2624 1816}%
\special{pa 2608 1846}%
\special{pa 2592 1876}%
\special{pa 2574 1906}%
\special{pa 2556 1936}%
\special{pa 2540 1966}%
\special{pa 2522 1996}%
\special{pa 2504 2026}%
\special{pa 2500 2030}%
\special{sp}%
%
\special{pn 13}%
\special{pa 3290 790}%
\special{pa 3272 820}%
\special{pa 3252 850}%
\special{pa 3232 878}%
\special{pa 3214 908}%
\special{pa 3196 936}%
\special{pa 3178 966}%
\special{pa 3160 996}%
\special{pa 3142 1024}%
\special{pa 3126 1054}%
\special{pa 3112 1084}%
\special{pa 3098 1114}%
\special{pa 3084 1142}%
\special{pa 3072 1172}%
\special{pa 3062 1202}%
\special{pa 3052 1230}%
\special{pa 3044 1260}%
\special{pa 3038 1290}%
\special{pa 3034 1320}%
\special{pa 3032 1350}%
\special{pa 3030 1380}%
\special{pa 3030 1408}%
\special{pa 3034 1438}%
\special{pa 3038 1468}%
\special{pa 3044 1498}%
\special{pa 3052 1528}%
\special{pa 3062 1558}%
\special{pa 3072 1588}%
\special{pa 3082 1618}%
\special{pa 3094 1648}%
\special{pa 3106 1680}%
\special{pa 3120 1710}%
\special{pa 3132 1740}%
\special{pa 3146 1770}%
\special{pa 3158 1802}%
\special{pa 3170 1832}%
\special{pa 3184 1862}%
\special{pa 3194 1894}%
\special{pa 3206 1924}%
\special{pa 3216 1956}%
\special{pa 3224 1988}%
\special{pa 3232 2018}%
\special{pa 3238 2050}%
\special{pa 3244 2082}%
\special{pa 3250 2112}%
\special{pa 3254 2144}%
\special{pa 3258 2176}%
\special{pa 3260 2208}%
\special{pa 3262 2240}%
\special{pa 3264 2272}%
\special{pa 3264 2304}%
\special{pa 3264 2336}%
\special{pa 3264 2368}%
\special{pa 3262 2400}%
\special{pa 3262 2432}%
\special{pa 3260 2464}%
\special{pa 3258 2498}%
\special{pa 3256 2530}%
\special{pa 3254 2562}%
\special{pa 3252 2594}%
\special{pa 3248 2626}%
\special{pa 3246 2658}%
\special{pa 3242 2692}%
\special{pa 3240 2710}%
\special{sp}%
%
\special{pn 13}%
\special{pa 3120 760}%
\special{pa 3138 790}%
\special{pa 3158 820}%
\special{pa 3176 848}%
\special{pa 3192 878}%
\special{pa 3210 908}%
\special{pa 3228 936}%
\special{pa 3244 966}%
\special{pa 3260 996}%
\special{pa 3276 1024}%
\special{pa 3290 1054}%
\special{pa 3304 1084}%
\special{pa 3318 1112}%
\special{pa 3330 1142}%
\special{pa 3340 1172}%
\special{pa 3350 1202}%
\special{pa 3358 1230}%
\special{pa 3366 1260}%
\special{pa 3372 1290}%
\special{pa 3376 1320}%
\special{pa 3380 1350}%
\special{pa 3380 1378}%
\special{pa 3380 1408}%
\special{pa 3378 1438}%
\special{pa 3376 1468}%
\special{pa 3370 1498}%
\special{pa 3364 1528}%
\special{pa 3356 1556}%
\special{pa 3348 1586}%
\special{pa 3338 1616}%
\special{pa 3328 1646}%
\special{pa 3316 1676}%
\special{pa 3302 1706}%
\special{pa 3290 1736}%
\special{pa 3274 1766}%
\special{pa 3260 1796}%
\special{pa 3244 1826}%
\special{pa 3228 1856}%
\special{pa 3212 1886}%
\special{pa 3194 1916}%
\special{pa 3176 1946}%
\special{pa 3160 1976}%
\special{pa 3142 2006}%
\special{pa 3124 2036}%
\special{pa 3120 2040}%
\special{sp}%
\put(20.1000,-19.3000){\makebox(0,0)[lb]{$\Theta_{2,1}$}}%
\put(19.9000,-22.5000){\makebox(0,0)[lb]{$\Theta_{2,0}$}}%
\put(26.6000,-24.9000){\makebox(0,0)[lb]{$\Theta_{3,0}$}}%
\put(26.4000,-19.9000){\makebox(0,0)[lb]{$\Theta_{3,1}$}}%
%
\special{pn 13}%
\special{pa 3850 760}%
\special{pa 3868 790}%
\special{pa 3888 820}%
\special{pa 3906 848}%
\special{pa 3922 878}%
\special{pa 3940 908}%
\special{pa 3958 936}%
\special{pa 3974 966}%
\special{pa 3990 996}%
\special{pa 4006 1024}%
\special{pa 4020 1054}%
\special{pa 4034 1084}%
\special{pa 4048 1112}%
\special{pa 4060 1142}%
\special{pa 4070 1172}%
\special{pa 4080 1202}%
\special{pa 4088 1230}%
\special{pa 4096 1260}%
\special{pa 4102 1290}%
\special{pa 4106 1320}%
\special{pa 4110 1350}%
\special{pa 4110 1378}%
\special{pa 4110 1408}%
\special{pa 4108 1438}%
\special{pa 4106 1468}%
\special{pa 4100 1498}%
\special{pa 4094 1528}%
\special{pa 4086 1556}%
\special{pa 4078 1586}%
\special{pa 4068 1616}%
\special{pa 4058 1646}%
\special{pa 4046 1676}%
\special{pa 4032 1706}%
\special{pa 4020 1736}%
\special{pa 4004 1766}%
\special{pa 3990 1796}%
\special{pa 3974 1826}%
\special{pa 3958 1856}%
\special{pa 3942 1886}%
\special{pa 3924 1916}%
\special{pa 3906 1946}%
\special{pa 3890 1976}%
\special{pa 3872 2006}%
\special{pa 3854 2036}%
\special{pa 3850 2040}%
\special{sp}%
%
\special{pn 13}%
\special{pa 3980 770}%
\special{pa 3962 800}%
\special{pa 3942 830}%
\special{pa 3922 858}%
\special{pa 3904 888}%
\special{pa 3886 916}%
\special{pa 3868 946}%
\special{pa 3850 976}%
\special{pa 3832 1004}%
\special{pa 3816 1034}%
\special{pa 3802 1064}%
\special{pa 3788 1094}%
\special{pa 3774 1122}%
\special{pa 3762 1152}%
\special{pa 3752 1182}%
\special{pa 3742 1210}%
\special{pa 3734 1240}%
\special{pa 3728 1270}%
\special{pa 3724 1300}%
\special{pa 3722 1330}%
\special{pa 3720 1360}%
\special{pa 3720 1388}%
\special{pa 3724 1418}%
\special{pa 3728 1448}%
\special{pa 3734 1478}%
\special{pa 3742 1508}%
\special{pa 3752 1538}%
\special{pa 3762 1568}%
\special{pa 3772 1598}%
\special{pa 3784 1628}%
\special{pa 3796 1660}%
\special{pa 3810 1690}%
\special{pa 3822 1720}%
\special{pa 3836 1750}%
\special{pa 3848 1782}%
\special{pa 3860 1812}%
\special{pa 3874 1842}%
\special{pa 3884 1874}%
\special{pa 3896 1904}%
\special{pa 3906 1936}%
\special{pa 3914 1968}%
\special{pa 3922 1998}%
\special{pa 3928 2030}%
\special{pa 3934 2062}%
\special{pa 3940 2092}%
\special{pa 3944 2124}%
\special{pa 3948 2156}%
\special{pa 3950 2188}%
\special{pa 3952 2220}%
\special{pa 3954 2252}%
\special{pa 3954 2284}%
\special{pa 3954 2316}%
\special{pa 3954 2348}%
\special{pa 3952 2380}%
\special{pa 3952 2412}%
\special{pa 3950 2444}%
\special{pa 3948 2478}%
\special{pa 3946 2510}%
\special{pa 3944 2542}%
\special{pa 3942 2574}%
\special{pa 3938 2606}%
\special{pa 3936 2638}%
\special{pa 3932 2672}%
\special{pa 3930 2690}%
\special{sp}%
\put(40.2000,-19.8000){\makebox(0,0)[lb]{$\Theta_{5,1}$}}%
\put(40.3000,-24.0000){\makebox(0,0)[lb]{$\Theta_{5,0}$}}%
%
\special{pn 13}%
\special{pa 4480 780}%
\special{pa 4498 810}%
\special{pa 4518 840}%
\special{pa 4536 868}%
\special{pa 4552 898}%
\special{pa 4570 928}%
\special{pa 4588 956}%
\special{pa 4604 986}%
\special{pa 4620 1016}%
\special{pa 4636 1044}%
\special{pa 4650 1074}%
\special{pa 4664 1104}%
\special{pa 4678 1132}%
\special{pa 4690 1162}%
\special{pa 4700 1192}%
\special{pa 4710 1222}%
\special{pa 4718 1250}%
\special{pa 4726 1280}%
\special{pa 4732 1310}%
\special{pa 4736 1340}%
\special{pa 4740 1370}%
\special{pa 4740 1398}%
\special{pa 4740 1428}%
\special{pa 4738 1458}%
\special{pa 4736 1488}%
\special{pa 4730 1518}%
\special{pa 4724 1548}%
\special{pa 4716 1576}%
\special{pa 4708 1606}%
\special{pa 4698 1636}%
\special{pa 4688 1666}%
\special{pa 4676 1696}%
\special{pa 4662 1726}%
\special{pa 4650 1756}%
\special{pa 4634 1786}%
\special{pa 4620 1816}%
\special{pa 4604 1846}%
\special{pa 4588 1876}%
\special{pa 4572 1906}%
\special{pa 4554 1936}%
\special{pa 4536 1966}%
\special{pa 4520 1996}%
\special{pa 4502 2026}%
\special{pa 4484 2056}%
\special{pa 4480 2060}%
\special{sp}%
%
\special{pn 13}%
\special{pa 4650 790}%
\special{pa 4632 820}%
\special{pa 4612 850}%
\special{pa 4592 878}%
\special{pa 4574 908}%
\special{pa 4556 936}%
\special{pa 4538 966}%
\special{pa 4520 996}%
\special{pa 4502 1024}%
\special{pa 4486 1054}%
\special{pa 4472 1084}%
\special{pa 4458 1114}%
\special{pa 4444 1142}%
\special{pa 4432 1172}%
\special{pa 4422 1202}%
\special{pa 4412 1230}%
\special{pa 4404 1260}%
\special{pa 4398 1290}%
\special{pa 4394 1320}%
\special{pa 4392 1350}%
\special{pa 4390 1380}%
\special{pa 4390 1408}%
\special{pa 4394 1438}%
\special{pa 4398 1468}%
\special{pa 4404 1498}%
\special{pa 4412 1528}%
\special{pa 4422 1558}%
\special{pa 4432 1588}%
\special{pa 4442 1618}%
\special{pa 4454 1648}%
\special{pa 4466 1680}%
\special{pa 4480 1710}%
\special{pa 4492 1740}%
\special{pa 4506 1770}%
\special{pa 4518 1802}%
\special{pa 4530 1832}%
\special{pa 4544 1862}%
\special{pa 4554 1894}%
\special{pa 4566 1924}%
\special{pa 4576 1956}%
\special{pa 4584 1988}%
\special{pa 4592 2018}%
\special{pa 4598 2050}%
\special{pa 4604 2082}%
\special{pa 4610 2112}%
\special{pa 4614 2144}%
\special{pa 4618 2176}%
\special{pa 4620 2208}%
\special{pa 4622 2240}%
\special{pa 4624 2272}%
\special{pa 4624 2304}%
\special{pa 4624 2336}%
\special{pa 4624 2368}%
\special{pa 4622 2400}%
\special{pa 4622 2432}%
\special{pa 4620 2464}%
\special{pa 4618 2498}%
\special{pa 4616 2530}%
\special{pa 4614 2562}%
\special{pa 4612 2594}%
\special{pa 4608 2626}%
\special{pa 4606 2658}%
\special{pa 4602 2692}%
\special{pa 4600 2710}%
\special{sp}%
%
\special{pn 13}%
\special{pa 3330 1010}%
\special{pa 3920 1000}%
\special{fp}%
%
\special{pn 13}%
\special{pa 4070 1000}%
\special{pa 4570 1010}%
\special{fp}%
%
\special{pn 13}%
\special{pa 1270 1020}%
\special{pa 1800 1030}%
\special{fp}%
%
\special{pn 13}%
\special{pa 1900 1030}%
\special{pa 2460 1030}%
\special{fp}%
%
\special{pn 13}%
\special{pa 2580 1030}%
\special{pa 3220 1020}%
\special{fp}%
%
\special{pn 13}%
\special{pa 1270 1390}%
\special{pa 1700 1390}%
\special{fp}%
%
\special{pn 13}%
\special{pa 1780 1390}%
\special{pa 2700 1380}%
\special{fp}%
%
\special{pn 13}%
\special{pa 2830 1390}%
\special{pa 2970 1390}%
\special{fp}%
%
\special{pn 13}%
\special{pa 3080 1390}%
\special{pa 4050 1390}%
\special{fp}%
%
\special{pn 13}%
\special{pa 4150 1390}%
\special{pa 4650 1390}%
\special{fp}%
\put(47.6000,-17.7000){\makebox(0,0)[lb]{$\Theta_{6,1}$}}%
\put(46.4000,-25.0000){\makebox(0,0)[lb]{$\Theta_{6,0}$}}%
\end{picture}%

\end{center}

Here we assume that $\Theta_{1, 0}$ and $O$ come from $z_o$. 
By the explicit formula of $\langle\, , \, \rangle$,  we have
\[
\langle s_{L_i}^{\pm}, s_{L_i}^{\pm} \rangle  = \frac 12,  (i = 3, 4) \quad \langle s_{L_3}^+, s_{L_4}^+\rangle = 0.
\]
By \cite{oguiso-shioda}, 
$\MW(S_{(\mcQ, z_o)}) \cong (A_1^*)^{\oplus 2} \oplus (\ZZ/2\ZZ)^{\oplus 2}$ and
 we may assume that
\[
(A_1^*)^{\oplus 2} \cong \ZZ s_{L_3}^+ \oplus \ZZ s_{L_4}^+,
\]
and that the $2$-torsions sections arise from $C_1, L_1$ and $L_2$.

As for $(\overline{q}\circ f)^*(C_2)$,  it also gives rise to two sections $s_{C_2}^{\pm}$. Since $C_2$ does not pass through
any singularities of $C_1 + L_1 + L_2$ and $s_{C_2}^{\pm}O = 0$, we have $\langle s_{C_2}^{\pm}, s_{C_2}^{\pm} \rangle = 2$.

On the other hand, any element $s \in \MW(S_{(\mcQ,z_o)})$ with  $\langle s, s \rangle = 2$ is of the form
\[
[2]s_{L_i}^{\pm} \dot{+} \tau, \quad (i = 3, 4) \quad \tau \in \MW(S_{(\mcQ,z_o)})_{\tor}.
\]
If $\tau \neq 0$, then $s_{C_2}^{\pm}$ meets $\Theta_{i, 1}$ for some $i$ by considering
the addition on singular fibers (see \cite[Theorem~9.1]{kodaira} or 
\cite[\S 1]{tokunaga11}). Hence, by the explicit formula
for $\langle \, , \, \rangle$, we have $s_{C_2}^{\pm}O \neq 0$. On the other hand,
 $s_{C_2}^{\pm}O = 0$ by our construction. Thus  we infer $\tau = 0$ and we may assume
that  $s_{C_2}^+ = [2]s_{L_3}^+$ after relabeling $\pm, L_3$ and $L_4$, if necessary.  Therefore
\[
s_{C_2}^+  \dot{+}  [p-2]s_{L_3}^+ \in [p]\MW(S_{(\mcQ,z_o)})
\]
for any odd prime $p$, while
\[
s_{C_2}^+ \dot{+}  [k]s_{L_4}^+ \not\in [p]\MW(S_{(\mcQ,z_o)})
\]
for any odd prime $p$ and $1 \le k \le p-1$. By Proposition~\ref{prop:res-1}, we infer that
$(C_1 + C_2 + L_1 + L_2 + L_3,  C_1 + C_2 + L_1 + L_2+ L_4)$ is a Zariski pair.  \proofend

\bigskip

{\sl Proof for Theorem~\ref{thm:zpair} (ii).}
Put $\mcQ = C_1 + C_2$ and choose 
a point $z_o \in C_1\cap C_3$ as the distinguished point. Let $f''_{\mcQ} : 
S''_{\mcQ} \to \PP^2$ be a double cover
with branch locus $\mcQ$ and let $\varphi_{z_o} : S_{(\mcQ, z_o)} \to \PP^1$ be
the rational elliptic surface as in \S 4. By our construction of $S_{\mcQ,z_o}$,  $L_0$, $L_1$ 
and $L_2$ give rise 
to sections, which we denote by $s_{L_i}^+$ and $s_{L_i}^- (= \sigma_f^*s_{L_i}^+ = [-1]s_{L_i})$ $(i =0,  1, 2)$, respectively.
By labeling singular fibers suitably, we may assume that $s_{L_i}^+$ $(i =0,  1, 2)$ and reducible singular fibers 
meet as in the following picture:

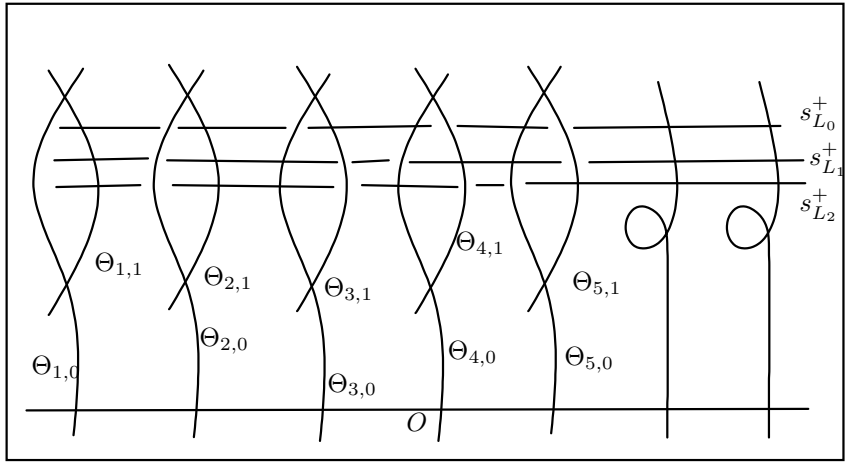
\begin{center}
\unitlength 0.1in
\begin{picture}( 43.8000, 23.9000)(  9.8000,-28.1000)
%
\special{pn 13}%
\special{pa 980 420}%
\special{pa 5360 420}%
\special{pa 5360 2810}%
\special{pa 980 2810}%
\special{pa 980 420}%
\special{fp}%
%
\special{pn 13}%
\special{pa 1084 2550}%
\special{pa 5176 2542}%
\special{fp}%
%
\special{pn 13}%
\special{pa 4390 830}%
\special{pa 4398 862}%
\special{pa 4406 894}%
\special{pa 4414 924}%
\special{pa 4422 956}%
\special{pa 4428 988}%
\special{pa 4436 1018}%
\special{pa 4444 1050}%
\special{pa 4450 1082}%
\special{pa 4456 1112}%
\special{pa 4462 1144}%
\special{pa 4468 1176}%
\special{pa 4472 1206}%
\special{pa 4478 1238}%
\special{pa 4482 1268}%
\special{pa 4486 1300}%
\special{pa 4488 1332}%
\special{pa 4490 1362}%
\special{pa 4490 1394}%
\special{pa 4488 1426}%
\special{pa 4484 1460}%
\special{pa 4480 1494}%
\special{pa 4472 1528}%
\special{pa 4462 1560}%
\special{pa 4450 1594}%
\special{pa 4434 1624}%
\special{pa 4416 1650}%
\special{pa 4392 1674}%
\special{pa 4366 1690}%
\special{pa 4336 1700}%
\special{pa 4302 1702}%
\special{pa 4272 1694}%
\special{pa 4248 1674}%
\special{pa 4230 1644}%
\special{pa 4220 1610}%
\special{pa 4220 1574}%
\special{pa 4228 1540}%
\special{pa 4246 1510}%
\special{pa 4270 1488}%
\special{pa 4300 1480}%
\special{pa 4330 1484}%
\special{pa 4362 1500}%
\special{pa 4390 1522}%
\special{pa 4410 1550}%
\special{pa 4424 1580}%
\special{pa 4432 1610}%
\special{pa 4436 1642}%
\special{pa 4438 1674}%
\special{pa 4438 1708}%
\special{pa 4440 1740}%
\special{pa 4440 1772}%
\special{pa 4440 1806}%
\special{pa 4440 1838}%
\special{pa 4440 1870}%
\special{pa 4440 1902}%
\special{pa 4440 1934}%
\special{pa 4440 1966}%
\special{pa 4440 1998}%
\special{pa 4440 2030}%
\special{pa 4440 2062}%
\special{pa 4440 2094}%
\special{pa 4440 2126}%
\special{pa 4440 2158}%
\special{pa 4440 2190}%
\special{pa 4440 2222}%
\special{pa 4440 2254}%
\special{pa 4440 2286}%
\special{pa 4440 2316}%
\special{pa 4440 2348}%
\special{pa 4440 2380}%
\special{pa 4440 2412}%
\special{pa 4440 2444}%
\special{pa 4440 2476}%
\special{pa 4438 2508}%
\special{pa 4438 2540}%
\special{pa 4438 2572}%
\special{pa 4438 2604}%
\special{pa 4438 2636}%
\special{pa 4438 2668}%
\special{pa 4438 2692}%
\special{sp}%
%
\special{pn 13}%
\special{pa 4920 830}%
\special{pa 4928 862}%
\special{pa 4936 894}%
\special{pa 4944 924}%
\special{pa 4952 956}%
\special{pa 4960 988}%
\special{pa 4966 1018}%
\special{pa 4974 1050}%
\special{pa 4980 1082}%
\special{pa 4986 1112}%
\special{pa 4992 1144}%
\special{pa 4998 1174}%
\special{pa 5004 1206}%
\special{pa 5008 1238}%
\special{pa 5012 1268}%
\special{pa 5016 1300}%
\special{pa 5018 1330}%
\special{pa 5020 1362}%
\special{pa 5020 1394}%
\special{pa 5018 1426}%
\special{pa 5016 1460}%
\special{pa 5010 1492}%
\special{pa 5004 1526}%
\special{pa 4994 1560}%
\special{pa 4980 1594}%
\special{pa 4964 1624}%
\special{pa 4946 1650}%
\special{pa 4924 1672}%
\special{pa 4896 1690}%
\special{pa 4866 1700}%
\special{pa 4834 1702}%
\special{pa 4802 1694}%
\special{pa 4778 1674}%
\special{pa 4760 1646}%
\special{pa 4750 1610}%
\special{pa 4750 1574}%
\special{pa 4758 1540}%
\special{pa 4776 1510}%
\special{pa 4800 1488}%
\special{pa 4828 1480}%
\special{pa 4860 1484}%
\special{pa 4890 1500}%
\special{pa 4918 1522}%
\special{pa 4940 1550}%
\special{pa 4954 1578}%
\special{pa 4964 1610}%
\special{pa 4968 1642}%
\special{pa 4970 1674}%
\special{pa 4970 1706}%
\special{pa 4970 1738}%
\special{pa 4972 1772}%
\special{pa 4972 1804}%
\special{pa 4972 1836}%
\special{pa 4972 1868}%
\special{pa 4972 1900}%
\special{pa 4972 1932}%
\special{pa 4972 1966}%
\special{pa 4972 1998}%
\special{pa 4972 2030}%
\special{pa 4972 2062}%
\special{pa 4972 2094}%
\special{pa 4972 2124}%
\special{pa 4972 2156}%
\special{pa 4972 2188}%
\special{pa 4972 2220}%
\special{pa 4972 2252}%
\special{pa 4972 2284}%
\special{pa 4972 2316}%
\special{pa 4972 2348}%
\special{pa 4970 2380}%
\special{pa 4970 2412}%
\special{pa 4970 2442}%
\special{pa 4970 2474}%
\special{pa 4970 2506}%
\special{pa 4970 2538}%
\special{pa 4970 2570}%
\special{pa 4970 2602}%
\special{pa 4970 2634}%
\special{pa 4970 2666}%
\special{pa 4970 2692}%
\special{sp}%
\put(51.8000,-13.2000){\makebox(0,0)[lb]{$s_{L_1}^+$}}%
\put(51.3000,-15.6000){\makebox(0,0)[lb]{$s_{L_2}^+$}}%
\put(30.7400,-26.5700){\makebox(0,0)[lb]{$O$}}%
\put(11.1000,-24.0000){\makebox(0,0)[lb]{$\Theta_{1,0}$}}%
\put(33.3000,-17.5000){\makebox(0,0)[lb]{$\Theta_{4,1}$}}%
\put(14.4000,-18.6000){\makebox(0,0)[lb]{$\Theta_{1,1}$}}%
\put(32.9000,-23.2000){\makebox(0,0)[lb]{$\Theta_{4,0}$}}%
%
\special{pn 13}%
\special{pa 2010 770}%
\special{pa 1992 800}%
\special{pa 1972 830}%
\special{pa 1952 858}%
\special{pa 1934 888}%
\special{pa 1916 916}%
\special{pa 1898 946}%
\special{pa 1880 976}%
\special{pa 1862 1004}%
\special{pa 1846 1034}%
\special{pa 1832 1064}%
\special{pa 1818 1094}%
\special{pa 1804 1122}%
\special{pa 1792 1152}%
\special{pa 1782 1182}%
\special{pa 1772 1210}%
\special{pa 1764 1240}%
\special{pa 1758 1270}%
\special{pa 1754 1300}%
\special{pa 1752 1330}%
\special{pa 1750 1360}%
\special{pa 1750 1388}%
\special{pa 1754 1418}%
\special{pa 1758 1448}%
\special{pa 1764 1478}%
\special{pa 1772 1508}%
\special{pa 1782 1538}%
\special{pa 1792 1568}%
\special{pa 1802 1598}%
\special{pa 1814 1628}%
\special{pa 1826 1660}%
\special{pa 1840 1690}%
\special{pa 1852 1720}%
\special{pa 1866 1750}%
\special{pa 1878 1782}%
\special{pa 1890 1812}%
\special{pa 1904 1842}%
\special{pa 1914 1874}%
\special{pa 1926 1904}%
\special{pa 1936 1936}%
\special{pa 1944 1968}%
\special{pa 1952 1998}%
\special{pa 1958 2030}%
\special{pa 1964 2062}%
\special{pa 1970 2092}%
\special{pa 1974 2124}%
\special{pa 1978 2156}%
\special{pa 1980 2188}%
\special{pa 1982 2220}%
\special{pa 1984 2252}%
\special{pa 1984 2284}%
\special{pa 1984 2316}%
\special{pa 1984 2348}%
\special{pa 1982 2380}%
\special{pa 1982 2412}%
\special{pa 1980 2444}%
\special{pa 1978 2478}%
\special{pa 1976 2510}%
\special{pa 1974 2542}%
\special{pa 1972 2574}%
\special{pa 1968 2606}%
\special{pa 1966 2638}%
\special{pa 1962 2672}%
\special{pa 1960 2690}%
\special{sp}%
%
\special{pn 13}%
\special{pa 1200 770}%
\special{pa 1218 800}%
\special{pa 1238 830}%
\special{pa 1256 858}%
\special{pa 1272 888}%
\special{pa 1290 918}%
\special{pa 1308 946}%
\special{pa 1324 976}%
\special{pa 1340 1006}%
\special{pa 1356 1034}%
\special{pa 1370 1064}%
\special{pa 1384 1094}%
\special{pa 1398 1122}%
\special{pa 1410 1152}%
\special{pa 1420 1182}%
\special{pa 1430 1212}%
\special{pa 1438 1240}%
\special{pa 1446 1270}%
\special{pa 1452 1300}%
\special{pa 1456 1330}%
\special{pa 1460 1360}%
\special{pa 1460 1388}%
\special{pa 1460 1418}%
\special{pa 1458 1448}%
\special{pa 1456 1478}%
\special{pa 1450 1508}%
\special{pa 1444 1538}%
\special{pa 1436 1566}%
\special{pa 1428 1596}%
\special{pa 1418 1626}%
\special{pa 1408 1656}%
\special{pa 1396 1686}%
\special{pa 1382 1716}%
\special{pa 1370 1746}%
\special{pa 1354 1776}%
\special{pa 1340 1806}%
\special{pa 1324 1836}%
\special{pa 1308 1866}%
\special{pa 1292 1896}%
\special{pa 1274 1926}%
\special{pa 1256 1956}%
\special{pa 1240 1986}%
\special{pa 1222 2016}%
\special{pa 1204 2046}%
\special{pa 1200 2050}%
\special{sp}%
%
\special{pn 13}%
\special{pa 1380 760}%
\special{pa 1362 790}%
\special{pa 1342 820}%
\special{pa 1322 848}%
\special{pa 1304 878}%
\special{pa 1286 906}%
\special{pa 1268 936}%
\special{pa 1250 966}%
\special{pa 1232 994}%
\special{pa 1216 1024}%
\special{pa 1202 1054}%
\special{pa 1188 1084}%
\special{pa 1174 1112}%
\special{pa 1162 1142}%
\special{pa 1152 1172}%
\special{pa 1142 1200}%
\special{pa 1134 1230}%
\special{pa 1128 1260}%
\special{pa 1124 1290}%
\special{pa 1122 1320}%
\special{pa 1120 1350}%
\special{pa 1120 1378}%
\special{pa 1124 1408}%
\special{pa 1128 1438}%
\special{pa 1134 1468}%
\special{pa 1142 1498}%
\special{pa 1152 1528}%
\special{pa 1162 1558}%
\special{pa 1172 1588}%
\special{pa 1184 1618}%
\special{pa 1196 1650}%
\special{pa 1210 1680}%
\special{pa 1222 1710}%
\special{pa 1236 1740}%
\special{pa 1248 1772}%
\special{pa 1260 1802}%
\special{pa 1274 1832}%
\special{pa 1284 1864}%
\special{pa 1296 1894}%
\special{pa 1306 1926}%
\special{pa 1314 1958}%
\special{pa 1322 1988}%
\special{pa 1328 2020}%
\special{pa 1334 2052}%
\special{pa 1340 2082}%
\special{pa 1344 2114}%
\special{pa 1348 2146}%
\special{pa 1350 2178}%
\special{pa 1352 2210}%
\special{pa 1354 2242}%
\special{pa 1354 2274}%
\special{pa 1354 2306}%
\special{pa 1354 2338}%
\special{pa 1352 2370}%
\special{pa 1352 2402}%
\special{pa 1350 2434}%
\special{pa 1348 2468}%
\special{pa 1346 2500}%
\special{pa 1344 2532}%
\special{pa 1342 2564}%
\special{pa 1338 2596}%
\special{pa 1336 2628}%
\special{pa 1332 2662}%
\special{pa 1330 2680}%
\special{sp}%
%
\special{pn 13}%
\special{pa 1830 740}%
\special{pa 1848 770}%
\special{pa 1868 800}%
\special{pa 1886 828}%
\special{pa 1902 858}%
\special{pa 1920 888}%
\special{pa 1938 916}%
\special{pa 1954 946}%
\special{pa 1970 976}%
\special{pa 1986 1004}%
\special{pa 2000 1034}%
\special{pa 2014 1064}%
\special{pa 2028 1092}%
\special{pa 2040 1122}%
\special{pa 2050 1152}%
\special{pa 2060 1182}%
\special{pa 2068 1210}%
\special{pa 2076 1240}%
\special{pa 2082 1270}%
\special{pa 2086 1300}%
\special{pa 2090 1330}%
\special{pa 2090 1358}%
\special{pa 2090 1388}%
\special{pa 2088 1418}%
\special{pa 2086 1448}%
\special{pa 2080 1478}%
\special{pa 2074 1508}%
\special{pa 2066 1536}%
\special{pa 2058 1566}%
\special{pa 2048 1596}%
\special{pa 2038 1626}%
\special{pa 2026 1656}%
\special{pa 2012 1686}%
\special{pa 2000 1716}%
\special{pa 1984 1746}%
\special{pa 1970 1776}%
\special{pa 1954 1806}%
\special{pa 1938 1836}%
\special{pa 1922 1866}%
\special{pa 1904 1896}%
\special{pa 1886 1926}%
\special{pa 1870 1956}%
\special{pa 1852 1986}%
\special{pa 1834 2016}%
\special{pa 1830 2020}%
\special{sp}%
%
\special{pn 13}%
\special{pa 2670 790}%
\special{pa 2652 820}%
\special{pa 2632 850}%
\special{pa 2612 878}%
\special{pa 2594 908}%
\special{pa 2576 936}%
\special{pa 2558 966}%
\special{pa 2540 996}%
\special{pa 2522 1024}%
\special{pa 2506 1054}%
\special{pa 2492 1084}%
\special{pa 2478 1114}%
\special{pa 2464 1142}%
\special{pa 2452 1172}%
\special{pa 2442 1202}%
\special{pa 2432 1230}%
\special{pa 2424 1260}%
\special{pa 2418 1290}%
\special{pa 2414 1320}%
\special{pa 2412 1350}%
\special{pa 2410 1380}%
\special{pa 2410 1408}%
\special{pa 2414 1438}%
\special{pa 2418 1468}%
\special{pa 2424 1498}%
\special{pa 2432 1528}%
\special{pa 2442 1558}%
\special{pa 2452 1588}%
\special{pa 2462 1618}%
\special{pa 2474 1648}%
\special{pa 2486 1680}%
\special{pa 2500 1710}%
\special{pa 2512 1740}%
\special{pa 2526 1770}%
\special{pa 2538 1802}%
\special{pa 2550 1832}%
\special{pa 2564 1862}%
\special{pa 2574 1894}%
\special{pa 2586 1924}%
\special{pa 2596 1956}%
\special{pa 2604 1988}%
\special{pa 2612 2018}%
\special{pa 2618 2050}%
\special{pa 2624 2082}%
\special{pa 2630 2112}%
\special{pa 2634 2144}%
\special{pa 2638 2176}%
\special{pa 2640 2208}%
\special{pa 2642 2240}%
\special{pa 2644 2272}%
\special{pa 2644 2304}%
\special{pa 2644 2336}%
\special{pa 2644 2368}%
\special{pa 2642 2400}%
\special{pa 2642 2432}%
\special{pa 2640 2464}%
\special{pa 2638 2498}%
\special{pa 2636 2530}%
\special{pa 2634 2562}%
\special{pa 2632 2594}%
\special{pa 2628 2626}%
\special{pa 2626 2658}%
\special{pa 2622 2692}%
\special{pa 2620 2710}%
\special{sp}%
%
\special{pn 13}%
\special{pa 2500 750}%
\special{pa 2518 780}%
\special{pa 2538 810}%
\special{pa 2556 838}%
\special{pa 2572 868}%
\special{pa 2590 898}%
\special{pa 2608 926}%
\special{pa 2624 956}%
\special{pa 2640 986}%
\special{pa 2656 1014}%
\special{pa 2670 1044}%
\special{pa 2684 1074}%
\special{pa 2698 1102}%
\special{pa 2710 1132}%
\special{pa 2720 1162}%
\special{pa 2730 1192}%
\special{pa 2738 1220}%
\special{pa 2746 1250}%
\special{pa 2752 1280}%
\special{pa 2756 1310}%
\special{pa 2760 1340}%
\special{pa 2760 1368}%
\special{pa 2760 1398}%
\special{pa 2758 1428}%
\special{pa 2756 1458}%
\special{pa 2750 1488}%
\special{pa 2744 1518}%
\special{pa 2736 1546}%
\special{pa 2728 1576}%
\special{pa 2718 1606}%
\special{pa 2708 1636}%
\special{pa 2696 1666}%
\special{pa 2682 1696}%
\special{pa 2670 1726}%
\special{pa 2654 1756}%
\special{pa 2640 1786}%
\special{pa 2624 1816}%
\special{pa 2608 1846}%
\special{pa 2592 1876}%
\special{pa 2574 1906}%
\special{pa 2556 1936}%
\special{pa 2540 1966}%
\special{pa 2522 1996}%
\special{pa 2504 2026}%
\special{pa 2500 2030}%
\special{sp}%
%
\special{pn 13}%
\special{pa 3290 790}%
\special{pa 3272 820}%
\special{pa 3252 850}%
\special{pa 3232 878}%
\special{pa 3214 908}%
\special{pa 3196 936}%
\special{pa 3178 966}%
\special{pa 3160 996}%
\special{pa 3142 1024}%
\special{pa 3126 1054}%
\special{pa 3112 1084}%
\special{pa 3098 1114}%
\special{pa 3084 1142}%
\special{pa 3072 1172}%
\special{pa 3062 1202}%
\special{pa 3052 1230}%
\special{pa 3044 1260}%
\special{pa 3038 1290}%
\special{pa 3034 1320}%
\special{pa 3032 1350}%
\special{pa 3030 1380}%
\special{pa 3030 1408}%
\special{pa 3034 1438}%
\special{pa 3038 1468}%
\special{pa 3044 1498}%
\special{pa 3052 1528}%
\special{pa 3062 1558}%
\special{pa 3072 1588}%
\special{pa 3082 1618}%
\special{pa 3094 1648}%
\special{pa 3106 1680}%
\special{pa 3120 1710}%
\special{pa 3132 1740}%
\special{pa 3146 1770}%
\special{pa 3158 1802}%
\special{pa 3170 1832}%
\special{pa 3184 1862}%
\special{pa 3194 1894}%
\special{pa 3206 1924}%
\special{pa 3216 1956}%
\special{pa 3224 1988}%
\special{pa 3232 2018}%
\special{pa 3238 2050}%
\special{pa 3244 2082}%
\special{pa 3250 2112}%
\special{pa 3254 2144}%
\special{pa 3258 2176}%
\special{pa 3260 2208}%
\special{pa 3262 2240}%
\special{pa 3264 2272}%
\special{pa 3264 2304}%
\special{pa 3264 2336}%
\special{pa 3264 2368}%
\special{pa 3262 2400}%
\special{pa 3262 2432}%
\special{pa 3260 2464}%
\special{pa 3258 2498}%
\special{pa 3256 2530}%
\special{pa 3254 2562}%
\special{pa 3252 2594}%
\special{pa 3248 2626}%
\special{pa 3246 2658}%
\special{pa 3242 2692}%
\special{pa 3240 2710}%
\special{sp}%
%
\special{pn 13}%
\special{pa 3120 760}%
\special{pa 3138 790}%
\special{pa 3158 820}%
\special{pa 3176 848}%
\special{pa 3192 878}%
\special{pa 3210 908}%
\special{pa 3228 936}%
\special{pa 3244 966}%
\special{pa 3260 996}%
\special{pa 3276 1024}%
\special{pa 3290 1054}%
\special{pa 3304 1084}%
\special{pa 3318 1112}%
\special{pa 3330 1142}%
\special{pa 3340 1172}%
\special{pa 3350 1202}%
\special{pa 3358 1230}%
\special{pa 3366 1260}%
\special{pa 3372 1290}%
\special{pa 3376 1320}%
\special{pa 3380 1350}%
\special{pa 3380 1378}%
\special{pa 3380 1408}%
\special{pa 3378 1438}%
\special{pa 3376 1468}%
\special{pa 3370 1498}%
\special{pa 3364 1528}%
\special{pa 3356 1556}%
\special{pa 3348 1586}%
\special{pa 3338 1616}%
\special{pa 3328 1646}%
\special{pa 3316 1676}%
\special{pa 3302 1706}%
\special{pa 3290 1736}%
\special{pa 3274 1766}%
\special{pa 3260 1796}%
\special{pa 3244 1826}%
\special{pa 3228 1856}%
\special{pa 3212 1886}%
\special{pa 3194 1916}%
\special{pa 3176 1946}%
\special{pa 3160 1976}%
\special{pa 3142 2006}%
\special{pa 3124 2036}%
\special{pa 3120 2040}%
\special{sp}%
\put(20.1000,-19.3000){\makebox(0,0)[lb]{$\Theta_{2,1}$}}%
\put(19.9000,-22.5000){\makebox(0,0)[lb]{$\Theta_{2,0}$}}%
\put(26.6000,-24.9000){\makebox(0,0)[lb]{$\Theta_{3,0}$}}%
\put(26.4000,-19.9000){\makebox(0,0)[lb]{$\Theta_{3,1}$}}%
%
\special{pn 13}%
\special{pa 3710 750}%
\special{pa 3728 780}%
\special{pa 3748 810}%
\special{pa 3766 838}%
\special{pa 3782 868}%
\special{pa 3800 898}%
\special{pa 3818 926}%
\special{pa 3834 956}%
\special{pa 3850 986}%
\special{pa 3866 1014}%
\special{pa 3880 1044}%
\special{pa 3894 1074}%
\special{pa 3908 1102}%
\special{pa 3920 1132}%
\special{pa 3930 1162}%
\special{pa 3940 1192}%
\special{pa 3948 1220}%
\special{pa 3956 1250}%
\special{pa 3962 1280}%
\special{pa 3966 1310}%
\special{pa 3970 1340}%
\special{pa 3970 1368}%
\special{pa 3970 1398}%
\special{pa 3968 1428}%
\special{pa 3966 1458}%
\special{pa 3960 1488}%
\special{pa 3954 1518}%
\special{pa 3946 1546}%
\special{pa 3938 1576}%
\special{pa 3928 1606}%
\special{pa 3918 1636}%
\special{pa 3906 1666}%
\special{pa 3892 1696}%
\special{pa 3880 1726}%
\special{pa 3864 1756}%
\special{pa 3850 1786}%
\special{pa 3834 1816}%
\special{pa 3818 1846}%
\special{pa 3802 1876}%
\special{pa 3784 1906}%
\special{pa 3766 1936}%
\special{pa 3750 1966}%
\special{pa 3732 1996}%
\special{pa 3714 2026}%
\special{pa 3710 2030}%
\special{sp}%
%
\special{pn 13}%
\special{pa 3880 750}%
\special{pa 3862 780}%
\special{pa 3842 810}%
\special{pa 3822 838}%
\special{pa 3804 868}%
\special{pa 3786 896}%
\special{pa 3768 926}%
\special{pa 3750 956}%
\special{pa 3732 984}%
\special{pa 3716 1014}%
\special{pa 3702 1044}%
\special{pa 3688 1074}%
\special{pa 3674 1102}%
\special{pa 3662 1132}%
\special{pa 3652 1162}%
\special{pa 3642 1190}%
\special{pa 3634 1220}%
\special{pa 3628 1250}%
\special{pa 3624 1280}%
\special{pa 3622 1310}%
\special{pa 3620 1340}%
\special{pa 3620 1368}%
\special{pa 3624 1398}%
\special{pa 3628 1428}%
\special{pa 3634 1458}%
\special{pa 3642 1488}%
\special{pa 3652 1518}%
\special{pa 3662 1548}%
\special{pa 3672 1578}%
\special{pa 3684 1608}%
\special{pa 3696 1640}%
\special{pa 3710 1670}%
\special{pa 3722 1700}%
\special{pa 3736 1730}%
\special{pa 3748 1762}%
\special{pa 3760 1792}%
\special{pa 3774 1822}%
\special{pa 3784 1854}%
\special{pa 3796 1884}%
\special{pa 3806 1916}%
\special{pa 3814 1948}%
\special{pa 3822 1978}%
\special{pa 3828 2010}%
\special{pa 3834 2042}%
\special{pa 3840 2072}%
\special{pa 3844 2104}%
\special{pa 3848 2136}%
\special{pa 3850 2168}%
\special{pa 3852 2200}%
\special{pa 3854 2232}%
\special{pa 3854 2264}%
\special{pa 3854 2296}%
\special{pa 3854 2328}%
\special{pa 3852 2360}%
\special{pa 3852 2392}%
\special{pa 3850 2424}%
\special{pa 3848 2458}%
\special{pa 3846 2490}%
\special{pa 3844 2522}%
\special{pa 3842 2554}%
\special{pa 3838 2586}%
\special{pa 3836 2618}%
\special{pa 3832 2652}%
\special{pa 3830 2670}%
\special{sp}%
\put(39.4000,-19.6000){\makebox(0,0)[lb]{$\Theta_{5,1}$}}%
\put(39.0000,-23.5000){\makebox(0,0)[lb]{$\Theta_{5,0}$}}%
%
\special{pn 13}%
\special{pa 1260 1070}%
\special{pa 1780 1070}%
\special{fp}%
%
\special{pn 13}%
\special{pa 1880 1070}%
\special{pa 2440 1070}%
\special{fp}%
%
\special{pn 13}%
\special{pa 2560 1070}%
\special{pa 3200 1060}%
\special{fp}%
%
\special{pn 13}%
\special{pa 3340 1060}%
\special{pa 3810 1070}%
\special{fp}%
%
\special{pn 13}%
\special{pa 3950 1070}%
\special{pa 5030 1060}%
\special{fp}%
%
\special{pn 13}%
\special{pa 1230 1240}%
\special{pa 1720 1230}%
\special{fp}%
%
\special{pn 13}%
\special{pa 1820 1240}%
\special{pa 2710 1250}%
\special{fp}%
%
\special{pn 13}%
\special{pa 3090 1250}%
\special{pa 3880 1250}%
\special{fp}%
%
\special{pn 13}%
\special{pa 4030 1250}%
\special{pa 5150 1240}%
\special{fp}%
%
\special{pn 13}%
\special{pa 1240 1380}%
\special{pa 1680 1370}%
\special{fp}%
%
\special{pn 13}%
\special{pa 1850 1370}%
\special{pa 2690 1380}%
\special{fp}%
%
\special{pn 13}%
\special{pa 2840 1370}%
\special{pa 3340 1380}%
\special{fp}%
%
\special{pn 13}%
\special{pa 3440 1370}%
\special{pa 3580 1370}%
\special{fp}%
%
\special{pn 13}%
\special{pa 3700 1360}%
\special{pa 5160 1360}%
\special{fp}%
\put(51.3000,-10.8000){\makebox(0,0)[lb]{$s_{L_0}^+$}}%
%
\special{pn 13}%
\special{pa 2790 1250}%
\special{pa 2980 1240}%
\special{fp}%
\end{picture}%

\end{center}
Here we assume that $\Theta_{1, 0}$ and $O$ come from $z_o$. 
By the explicit formula of $\langle\, , \, \rangle$,  we have
\[
\langle s_{L_i}^{\pm}, s_{L_i}^{\pm} \rangle  = \frac 12, \, \,  (i = 0, 1, 2) \quad \langle s_{L_i}^+, s_{L_j}^+\rangle = 0. \, \, (i, j = 0, 1, 2, i \neq j)
\]
By \cite{oguiso-shioda}, $\MW(S_{(\mcQ, z_o)}) \cong (A_1^*)^{\oplus 3} \oplus (\ZZ/2\ZZ)$ and we may assume that
\[
(A_1^*)^{\oplus 2} \cong \ZZ s_{L_0}^+ \oplus \ZZ s_{L_1}^+ \oplus \ZZ s_{L_2}^+,
\]
and that the unique $2$-torsion section arises from $C_1$.

By \cite[Theorem~9.1]{kodaira},  
$[2]s_{L_i}^{\pm}$ $(i = 0, 1, 2)$ meet
the identity component at each singular fiber. Hence  by
the explicit formula for $\langle\, , \, \rangle$, we have $[2]s_{L_i}^{\pm}O = 0$ for each $i$. This implies that, for each $i$, 
$C_{L_i}: = \overline{q}\circ f([2]s_{L_i}^{\pm})$ is a conic not passing through $P_j$ 
$(j = 1, 2, 3, 4)$.  
If $C_{L_i}$ and  $C_1 + C_3$ has an intersection point at which
intersection multiplicity is odd, then we easily see that the closure of
$(\overline{q}\circ f)^{-1}(C_{L_i}\setminus z_o)$ is irreducible. This is impossible, as $C_{L_i}$
is the image of $[2]s_{L_i}^\pm$.
Hence we have three conic satisfying the first two conditions. 

Conversely, suppose that there exists a conic $C_o$ satisfying the first two conditions.
We infer that $C_o$ gives rise to two sections $s_{C_o}^{\pm}$. Since $C_o$ does not pass through
any singularities of $C_1 + C_2$ and $s_{C_o}^{\pm}O = 0$, we have $\langle s_{C_o}^{\pm}, s_{C_o}^{\pm} \rangle = 2$.
On the other hand, any element $s \in \MW(S_{(\mcQ,z_o)})$ with  $\langle s, s \rangle = 2$ is of the form
\[
[2]s_{L_i}^{\pm} \dot{+} \tau, \quad (i = 0, 1, 2) \quad \tau \in \MW(S_{(\mcQ,z_o)})_{\tor}.
\]
By a similar argument to that in the case of Line-conic arrangement 1,   we infer $\tau = 0$. Hence
we infer that $C_{L_i}$ are all  conics satisfying the first two conditions.
Now put $C_3^{(i)} := C_{L_i},  s_{C_3^{(i)}}:= [2]s_{L_i}^+$ $(i = 1, 2)$.
For $(i, j) = (1, 2), (2, 1)$,  we have 
\[
s_{C_3^{(i)}}  \dot{+} [p-2]s_{L_i}^+ \in [p]\MW(S_{(\mcQ,z_o)})
\]
for any odd prime $p$, while
\[
s_{C_3^{(i)}} \dot{+} [k]s_{L_j}^+ \not\in [p]\MW(S_{(\mcQ,z_o)})
\]
for any odd prime $p$ and $1 \le k \le p-1$.

 By Proposition~\ref{prop:res-1}, 
 if both of $C_1 + C_2 + C_3^{(i)} + L_1$ and 
$C_1 + C_2 + C_3^{(i)} + L_2$ have the combinatorics for Line-conic arrangement 2 of the
same type, then 
$(C_1+ C_2 +C_3^{(i)} + L_1, C_1+ C_2 +C_3^{(i)} + L_2)$ is a Zariski pair for each
$i = 1, 2$. \proofend

\bigskip

\begin{rem}\label{rem:fundamental}{\rm Let $B$ be the line-conic arrangement as in Theorem~\ref{thm:zpair}. By Corollary~\ref{cor:split}, if
there exists a $D_{2p}$-cover $\pi: X \to \PP^2$ with branch locus $B$,
then $\Delta_{\beta_1(\pi)} = L_1 + L_2 + C_1$ (resp. $C_1 + C_2$) 
for Line-conic arrangement 1 (resp. 2). This means that the $D_{2p}$-covers in our proof of Theorem~\ref{thm:zpair} are only possible ones.
Therefore for the fundamental group $\pi_1(\PP^2\setminus B, \ast)$, 
we infer that 
\[
\pi_1(\PP^2\setminus(L_1 + L_2 + L_3 + C_1 + C_2), \ast)
\not \cong 
\pi_1(\PP^2\setminus(L_1 + L_2 + L_4 + C_1 + C_2), \ast)
\]
for Line-conic arrangement 1, and 
\[
\pi_1(\PP^2\setminus(L_1  + C_1 + C_2 + C_3^{(i)}), \ast)
\not \cong 
\pi_1(\PP^2\setminus(L_2 + C_1 + C_2+ C_3^{(i)}), \ast)
\]
for Line-conic arrangement 2.
}
\end{rem}

\begin{exmple}\label{eg:line-conic1}{\rm Let $[T, X, Z]$ be  homogeneous
coordinates of $\PP^2$ and 
let $(t, x):= (T/Z, X/Z)$ be  affine coordinates of $\PP^2$ and 
consider a conic and four lines as follows:
\[
\begin{array}{lll}
C_1: x - t^2 = 0,  & L_1: x - 3t + 2 = 0, & L_2: x + 3t + 2 = 0, \\
L_3: x - t -2 = 0,  & L_4: x -1 = 0. &  \\
\end{array}
\]
Note that $C_1\cap(L_1\cup L_2) = \{[\pm1,1,1], [\pm 2, 4, 1]\}$.  Put $\mcQ = C_1 + L_1 + L_2$ and choose $[0, 1, 0]$ as the distinguished point 
$z_o$. Let $S_{(\mcQ, z_o)}$ be the rational elliptic surface obtained as in \S 4. Then its generic fiber is 
an elliptic curve over $\CC(t)$ given by the Weierstrass equation:
\[
y^2 = (x - t^2)(x - 3t + 2)(x + 3t+ 2).
\]
Under this setting, we may assume that  the sections $s_{L_i}^{\pm} (i = 3, 4)$ are as follows:
\[
\begin{array}{cc}
s_{L_3}^{\pm} = (t+2, \pm2\sqrt{2}(t-2)(t+1)) & s_{L_4}^{\pm} = (1, \pm 3(t+1)(t-1)).
\end{array}
\]
Hence we have
\[
\begin{array}{cc}
[2]s_{L_3}^+ = (\frac 98 t^2, \frac 1{32} t \sqrt{2} (9 t^2-16)),  & [2]s_{L_4}^+ = (t^2+\frac 14, \frac12 t^2- \frac 98)
\end{array}
\]
Now put 
\[
\begin{array}{cc}
C_2 : x - \frac 98 t^2 = 0, & C'_2 : x - t^2 - \frac 14 = 0. 
\end{array}
\]
Then $(\mcQ + C_2 + L_3, Q + C_2 + L_4)$ is a Zariski pair for Line-conic arrangement 1 of
type $(a)$, and $(Q + C'_2 + L_3,   Q  + C'_2 + L_4)$ is a Zariski pair for Line-conic
arrangement 1 of type $(b)$.
}
\end{exmple}

\begin{exmple}\label{eg:line-conic2}{\rm We keep the same coordinates as Example~\ref{eg:line-conic1}. 

{\sl Line-conic arrangement 2 of type $(a)$.} 
Consider
two conics and two lines:
\[
\begin{array}{ll}
C_1:  x - t^2 + 2 =  0,  & C_2: x^2 -2x + t^2 -4 = 0, \\ 
L_1: x - t = 0,  & L_2: x - 3t + 4 = 0.
\end{array}
\]
Note that $C_1\cap C_2 = \{[\pm 2, 2, 1], [\pm 1, -1, 1]\}$.  Put $\mcQ = C_1 + C_2$ and choose $[0, 1, 0]$ as the distinguished point 
$z_o$. Let $S_{(\mcQ, z_o)}$ be the rational elliptic surface obtained as before. Then its generic fiber is 
an elliptic curve over $\CC(t)$ given by the Weierstrass equation:
\[
y^2 = (x - t^2+2)(x^2 - 2x  + t^2 -4).
\]
Then we may assume that the sections $s_{L_2}^{\pm} (i = 1, 2)$ are as follows:
\[
\begin{array}{cc}
s_{L_1}^{\pm} = (t, \pm \sqrt {-2} (t + 1) (t-2)),  & s_{L_2}^{\pm} = (3t - 4, \pm \sqrt{-10}(t-1)(t-2)).
\end{array}
\]
Thus we have
\[
\begin{array}{cc}
[2]s_{L_1}^+ = (\frac 12 t^2-2, -\frac14\sqrt{-2} t (t^2-4)), & 
[2]s_{L_2}^+ = (\frac 1{10}t^2 -2,  -\frac 3{100}\sqrt{-10}t(t^2+20)).
\end{array}
\]
Now we put
\[
\begin{array}{cc}
C_3 : x - \frac 12 t^2 + 2 = 0 , & C'_3: x - \frac 1{10}t^2 + 2 = 0
\end{array}
\]
Then both $(\mcQ + C_3 + L_1, \mcQ + C_3 + L_2)$ and $(\mcQ + C'_3 + L_1, Q +C'_3 + L_2)$ are Zariski pairs for
Line-conic arrangement 2 of type $(a)$.

\bigskip

{\sl Line-conic arrangement 2 of type $(b)$.}
Consider
two conics and two lines:
\[
\begin{array}{ll}
C_1:  x - t^2 + 2 =  0,  & C_2: x^2 -2x + t^2 -4 = 0, \\ 
L_1: x - t = 0,  & L_2: x  +1 = 0.
\end{array}
\]
Put $\mcQ = C_1 + C_2$ and choose $[0, 1, 0]$ as the distinguished point 
$z_o$. Let $S_{(\mcQ, z_o)}$ be the rational elliptic surface obtained as before. Then its generic fiber is 
an elliptic curve over $\CC(t)$ given by the Weierstrass equation:
\[
y^2 = (x - t^2+2)(x^2 - 2x  + t^2 -4).
\]
Then we may assume that the sections $s_{L_2}^{\pm} (i = 1, 2)$ is as follows:
\[
 s_{L_2}^{\pm} = (-1, \pm \sqrt{-1}(t-1)(t+1)).
 \]
Thus we have
\[
[2]s_{L_2}^+ = \left (t^2 - \frac{17}4,  \frac 38\sqrt{-1}(4t^2 - 19)\right ).
\]
Now we put
\[
C_3 : x - t^2 + \frac{17}4 = 0.
\]
As $C_3$ is tangent to $C_1$ (resp. $C_2$) at one point (resp. two distinct points), we infer that
 $(\mcQ + C_3 + L_1, \mcQ + C_3 + L_2)$ is a  Zariski pair for
Line-conic arrangement 2 of type $(b)$.

\bigskip

{\sl Line-conic arrangement 2 of type $(c)$.}
Consider
two conics and two lines:
\[
\begin{array}{ll}
C_1:  x - t^2 + \frac 12 =  0,  & C_2: x^2 -x + t^2  = 0, \\ 
L_1: x = \frac 1{\sqrt 2},  & L_2: \frac {\sqrt 2}4\left ( \sqrt {-1}c_1- c_2\right ) x
 + t - \frac{1}4\left ( \sqrt{-1}c_1+ c_2 \right )= 0, 
\end{array}
\]
where $c_1 = \sqrt {2 + 2\sqrt 2}, c_2 = \sqrt{-2 + 2\sqrt 2}$.
Note that 
\[
C_1\cap  C_2  = \left \{
\left [\pm \sqrt{-1/2 + 1/\sqrt{2}}, 1/\sqrt{2}, 1\right ], \left [\pm \sqrt{-1/2 - 1/\sqrt{2}}, - 1/\sqrt{2}, 1\right] \right \}.
\]
Put $\mcQ = C_1 + C_2$ and choose $[0, 1, 0]$ as the distinguished point 
$z_o$. Let $S_{(\mcQ, z_o)}$ be the rational elliptic surface obtained as before. Then its generic fiber is 
an elliptic curve over $\CC(t)$ given by the Weierstrass equation:
\[
y^2 = \left (x - t^2- \frac 12\right)(x^2 - x  + t^2).
\]
Then we may assume that the sections $s_{L_1}^{\pm}$ are as follows:
\[
s_{L_1}^{\pm} = \left (\frac 1{\sqrt 2},  \pm \frac {\sqrt{-1}}2(-2t^2 - 1 +\sqrt {2}) \right ).
\]
Thus we have
\[
[2]s_{L_1}^+ = \left ( t^2, \sqrt{- \frac12}t^2\right ).
\]
 Now we put
\[
C_3 : x - t^2 = 0 
\]
Then  $(\mcQ + C_3 + L_1, \mcQ + C_3 + L_2)$ is a Zariski pair for
Line-conic arrangement 2 of type $(c)$.
}
\end{exmple}

We end this section by giving another example:

\begin{prop}\label{prop:eg-3}{Let ${\mathcal Q}$ be an irreducible quartic with a ${\mathbb D}_4$ singularity, $P$. Let
$z_o$ be a point on  ${\mathcal Q}$ such that the tangent line $L_{z_o}$ at 
$z_o$ meets ${\mathcal Q}$ with two other 
distinct points. Let $L_1, L_2$ and $L_3$ be three tangent lines which meet ${\mathcal Q}$ at $P$ with multiplicity $4$. 
Then the following statements hold:

\begin{enumerate}

\item[(i)] For each $L_i$, there exists a unique conic $C_i$ such that $(a)$ $z_o \in C_i$, $(b)$ $P \not\in C_i$ and $(c)$ for $\forall x \in C_i\cap {\mathcal Q}$, $I_x(C_i, \mcQ)$ is even.

\item[(ii)] For any odd prime $p$, there exists a $D_{2p}$-cover of $\PP^2$ branched at $2{\mathcal Q} + p(C_i + L_i)$ for 
each $i = 1, 2,3$,
while there exists no $D_{2p}$-cover of $\PP^2$ branched at $2{\mathcal Q} + p(C_i + L_j)$ for any $i, j$ $(i \neq j)$.

\end{enumerate}
}
\end{prop}

\proof (i) Let $f''_{\mcQ} \to \PP^2$ be  a double cover with branch locus $\mcQ$ and let $\varphi_{z_o} : S_{(\mcQ, z_o)}
\to \PP^1$ be the rational elliptic surface obtained as in \S 4. By our assumption on $\mcQ$ and $z_o$, the configuration
of reducible singular fiber of $\varphi_{z_o}$ is $\I_0^*, \I_2$ and three lines $L_i$ $(i = 1, 2,3)$ give rise to sections $s_{L_i}^{\pm}$
($i = 1, 2,3$), respectively. By labeling irreducible components of singular fibers suitably, we have the following picture for
$s_{L_i}^+ (i = 1, 2, 3)$:

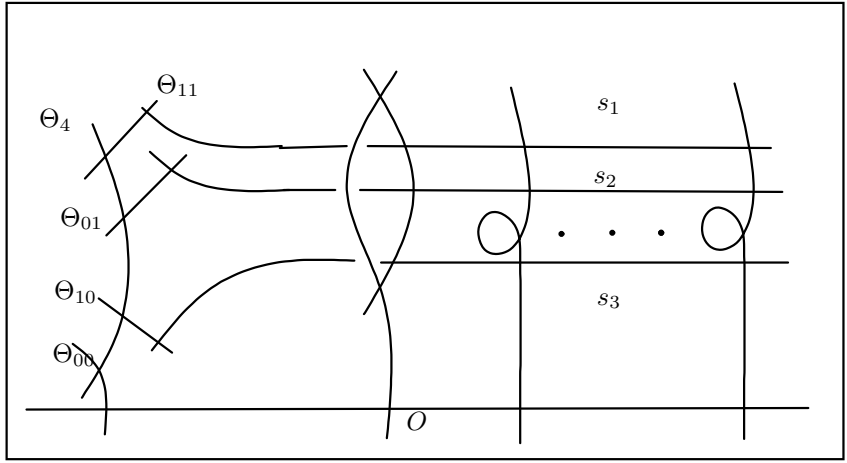
\begin{center}
\unitlength 0.1in
\begin{picture}( 43.8000, 23.9000)(  9.8000,-28.1000)
%
\special{pn 13}%
\special{pa 980 420}%
\special{pa 5360 420}%
\special{pa 5360 2810}%
\special{pa 980 2810}%
\special{pa 980 420}%
\special{fp}%
%
\special{pn 13}%
\special{pa 1084 2550}%
\special{pa 5176 2542}%
\special{fp}%
%
\special{pn 13}%
\special{pa 1326 2206}%
\special{pa 1352 2226}%
\special{pa 1378 2246}%
\special{pa 1402 2268}%
\special{pa 1424 2290}%
\special{pa 1444 2314}%
\special{pa 1462 2340}%
\special{pa 1474 2366}%
\special{pa 1484 2396}%
\special{pa 1492 2424}%
\special{pa 1498 2456}%
\special{pa 1500 2488}%
\special{pa 1502 2520}%
\special{pa 1502 2554}%
\special{pa 1500 2588}%
\special{pa 1498 2624}%
\special{pa 1496 2658}%
\special{pa 1494 2680}%
\special{sp}%
%
\special{pn 13}%
\special{pa 1430 1056}%
\special{pa 1442 1086}%
\special{pa 1454 1116}%
\special{pa 1466 1146}%
\special{pa 1478 1176}%
\special{pa 1490 1204}%
\special{pa 1500 1234}%
\special{pa 1512 1264}%
\special{pa 1522 1294}%
\special{pa 1534 1326}%
\special{pa 1544 1356}%
\special{pa 1552 1386}%
\special{pa 1562 1416}%
\special{pa 1570 1448}%
\special{pa 1578 1478}%
\special{pa 1586 1510}%
\special{pa 1592 1542}%
\special{pa 1598 1574}%
\special{pa 1604 1606}%
\special{pa 1608 1638}%
\special{pa 1612 1670}%
\special{pa 1614 1704}%
\special{pa 1616 1736}%
\special{pa 1618 1768}%
\special{pa 1618 1802}%
\special{pa 1618 1834}%
\special{pa 1616 1866}%
\special{pa 1614 1898}%
\special{pa 1610 1930}%
\special{pa 1606 1964}%
\special{pa 1600 1994}%
\special{pa 1594 2026}%
\special{pa 1588 2058}%
\special{pa 1580 2088}%
\special{pa 1570 2120}%
\special{pa 1560 2150}%
\special{pa 1548 2178}%
\special{pa 1536 2208}%
\special{pa 1522 2238}%
\special{pa 1508 2266}%
\special{pa 1492 2294}%
\special{pa 1478 2322}%
\special{pa 1460 2350}%
\special{pa 1444 2378}%
\special{pa 1428 2404}%
\special{pa 1410 2432}%
\special{pa 1392 2458}%
\special{pa 1376 2486}%
\special{pa 1374 2488}%
\special{sp}%
%
\special{pn 13}%
\special{pa 1766 934}%
\special{pa 1390 1340}%
\special{fp}%
%
\special{pn 13}%
\special{pa 1920 1218}%
\special{pa 1502 1638}%
\special{fp}%
%
\special{pn 13}%
\special{pa 1462 1968}%
\special{pa 1846 2252}%
\special{fp}%
%
\special{pn 13}%
\special{pa 3620 864}%
\special{pa 3628 896}%
\special{pa 3636 928}%
\special{pa 3642 958}%
\special{pa 3650 990}%
\special{pa 3658 1022}%
\special{pa 3666 1052}%
\special{pa 3672 1084}%
\special{pa 3678 1116}%
\special{pa 3686 1146}%
\special{pa 3692 1178}%
\special{pa 3696 1210}%
\special{pa 3702 1240}%
\special{pa 3706 1272}%
\special{pa 3710 1302}%
\special{pa 3714 1334}%
\special{pa 3716 1366}%
\special{pa 3718 1396}%
\special{pa 3718 1428}%
\special{pa 3716 1460}%
\special{pa 3714 1494}%
\special{pa 3708 1528}%
\special{pa 3702 1562}%
\special{pa 3692 1594}%
\special{pa 3678 1628}%
\special{pa 3664 1658}%
\special{pa 3644 1684}%
\special{pa 3622 1708}%
\special{pa 3596 1724}%
\special{pa 3564 1734}%
\special{pa 3532 1736}%
\special{pa 3502 1728}%
\special{pa 3476 1708}%
\special{pa 3460 1678}%
\special{pa 3450 1644}%
\special{pa 3448 1608}%
\special{pa 3458 1574}%
\special{pa 3476 1544}%
\special{pa 3500 1522}%
\special{pa 3528 1514}%
\special{pa 3560 1518}%
\special{pa 3590 1534}%
\special{pa 3618 1556}%
\special{pa 3640 1584}%
\special{pa 3654 1614}%
\special{pa 3662 1644}%
\special{pa 3666 1676}%
\special{pa 3668 1708}%
\special{pa 3668 1742}%
\special{pa 3668 1774}%
\special{pa 3668 1806}%
\special{pa 3668 1840}%
\special{pa 3670 1872}%
\special{pa 3670 1904}%
\special{pa 3670 1936}%
\special{pa 3670 1968}%
\special{pa 3670 2000}%
\special{pa 3670 2032}%
\special{pa 3670 2064}%
\special{pa 3670 2096}%
\special{pa 3670 2128}%
\special{pa 3670 2160}%
\special{pa 3670 2192}%
\special{pa 3670 2224}%
\special{pa 3670 2256}%
\special{pa 3670 2288}%
\special{pa 3670 2320}%
\special{pa 3668 2350}%
\special{pa 3668 2382}%
\special{pa 3668 2414}%
\special{pa 3668 2446}%
\special{pa 3668 2478}%
\special{pa 3668 2510}%
\special{pa 3668 2542}%
\special{pa 3668 2574}%
\special{pa 3668 2606}%
\special{pa 3668 2638}%
\special{pa 3668 2670}%
\special{pa 3668 2702}%
\special{pa 3668 2726}%
\special{sp}%
%
\special{pn 13}%
\special{pa 4790 842}%
\special{pa 4798 872}%
\special{pa 4806 904}%
\special{pa 4814 936}%
\special{pa 4822 966}%
\special{pa 4830 998}%
\special{pa 4836 1030}%
\special{pa 4844 1060}%
\special{pa 4850 1092}%
\special{pa 4856 1124}%
\special{pa 4862 1154}%
\special{pa 4868 1186}%
\special{pa 4874 1218}%
\special{pa 4878 1248}%
\special{pa 4882 1280}%
\special{pa 4886 1310}%
\special{pa 4888 1342}%
\special{pa 4890 1374}%
\special{pa 4890 1406}%
\special{pa 4888 1438}%
\special{pa 4886 1470}%
\special{pa 4880 1504}%
\special{pa 4874 1538}%
\special{pa 4864 1572}%
\special{pa 4850 1604}%
\special{pa 4834 1634}%
\special{pa 4816 1662}%
\special{pa 4794 1684}%
\special{pa 4766 1700}%
\special{pa 4736 1712}%
\special{pa 4704 1714}%
\special{pa 4672 1704}%
\special{pa 4648 1684}%
\special{pa 4630 1656}%
\special{pa 4620 1622}%
\special{pa 4620 1586}%
\special{pa 4628 1550}%
\special{pa 4646 1520}%
\special{pa 4670 1500}%
\special{pa 4698 1490}%
\special{pa 4730 1496}%
\special{pa 4760 1510}%
\special{pa 4788 1532}%
\special{pa 4810 1560}%
\special{pa 4824 1590}%
\special{pa 4834 1620}%
\special{pa 4838 1652}%
\special{pa 4840 1684}%
\special{pa 4840 1718}%
\special{pa 4840 1750}%
\special{pa 4842 1782}%
\special{pa 4842 1814}%
\special{pa 4842 1848}%
\special{pa 4842 1880}%
\special{pa 4842 1912}%
\special{pa 4842 1944}%
\special{pa 4842 1976}%
\special{pa 4842 2008}%
\special{pa 4842 2040}%
\special{pa 4842 2072}%
\special{pa 4842 2104}%
\special{pa 4842 2136}%
\special{pa 4842 2168}%
\special{pa 4842 2200}%
\special{pa 4842 2232}%
\special{pa 4842 2264}%
\special{pa 4842 2296}%
\special{pa 4842 2326}%
\special{pa 4842 2358}%
\special{pa 4840 2390}%
\special{pa 4840 2422}%
\special{pa 4840 2454}%
\special{pa 4840 2486}%
\special{pa 4840 2518}%
\special{pa 4840 2550}%
\special{pa 4840 2582}%
\special{pa 4840 2614}%
\special{pa 4840 2646}%
\special{pa 4840 2678}%
\special{pa 4840 2704}%
\special{sp}%
%
\special{pn 13}%
\special{sh 1}%
\special{ar 3884 1630 10 10 0  6.28318530717959E+0000}%
\special{sh 1}%
\special{ar 4150 1624 10 10 0  6.28318530717959E+0000}%
\special{sh 1}%
\special{ar 4406 1624 10 10 0  6.28318530717959E+0000}%
\special{sh 1}%
\special{ar 4406 1624 10 10 0  6.28318530717959E+0000}%
\put(40.6800,-10.0200){\makebox(0,0)[lb]{$s_1$}}%
\put(40.5200,-13.8500){\makebox(0,0)[lb]{$s_2$}}%
\put(40.7000,-20.2000){\makebox(0,0)[lb]{$s_3$}}%
\put(30.7400,-26.5700){\makebox(0,0)[lb]{$O$}}%
\put(12.2100,-23.2000){\makebox(0,0)[lb]{$\Theta_{00}$}}%
\put(12.2900,-19.9000){\makebox(0,0)[lb]{$\Theta_{10}$}}%
\put(12.6100,-16.0700){\makebox(0,0)[lb]{$\Theta_{01}$}}%
\put(17.6600,-9.1000){\makebox(0,0)[lb]{$\Theta_{11}$}}%
\put(11.5000,-10.9000){\makebox(0,0)[lb]{$\Theta_{4}$}}%
%
\special{pn 13}%
\special{pa 3020 780}%
\special{pa 3002 810}%
\special{pa 2982 840}%
\special{pa 2962 868}%
\special{pa 2944 898}%
\special{pa 2926 926}%
\special{pa 2908 956}%
\special{pa 2890 986}%
\special{pa 2872 1014}%
\special{pa 2856 1044}%
\special{pa 2842 1074}%
\special{pa 2828 1104}%
\special{pa 2814 1132}%
\special{pa 2802 1162}%
\special{pa 2792 1192}%
\special{pa 2782 1220}%
\special{pa 2774 1250}%
\special{pa 2768 1280}%
\special{pa 2764 1310}%
\special{pa 2762 1340}%
\special{pa 2760 1370}%
\special{pa 2760 1398}%
\special{pa 2764 1428}%
\special{pa 2768 1458}%
\special{pa 2774 1488}%
\special{pa 2782 1518}%
\special{pa 2792 1548}%
\special{pa 2802 1578}%
\special{pa 2812 1608}%
\special{pa 2824 1638}%
\special{pa 2836 1670}%
\special{pa 2850 1700}%
\special{pa 2862 1730}%
\special{pa 2876 1760}%
\special{pa 2888 1792}%
\special{pa 2900 1822}%
\special{pa 2914 1852}%
\special{pa 2924 1884}%
\special{pa 2936 1914}%
\special{pa 2946 1946}%
\special{pa 2954 1978}%
\special{pa 2962 2008}%
\special{pa 2968 2040}%
\special{pa 2974 2072}%
\special{pa 2980 2102}%
\special{pa 2984 2134}%
\special{pa 2988 2166}%
\special{pa 2990 2198}%
\special{pa 2992 2230}%
\special{pa 2994 2262}%
\special{pa 2994 2294}%
\special{pa 2994 2326}%
\special{pa 2994 2358}%
\special{pa 2992 2390}%
\special{pa 2992 2422}%
\special{pa 2990 2454}%
\special{pa 2988 2488}%
\special{pa 2986 2520}%
\special{pa 2984 2552}%
\special{pa 2982 2584}%
\special{pa 2978 2616}%
\special{pa 2976 2648}%
\special{pa 2972 2682}%
\special{pa 2970 2700}%
\special{sp}%
%
\special{pn 13}%
\special{pa 2850 770}%
\special{pa 2868 800}%
\special{pa 2888 830}%
\special{pa 2906 858}%
\special{pa 2922 888}%
\special{pa 2940 918}%
\special{pa 2958 946}%
\special{pa 2974 976}%
\special{pa 2990 1006}%
\special{pa 3006 1034}%
\special{pa 3020 1064}%
\special{pa 3034 1094}%
\special{pa 3048 1122}%
\special{pa 3060 1152}%
\special{pa 3070 1182}%
\special{pa 3080 1212}%
\special{pa 3088 1240}%
\special{pa 3096 1270}%
\special{pa 3102 1300}%
\special{pa 3106 1330}%
\special{pa 3110 1360}%
\special{pa 3110 1388}%
\special{pa 3110 1418}%
\special{pa 3108 1448}%
\special{pa 3106 1478}%
\special{pa 3100 1508}%
\special{pa 3094 1538}%
\special{pa 3086 1566}%
\special{pa 3078 1596}%
\special{pa 3068 1626}%
\special{pa 3058 1656}%
\special{pa 3046 1686}%
\special{pa 3032 1716}%
\special{pa 3020 1746}%
\special{pa 3004 1776}%
\special{pa 2990 1806}%
\special{pa 2974 1836}%
\special{pa 2958 1866}%
\special{pa 2942 1896}%
\special{pa 2924 1926}%
\special{pa 2906 1956}%
\special{pa 2890 1986}%
\special{pa 2872 2016}%
\special{pa 2854 2046}%
\special{pa 2850 2050}%
\special{sp}%
%
\special{pn 13}%
\special{pa 1730 1200}%
\special{pa 1754 1222}%
\special{pa 1780 1244}%
\special{pa 1804 1264}%
\special{pa 1828 1284}%
\special{pa 1854 1304}%
\special{pa 1880 1320}%
\special{pa 1908 1336}%
\special{pa 1936 1350}%
\special{pa 1964 1360}%
\special{pa 1994 1370}%
\special{pa 2024 1380}%
\special{pa 2054 1386}%
\special{pa 2084 1392}%
\special{pa 2116 1396}%
\special{pa 2148 1400}%
\special{pa 2180 1402}%
\special{pa 2214 1404}%
\special{pa 2246 1406}%
\special{pa 2280 1406}%
\special{pa 2314 1406}%
\special{pa 2348 1404}%
\special{pa 2382 1404}%
\special{pa 2416 1402}%
\special{pa 2450 1400}%
\special{pa 2460 1400}%
\special{sp}%
%
\special{pn 13}%
\special{pa 1690 970}%
\special{pa 1714 992}%
\special{pa 1740 1014}%
\special{pa 1764 1034}%
\special{pa 1788 1054}%
\special{pa 1814 1074}%
\special{pa 1840 1090}%
\special{pa 1868 1106}%
\special{pa 1896 1120}%
\special{pa 1924 1130}%
\special{pa 1954 1140}%
\special{pa 1984 1150}%
\special{pa 2014 1156}%
\special{pa 2044 1162}%
\special{pa 2076 1166}%
\special{pa 2108 1170}%
\special{pa 2140 1172}%
\special{pa 2174 1174}%
\special{pa 2206 1176}%
\special{pa 2240 1176}%
\special{pa 2274 1176}%
\special{pa 2308 1174}%
\special{pa 2342 1174}%
\special{pa 2376 1172}%
\special{pa 2410 1170}%
\special{pa 2420 1170}%
\special{sp}%
%
\special{pn 13}%
\special{pa 2410 1180}%
\special{pa 2760 1170}%
\special{fp}%
\special{pa 2870 1170}%
\special{pa 4980 1180}%
\special{fp}%
%
\special{pn 13}%
\special{pa 2430 1400}%
\special{pa 2700 1400}%
\special{fp}%
\special{pa 2830 1400}%
\special{pa 5040 1410}%
\special{fp}%
%
\special{pn 13}%
\special{pa 2940 1780}%
\special{pa 5070 1780}%
\special{fp}%
%
\special{pn 13}%
\special{pa 1740 2240}%
\special{pa 1760 2214}%
\special{pa 1780 2190}%
\special{pa 1800 2164}%
\special{pa 1822 2138}%
\special{pa 1842 2114}%
\special{pa 1864 2088}%
\special{pa 1884 2064}%
\special{pa 1906 2040}%
\special{pa 1928 2018}%
\special{pa 1952 1996}%
\special{pa 1974 1974}%
\special{pa 1998 1954}%
\special{pa 2024 1934}%
\special{pa 2050 1914}%
\special{pa 2076 1898}%
\special{pa 2102 1880}%
\special{pa 2132 1866}%
\special{pa 2160 1852}%
\special{pa 2190 1838}%
\special{pa 2220 1826}%
\special{pa 2250 1816}%
\special{pa 2280 1806}%
\special{pa 2312 1798}%
\special{pa 2344 1790}%
\special{pa 2374 1784}%
\special{pa 2406 1778}%
\special{pa 2438 1774}%
\special{pa 2470 1772}%
\special{pa 2502 1768}%
\special{pa 2534 1766}%
\special{pa 2566 1766}%
\special{pa 2598 1766}%
\special{pa 2630 1766}%
\special{pa 2662 1766}%
\special{pa 2694 1766}%
\special{pa 2726 1768}%
\special{pa 2758 1768}%
\special{pa 2790 1770}%
\special{pa 2800 1770}%
\special{sp}%
\end{picture}%

\end{center}

By the explicit formula for $\langle \, , \, \rangle$, we have
\[
\langle s_{L_i}^+, s_{L_i}^+ \rangle = \frac 12 (i = 1, 2, 3), \quad \langle s_{L_i}^+, s_{L_j}^+ \rangle = 0\,\,  (i \neq j).
\]
By \cite{oguiso-shioda}, we have $\MW(S_{(\mcQ, z_o)}) \cong (A_1^*)^{\oplus 3}$. Hence we may assume that
\[
\MW(S_{\mcQ, z_o}) \cong \ZZ s_{L_1}^+ \oplus \ZZ s_{L_2}^+ \oplus \ZZ s_{L_3}^+.
\]
By the lattice structure of $\MW((S_{\mcQ, z_o})$, all elements $s \in \MW(S_{\mcQ, z_o}) $ with $\langle s, s \rangle =2$
given by $[2]s_{L_i}^{\pm}$ ($i = 1, 2, 3$).  By \cite[Theorem 9.1]{kodaira}, $[2]s_{L_i}^{\pm}$ 
($i = 1, 2, 3$) meet the identity component at each singular fiber. Hence, $[2]s_{L_i}^{\pm}O = 0$ $(i =1, 2, 3)$ by
the explicit formula for $\langle \, , \, \rangle$. By our construction of $S_{(\mcQ, z_o)}$, $\Delta_i := q\circ f([2]s_{L_i}^{\pm})
\sim \Delta_0 + 2{\mathfrak f}$ ($i = 1, 2, 3$). Hence $C_i:= \overline{q}\circ f([2]s_{L_i}^{\pm})$ ($i = 1, 2, 3$) are
all conic and $z_o \in C_i, P \not\in C_i$.  Moreover as $[2]s_{L_i}^+ \neq [2]s_{L_I}^-$ ($i = 1, 2, 3$), our assertion
for the intersection multiplicities  follows.

(ii) By Theorem~\ref{thm:main1} and Lemma~\ref{lem:equiv}, our statement follows.
\proofend

\begin{cor}\label{cor:eg-3}{If $L_i$ and $L_j$ $(i \neq j)$ meet $C_i$ transversely, then $({\mathcal Q} + L_i + C_i, 
{\mathcal Q} + L_j + C_i)$ is a Zariski pair.
}
\end{cor}

\proof Since the combinatorics of  ${\mathcal Q} + L_i + C_i$ and 
${\mathcal Q} + L_j + C_i$ are the same, our assertion follows from 
Proposition~\ref{prop:eg-3}.  \proofend

\begin{exmple}\label{eg:eg-3}{\rm 

We keep the same coordinates as in 
Examples~\ref{eg:line-conic1}
and~\ref{eg:line-conic2}. Consider $Q, L_1$ and $L_2$ as follows:
\[
\mcQ: f_{\mcQ} (t, x) :=  x^3 + \frac{343}{64}\left (\frac {121}{49}t^2 + \frac {768}{2401}t\right )x^2
+ \frac{343}{64}\left (\frac{384}{2401}t^2 + \frac{92}{49}t^3\right ) u + 
\frac {35}{16}t^4 + \frac 17 t^3 = 0
\]
\[
L_1: x + t = 0, \quad \quad  L_2: x - \frac {\zeta_3 - 2}7 t = 0, \,\,
\zeta_3 = \exp(2\pi i/3)
\]
$\mcQ$ is irreducible and has a ${\mathbb D}_4$ singularity at $(0, 0)$. 
Both $L_1$ and $L_2$ meet $\mcQ$ at $(0, 0)$ with multiplicity $4$.
Choose $[0, 1, 0]$ as the distinguished point $z_o$. Let $S_{(\mcQ, z_o)}$ be the rational elliptic surface obtained as before. Then its generic fiber is 
an elliptic curve over $\CC(t)$ given by the Weierstrass equation $y^2 = f_{\mcQ}(t, x)$.
Under these circumstances,  we have
\[
s_{L_1}^{\pm} = \left (-t, \pm \frac{\sqrt{343}}8 t^2\right ), \quad\quad 
s_{L_2}^{\pm}  =\left ( \frac {\zeta_3 -2}7 t, 
\pm \frac{\sqrt{71 + 39\sqrt{-3}}}{8\sqrt{14}}t^2 \right )
\]
Then we have
\[
[2]s_{L_1}^+ = \left (\frac{144}{16807} - \frac{127}{343} t - \frac{19}{28}t^2, 
- \frac {\sqrt{7}(55296 + 1947456t + 1450204t^2 + 167649825t^3)}{184473632}
 \right ).
\]
Now put
\[
C : x - \frac{144}{16807} + \frac{127}{343} t + \frac{19}{28}t^2.
\]
Since one can see  that  both of $L_1$ and $L_2$ meet $C$ with two distinct points,
 $\mcQ + C + L_1$ and $\mcQ + C + L_2$ have the same
combinatorics. By Corollary~\ref{cor:eg-3}, $(\mcQ + C + L_1, \mcQ + C + L_2)$
is a Zariski pair.
}
\end{exmple}


\noindent Hiro-o TOKUNAGA\\
Department of Mathematics and Information Sciences\\
Graduate School of Science and Engineering,\\
Tokyo Metropolitan University\\
1-1 Minami-Ohsawa, Hachiohji 192-0397 JAPAN \\
{\tt tokunaga@tmu.ac.jp}
      
  \end{document}